\documentclass[a4paper,12pt]{article}
\usepackage[latin1]{inputenc}
\usepackage[T1]{fontenc}
\usepackage[francais]{babel}
\usepackage{amssymb}
\usepackage{amsfonts, amsmath}
\newtheorem{defn}{D\'efinition}[section]
\newtheorem{thm}[defn]{Th\'eor\`eme}
\newtheorem{prop}[defn]{Proposition}

\newtheorem{rem}{Remarque}[section]

\newtheorem{Exmp}{Exemple}[section]

\numberwithin{equation}{section}
\title{Surfaces de Riemann compactes, courbes alg\'{e}briques complexes et leurs Jacobiennes}

\author{\textbf{A. Lesfari}
\\\emph{B.P. 271, Poste principale, 24000 El-Jadida, Maroc.} \\\emph{E.
mail} : lesfariahmed@yahoo.fr, lesfari@gmail.com}

\date{}

\begin{document}
\maketitle

\begin{abstract}
Topologically, a compact Riemann surface $X$ of genus
$g$ is a $g$-holed torus (a sphere with $g$ handles). This paper
is an introduction to the theory of compact Riemann surfaces and
algebraic curves. It presents the basic ideas and properties as an
expository essay, explores some of their numerous consequences and
gives a concise account of the elementary aspects of different
viewpoints in curve theory. We discuss and prove most intuitively
some geometric-topological aspects of the algebraic functions and
the associated Riemann surfaces. Abelian and normalized
differentials, Riemann's bilinear relations and the period matrix
for $X$ are defined and some consequences drawn. The space of
holomorphic $1$-forms on $X$ has dimension $g$ as a complex vector
space. Fundamental results on divisors on compact Riemann surfaces
are stated and proved. The Riemann-Roch theorem is of utmost
importance in the algebraic geometric theory of compact Riemann
surfaces. It tells us how many linearly independent meromorphic
functions there are having certain restrictions on their poles. We
present a simple direct proof of this theorem and explore some of
its numerous consequences. We also give an analytic proof of the
Riemann-Hurwitz formula. As an application, we compute the genus
of some interesting algebraic curves. Abel's theorem classifies
divisors by their images in the jacobian. The Jacobi inversion
problem askes whether we can find a divisor that is the preimage
for an arbitrary point in the jacobian. In the first appendix, we
introduced intuitively and explicitly elliptic and hyperelliptic
Riemann surfaces. In the second appendix, we study some results of
resultant and discriminant as needed in the paper.
\end{abstract}

\vskip0.6cm \tableofcontents

\section{Introduction}

Les surfaces de Riemann interviennent souvent lors de la
r\'{e}solution de probl\`{e}mes aussi bien th\'{e}oriques que
pratiques et sont la source de plusieurs domaines de la recherche
contemporaine.

Dans ce travail, on \'{e}tudie les surfaces de Riemann compactes
$X$ ou courbes alg\'{e}briques complexes. Ce sont des
vari\'{e}t\'{e}s analytiques de dimension $1$ complexe ($2$
r\'{e}elle) munies d'atlas dont les changements de cartes sont
holomorphes.  On montre que la courbe $X$ est hom\'{e}omorphe
\`{a} un tore \`{a} $g$ trous (ou sph\`{e}re \`{a} $g$ anses) pour
un certain entier $g \geq 0$. Le nombre $g$ est le genre de $X$.
Celui-ci est la dimension de l'espace vectoriel complexe
$H^{1}\left( X,\mathcal{O}_{X}\right) $ ($1^{er}$groupe de
cohomologie \`{a} coefficients dans le faisceau $\mathcal{O}_{X}$
des fonctions holomorphes sur $X$). On montre aussi que c'est le
nombre des int\'{e}grales ab\'{e}liennes de
$1^{\grave{e}re}$esp\`{e}ce attach\'{e}es \`{a} la courbe $X$,
lin\'{e}airement ind\'{e}pendants. Un cas particulier important
est repr\'{e}sent\'{e} par les courbes hyperelliptiques de genre
$g$ ainsi que les courbes elliptiques ($g=1$). On \'{e}tudie
ensuite les formes diff\'{e}rentielles, les relations
bilin\'{e}aires de Riemann et la matrice des p\'{e}riodes.
Apr\`{e}s avoir rappel\'{e} les d\'{e}finitions et
propri\'{e}t\'{e}s des diviseurs et des fibr\'{e}s en droites
n\'{e}cessaires \`{a} la compr\'{e}hension des r\'{e}sultats
principaux de ce travail, on aborde le th\'{e}or\`{e}me de
Riemann-Roch. Ce dernier est un r\'{e}sultat central de la
th\'{e}orie des surfaces de Riemann compactes. Il permet, entre
autres, de d\'{e}finir le genre d'une surface de Riemann qui est
un invariant fondamental. Il s'agit d'un th\'{e}or\`{e}me
d'existence efficace qui permet, entre autres, de d\'{e}terminer
le nombre de fonctions m\'{e}romorphes lin\'{e}airement
ind\'{e}pendantes ayant certaines restrictions sur leurs
p\^{o}les. A cause de l'importance de ce th\'{e}or\`{e}me, nous
donnons une preuve d\'{e}taill\'{e}e constructive bien qu'un peu
technique. Nous mentionnons quelques cons\'{e}quences de ce
th\'{e}or\`{e}me et nous donnons \'{e}galement une preuve
analytique de l'importante formule de Riemann-Hurwitz. Elle
exprime le genre d'une surface de Riemann \`{a} l'aide du nombre
de ses points de ramifications et du nombre de ses feuillets. Nous
montrons que cette formule fournit un moyen efficace pour
d\'{e}terminer le genre d'une surface de Riemann donn\'{e}e. En
outre plusieurs exemples int\'{e}ressants seront \'{e}tudi\'{e}s.
Deux autres th\'{e}or\`{e}mes, celui d'Abel et celui de Jacobi, de
nature transcendante et consid\'{e}r\'{e}s comme importants de la
th\'{e}orie des surfaces de Riemann compactes, sont
\'{e}tudi\'{e}s en d\'{e}tail. Le th\'{e}or\`{e}me d'Abel
classifie les diviseurs par leurs images dans la vari\'{e}t\'{e}
jacobienne (tore complexe alg\'{e}brique) tandis que le
probl\`{e}me d'inversion de Jacobi concerne l'existence d'un
diviseur qui soit l'image inverse d'un point arbitraire sur la
vari\'{e}t\'{e} jacobienne.

Deux appendices enfin, expliquent certaines notions utilis\'{e}es
dans les sections pr\'{e}c\'{e}dentes et qui auraient autrement
alourdi le texte. Dans le premier appendice, on introduit de
mani\`{e}re intuitive et on construit explicitement les surfaces
de Riemann dans les cas elliptique et hyperelliptique. Dans le
second appendice, on \'{e}tudie quelques r\'{e}sultats
d'alg\`{e}bre concernant les r\'{e}sultants et discriminants que
l'on utilise dans la preuve de la connexit\'{e} de la surface de
Riemann construite.

On ne trouvera pas dans ces notes des figures qui aident à motiver
et comprendre de mani\`{e}re intuitive les diff\'{e}rentes
configurations g\'{e}om\'{e}triques. Ces lacunes peuvent \^{e}tre
ais\'{e}ment combl\'{e}es en pratique par l'expos\'{e} oral comme
nous l'avons fait pour nos \'{e}tudiants.

\section{Etude g\'{e}om\'{e}trique et topologique}

Dans cette partie nous allons \'{e}tudier les courbes
alg\'{e}briques complexes $X$ ou surfaces de Riemann compactes. Ce
sont des vari\'{e}t\'{e}s analytiques de dimension $1$ complexe
($2$ r\'{e}elle) munies d'atlas dont les changements de cartes
sont holomorphes. On les d\'{e}finit par
$$
X=\{(w,z)\in \mathbb{C}^2 : F(w,z)=0\},
$$
o\`{u}
\begin{equation}\label{eqn:euler}
F(w,z) \equiv p_{0}\left( z\right) w^{n}+p_{1}\left( z\right)
w^{n-1}+\cdots +p_{n}\left( z\right),
\end{equation}
est un polyn\^{o}me \`{a} deux variables complexes $w$ et $z$, de
degr\'{e} $n$ en $w$ et irr\'{e}ductible (i.e., sans facteurs
multiples ou encore ne soit pas le produit de deux autres
polyn\^{o}mes en $w$ et $z$). Ici $p_{0}\left( z\right) \neq 0$,
$p_{1}\left( z\right) ,\ldots ,$ $p_{n}\left( z\right) $ sont des
polyn\^{o}mes en $z$. Nous montrerons que la courbe $X$ est
hom\'{e}omorphe \`{a} un tore \`{a} $g$ trous (ou sph\`{e}re \`{a}
$g$ anses) pour un certain entier $g \geq 0$. Le nombre $g$
s'appelle genre de la courbe $X$. Nous verrons aussi qu'il est
\'{e}quivalent de dire que c'est la dimension de l'espace
vectoriel complexe $H^{1}\left( X,\mathcal{O}_{X}\right)$
($1^{er}$groupe de cohomologie \`{a} coefficients dans le faisceau
$\mathcal{O}_{X}$ des fonctions holomorphes sur $X$). On montrera
aussi que c'est le nombre des int\'{e}grales ab\'{e}liennes de
$1^{\grave{e}re}$esp\`{e}ce (voir plus loin pour les
d\'{e}finitions) attach\'{e}es \`{a} la courbe $X$,
lin\'{e}airement ind\'{e}pendants. Un cas particulier important
est repr\'{e}sent\'{e} par les courbes hyperelliptiques de genre
$g$ d'\'{e}quation
$$F\left( w,z\right) =w^{2}-p_{n}\left( z\right) =0,$$
o\`{u} $p_{n}\left( z\right) $ est un polyn\^{o}me sans racines
multiples, de degr\'{e} $n=2g+1$ ou $2g+2.$ Lorsque $g=1,$ on dit
courbes elliptiques.

Tout au long de cette partie, $X$ d\'{e}signe une courbe
alg\'{e}brique complexe (surface de Riemann compacte) de genre $g$
non-singuli\`{e}re. Celà signifie que les fonctions
$\frac{\partial F}{\partial w}$ ou $\frac{\partial F}{\partial z}$
ne s'annulent identiquement sur aucune composante de $X$ ou encore
que
$$\mbox{grad } F\equiv
\left( \frac{\partial F}{\partial w},\frac{\partial F}{\partial
z}\right) \neq 0.
$$

Consid\'{e}rons donc l'\'{e}quation (2.1). A chaque valeur de $z$
correspond $n$ valeurs de $w$. Notre probl\`{e}me consiste \`{a}
trouver un domaine pour lequel
$$\mathbb{C}\longrightarrow \mathbb{C},\quad z\longmapsto w
:F(w,z)=0,$$ soit une fonction uniforme. Autrement dit, on cherche
\`{a} construire la surface de Riemann associ\'{e}e \`{a}
l'\'{e}quation $F(w,z)=0$. D\'{e}signons par $z_1,...,z_m$ les
z\'{e}ros de $p_0(z)$ et les z\'{e}ros communs de $F(w,z)=0$ et
$\frac{\partial F}{\partial w}(w,z)=0$. Ce sont les valeurs pour
lesquelles $F(w,z)$ a un z\'{e}ros double en $w$. Soit
$\overline{\mathbb{C}}=\mathbb{C}\cup \{\infty\}$ le plan complexe
compactifi\'{e} ( ou sph\`{e}re de Riemann puisqu'ils sont
hom\'{e}omorphes). En r\'{e}solvant l'\'{e}quation $F(w,z)=0$ pour
$z\in \overline{\mathbb{C}}\backslash\{z_1,...,z_m\}$, on obtient
$n$ solutions $w_k(z)$, $k=1, 2,...,n$.

\begin{thm}
Les solutions $w_k(z)$ sont localement analytiques.
\end{thm}
\emph{D\'{e}monstration}: Posons $z=x+iy$ et $w=u+iv$ o\`{u}
$u=u(x,y)$, $v=v(x,y)$. Donc
$$F(w,z)=G(u,v,x,y)+iH(u,v,x,y),$$
o\`{u} $G$ et $H$ sont des polyn\^{o}mes. Par cons\'{e}quent
$$F(w,z)=0\quad \Longleftrightarrow\quad
\left\{\begin{array}{rl}
G(u,v,x,y)=0,&\\
H(u,v,x,y)=0.&
\end{array}\right.
$$
Fixons un point $z_0=x_0+iy_0\in
\overline{\mathbb{C}}\backslash\{z_1,...,z_m\}$. L'\'{e}quation
$F(w,z)=0$ a exactement $n$ racines distinctes : $w=w_{01},
w_{02},...,w_{0n}$. Soit $w_0=u_0+iv_0$ l'une de ces racines. Pour
pouvoir appliquer le th\'{e}or\`{e}me des fonctions implicites, il
suffit de v\'{e}rifier que
$$
\det\left(\begin{array}{cc}
\frac{\partial G}{\partial u}&\frac{\partial G}{\partial v}\\
\frac{\partial H}{\partial u}&\frac{\partial H}{\partial v}
\end{array}\right),
$$
est non nul en $(u_0,v_0)$. En effet, pour $z$ fix\'{e},
$F(w,z_0)$ est un polyn\^{o}me en $w$ et donc $F(w,z_0)$ est
analytique. D'apr\`{e}s les \'{e}quations de Cauchy-Riemann, on a
\begin{eqnarray}
\frac{\partial G}{\partial u}&=&\frac{\partial H}{\partial v},\nonumber\\
\frac{\partial G}{\partial v}&=&-\frac{\partial H}{\partial
u},\nonumber
\end{eqnarray}
et d\`{e}s lors
\begin{eqnarray}
\det\left(\begin{array}{cc}
\frac{\partial G}{\partial u}&\frac{\partial G}{\partial v}\\
\frac{\partial H}{\partial u}&\frac{\partial H}{\partial v}
\end{array}\right)&=&\frac{\partial G}{\partial u}\frac{\partial H}{\partial
v}-\frac{\partial G}{\partial v}\frac{\partial H}{\partial
u},\nonumber\\
&=&\left(\frac{\partial G}{\partial
u}\right)^2+\left(\frac{\partial H}{\partial
u}\right)^2,\nonumber\\
&=&\left| \frac{\partial G}{\partial u}+i\frac{\partial
H}{\partial
u}\right| ^2,\nonumber\\
&=&\left| \frac{\partial F}{\partial u}\right| ^2.\nonumber
\end{eqnarray}
De m\^{e}me, on a
$$
\det\left(\begin{array}{cc}
\frac{\partial G}{\partial u}&\frac{\partial G}{\partial v}\\
\frac{\partial H}{\partial u}&\frac{\partial H}{\partial v}
\end{array}\right)=\left| \frac{\partial F}{\partial v}\right| ^2.
$$
Donc pour que
$$\det\left(\begin{array}{cc}
\frac{\partial G}{\partial u}&\frac{\partial G}{\partial v}\\
\frac{\partial H}{\partial u}&\frac{\partial H}{\partial v}
\end{array}\right)\neq 0,$$
il suffit que
$$\left| \frac{\partial F}{\partial u}\right| ^2\neq0\quad
\Longleftrightarrow \quad\left| \frac{\partial F}{\partial
v}\right| ^2\neq 0,$$ ou encore
$$\left| \frac{\partial F}{\partial w}\right| ^2\neq 0.$$
Or par hypoth\`{e}se, l'\'{e}quation $F(w,z)=0$ n'a pas de racines
double en $w$ pour $z$ fix\'{e}, i.e.,
$$\frac{\partial F}{\partial w}(w_0,z_0)\neq 0.$$
Par le th\'{e}or\`{e}me des fontions implicites, on peut
r\'{e}soudre $w$ en fonction de $z$ et exprimer que $w$ est
diff\'{e}rentiable dans un voisinage de $z$. Pour montrer que
$w=w(z)$ est analytique, on va v\'{e}rifier que les \'{e}quations
de Cauchy-Riemann
\begin{eqnarray}
\frac{\partial u}{\partial x}&=&\frac{\partial v}{\partial y},\nonumber\\
\frac{\partial v}{\partial x}&=&-\frac{\partial u}{\partial
y},\nonumber
\end{eqnarray}
sont satisfaites. En effet, comme $F(w,z)=0$, alors
\begin{eqnarray}
\frac{\partial F}{\partial x}&=&\frac{\partial F}{\partial
w}\frac{\partial w}{\partial x}
+\frac{\partial F}{\partial z}=0,\nonumber\\
\frac{\partial F}{\partial y}&=&-\frac{\partial F}{\partial
w}\frac{\partial w}{\partial y}+i\frac{\partial F}{\partial
z}=0.\nonumber
\end{eqnarray}
On multiplie la deuxi\`{e}me \'{e}quation par $i$ et on fait la
somme avec la premi\`{e}re. On obtient
$$\frac{\partial F}{\partial w}\left(\frac{\partial w}{\partial x}+
i\frac{\partial w}{\partial y}\right)=0.$$ Or $\frac{\partial
F}{\partial w}(w_0,z_0)\neq 0$, donc la seule possibilit\'{e} qui
reste est
$$\frac{\partial w}{\partial x}=-i\frac{\partial w}{\partial y},$$
i.e.,
$$\frac{\partial u}{\partial x}+
i\frac{\partial v}{\partial x}=-i\frac{\partial u}{\partial y}+
\frac{\partial v}{\partial y},$$ ce qui ach\`{e}ve la
d\'{e}monstration. $\square$

Nous allons montrer que l'on peut prolonger analytiquement
$w_k=w_k(z)$ sur tout
$\overline{\mathbb{C}}\backslash\{z_1,...,z_m\}$ et chaque
fonction ainsi obtenue satisfait \`{a} l'\'{e}quation $F(w,z)=0$.
Mais auparavant nous aurons besoin de quelques pr\'{e}liminaires.

Soit $D(a,r_a)$ un disque de centre $a$ et de rayon $r_a$ et soit
$f$ une fonction analytique sur $D(a,r_a)$. Cette fonction admet
un d\'{e}veloppement en s\'{e}rie enti\`{e}re convergente de la
forme
$$f(z)=f(a)+a_1(z-a)+a_2(z-a)^2+\cdots$$

\begin{prop}
Soit $b\in D(a,r_a)$. La fonction $f$ admet un d\'{e}veloppement
en s\'{e}rie enti\`{e}re convergente autour de $b$ de la forme
$$f(z)=f(b)+b_1(z-b)+b_2(z-b)^2+\cdots$$
\end{prop}
\emph{D\'{e}monstration}: Posons
$$z-a=(z-b)+(b-a),$$
d'o\`{u}
$$f(z)=f(a)+a_1(b-a)+a_2(b-a)^2+\cdots+a_1(z-b)+2a_2(b-a)(z-b)+\cdots+a_2(z-b)^2+\cdots$$
On en d\'{e}duit que
\begin{eqnarray}
f(b)&=&f(a)+a_1(b-a)+a_2(b-a)^2+\cdots\nonumber\\
b_1&=&a_1+2a_2(b-a)+\cdots\nonumber\\
b_2&=&a_2+\cdots\nonumber
\end{eqnarray}
et
\begin{equation}\label{eqn:euler}
f(z)=f(b)+b_1(z-a)+b_2(z-a)^2+\cdots
\end{equation}
Cette s\'{e}rie converge en tout point de $D(a,r_a)$. Cherchons
maintenant le disque $D(b,r_b)$ de centre $b$ et de rayon $r_b$,
dans lequel la s\'{e}rie (2.2) converge. On sait que
$$r_b\geq r_a-\mid b-a\mid>0,$$
et il se peut aussi qu'on ait
$$r_b> r_a-\mid b-a\mid>0.$$
Si tel est le cas, la s\'{e}rie (2.2) convergera aussi \`{a}
l'ext\'{e}rieur du disque $D(a,r_a)$, \`{a} savoir dans le domaine
du disque $D(b,r_b)$ ext\'{e}rieur au disque $D(a,r_a)$ et la
proposition est d\'{e}montr\'{e}e. $\square$

Soit $f$ une fonction analytique sur $D(a,r_a)$ et soit $b$ un
point en dehors de $D(a,r_a)$. On veut construire un prolongement
analytique de $f$ au point $b$. Du point $a$ au point $b$, traçons
un chemin $\mathcal{C}$. Soit $c_1\in \mathcal{C}\cap D(a,r_a)$.
On sait que la fonction $f$ peut-\^{e}tre d\'{e}velopp\'{e}e en
s\'{e}rie enti\`{e}re de $z-c_1$. Soit $D_1$, sur le chemin entre
$c_1$ et $b$. En ce point, $f$ admet un prolongement analytique.
Soit $D_2$, le disque de centre $c_2$ dans lequel le
d\'{e}veloppement obtenu est convergent. De proche en proche, on
avance progressivement sur $\mathcal{C}$, vers le point $b$. Quand
$b$ sera dans un disque $D_n$ de centre $c_n$, on prendra $b$ pour
$c_{n+1}$ et ainsi on obtiendra le prolongement analytique
cherch\'{e}. Signalons que cette construction ne d\'{e}montre pas
l'existence du prolongement analytique.

\begin{thm} (de monodromie):
Soit $f$ une fonction analytique dans un voisinage de $a$ et soit
$D$ un domaine simplement connexe. On suppose que pour tout $x\in
D$, il existe un chemin de $a$ vers $x$ tel que $f$ peut-\^{e}tre
prolong\'{e}e analytiquement en $x$. Alors ce prolongement ne
d\'{e}pend pas du chemin suivi.
\end{thm}
\emph{D\'{e}monstration}: Soient $\mathcal{C}$ et $\mathcal{L}$
deux chemins de $a$ vers $x$. On prolonge $f$ le long du chemin
$\mathcal{C}$. On d\'{e}signe par $\mathcal{C}_1,
\mathcal{C}_2,...,\mathcal{C}_n$ les points sur $\mathcal{C}$
interm\'{e}diaires de prolongement analytique. On a d\'{e}j\`{a}
montr\'{e} que : quel que soit $\mathcal{C}_i$, il existe au moins
$r_i>0$ tel que le prolongement analytique de $f$ en
$\mathcal{C}$, converge dans $D_i(\mathcal{C}_i,r_i)$. Notons que
$\bigcup_{i=1}^nD_i$ recouvre $\mathcal{C}$. Soit $\mathcal{C}_1$
un autre chemin de $a$ vers $x$ tel que : $\mathcal{C}_1\subset
\bigcup_{i=1}^nD_i$. Soient $c_{11},...,c_{1m}$, les points
interm\'{e}diaires sur $\mathcal{C}_1$ tels que : quel que soit
$j$, il existe au moins $i$, $c_{1j}\in D_i$ et le prolongement
sera valable sur un disque $D_{1j}(c_{1j},r_{1j})$. Donc, quel que
soit $j$, le prolongement en $c_{1j}$ n'est rien d'autre qu'un
r\'{e}arrangement des termes du prolongement en $c_i$. Ainsi,
lorsqu'on arrive en $x=c_{n+1}=c_{1,m+1}$, les prolongements
coincident sur $D_{n+1}\cap D_{1,m+1}$, et \`{a} fortiori en $x$.
Notons que $\bigcup_{j=1}^mD_{1j}$ recouvre $\mathcal{C}_1$. Comme
pr\'{e}c\'{e}demment, choisissons $\mathcal{C}_2$, un chemin de
$a$ vers $x$ tel que : $\mathcal{C}_2\subset
\bigcup_{j=1}^mD_{1j}$. Le prolongement en $x$ le long de ce
chemin coincidera avec celle en $x$ le long de $\mathcal{C}_1$.
Donc le prolongement le long de $\mathcal{C}_2$ coincide avec le
prolongement le long de $\mathcal{C}$, au point $x$. On peut
continuer la proc\'{e}d\'{e} jusqu'\`{a} ce que le chemin obtenu
soit $\mathcal{L}$. Montrons que l'on peut atteindre $\mathcal{L}$
par ce proc\'{e}d\'{e}. En effet, nous avons suppos\'{e} que le
domaine $D$ est simplement connexe car sinon le proc\'{e}d\'{e}
utilis\'{e} peut ne pas marcher. Si le trou du domaine se trouve
entre $\mathcal{L}$ et $\mathcal{C}$, on ne pourra pas passer
continument de l'un \`{a} l'autre car on ne peut pas traverser le
trou. Supposons maintenat qu'il soit impossible de se rapprocher
de $\mathcal{L}$. Supposons donc qu'il y ait un chemin limite
$\mathcal{K}$ qu'on ne puisse atteindre entre $\mathcal{C}$ et
$\mathcal{L}$. Dans ce cas, il suffit de prendre le prolongement
analytique le long de ce chemin $\mathcal{K}$. Celui-ci est
recouvert par une suite de disques. Comme on peut se rapprocher de
$\mathcal{K}$ aussi pr\`{e}s que l'on veut, il est possible de
trouver un chemin qui soit recouvert par ces disques et qui
conduise ainsi au m\^{e}me prolongement en $x$. On peut alors
continuer de l'autre c\^{o}t\'{e} \`{a} l'int\'{e}rieur des
disques, le proc\'{e}d\'{e} d\'{e}j\`{a} commenc\'{e}. Notons
enfin que si c'est un point $p$ qui nous barrait la route, alors
celui-ci jouerait le m\^{e}me r\^{o}le que le trou. Or par
hypoth\`{e}se, la fonction $f$ peut-\^{e}tre prolong\'{e}e
analytiquement \`{a} tous les points de $D$, donc \`{a}
celui-l\`{a} aussi. Soit donc $\mathcal{K}'$ un chemin de $a$ vers
$p$, le long duquel on peut prolonger $f$ analytiquement.
Prolongeons $\mathcal{K}'$ jusqu'\`{a} $x$, on obtient le chemin
$\mathcal{K}"$. Celui-ci joue alors le r\^{o}le du chemin
$\mathcal{K}$ du cas pr\'{e}c\'{e}dent. La solution est donc la
m\^{e}me et la d\'{e}monstration est compl\`{e}te. $\square$

Revenons maintenant \`{a} notre probl\`{e}me initial. Nous avons
montr\'{e} qu'en tout point $z\in
\overline{\mathbb{C}}\backslash\{z_1,...,z_m\}$, les solutions
$w_k(z)$, $k=1, 2,...,n$ de l'\'{e}quation $F(w,z)=0$ sont
localement analytiques. En outre, on a

\begin{thm}
On peut prolonger analytiquement les fonctions  $w_k=w_k(z)$ sur
tout $\overline{\mathbb{C}}\backslash\{z_1,...,z_m\}$ et chaque
fonction ainsi obtenue satisfait à l'\'{e}quation : $F(w,z)=0$.
\end{thm}
\emph{D\'{e}monstration}: Soit $z_0\in
\overline{\mathbb{C}}\backslash\{z_1,...,z_m\}$. D'apr\`{e}s les
r\'{e}sultats pr\'{e}c\'{e}dents, on peut prolonger analytiquement
$w_k$ le long de tous les chemins contenus dans un voisinage
suffisamment petit de $z_0$. Pour le reste, on va utiliser un
raisonnement par l'absurde. Supposons qu'il existe un point $a$ et
un chemin $\mathcal{C}$ de $z_0$ vers $a$, de sorte que l'on
puisse prolonger $w_k$ le long de $\mathcal{C}$ jusqu'\`{a} tous
les points $b$ avant $a$, mais pas jusqu'au $a$. Dans $D(a,r_a)$,
il existent $n$ s\'{e}ries enti\`{e}res $w_k(z)$ qui satisfont
\`{a} l'\'{e}quation $F(w_k(z),z)=0$ et qui en donnent toutes les
solutions. Choisissons $b$, de sorte que la partie du chemin
$\mathcal{C}$ comprise entre $a$ et $b$, soit enti\`{e}rement
incluse dans le disque $D(a,r_a)$. Consid\'{e}rons la s\'{e}rie
$w(z)$ r\'{e}sultant du prolongement analytique de $w_k(z)$ de
$z_0$ jusqu'au $b$, le long du chemin $\mathcal{C}$. D'apr\`{e}s
ce qui pr\'{e}c\`{e}de, on a $F(w(z),z)=0$ dans un disque
$D(b,r_b)$ autour du point $b$. Donc dans un disque $D(b,r)$ de
rayon $r=\min(r_a-\mid a-b\mid,r_b)$, $w(z)$ doit coincider avec
l'une des $w_k(z)$ puisque dans ce disque toutes les solutions de
l'\'{e}quation $F(w,z)=0$ sont donn\'{e}es par les fonctions
$w_k(z)$ et que $F(w(z),z)=0$. Soit $w_l$ qui coincide avec $w$
dans $D(b,r)$. Donc $w_k(z)$ peut-\^{e}tre prolong\'{e} de $b$
vers $a$ le long de $\mathcal{C}$. Ceci contredit l'hypoth\`{e}se
de d\'{e}part et d\'{e}montre le th\'{e}or\`{e}me. $\square$

On veut que le prolongement ne d\'{e}pend pas du chemin. On
utilise \`{a} cette fin le th\'{e}or\`{e}me de monodromie. Nous
allons tout d'abord modifier
$\overline{\mathbb{C}}\backslash\{z_1,...,z_m\}$ pour obtenir une
surface simplement connexe. On proc\`{e}de comme suit : Soit
$z^*\in \overline{\mathbb{C}}\backslash\{z_1,...,z_m\}$ un point
arbitraire et consid\'{e}rons $m+1$ chemins
$\mathcal{L}_1,...,\mathcal{L}_{m+1}$ de $z^*$ jusque
$z_1,...,z_{m+1}=\infty$ respectivement. On suppose que chaque
$\mathcal{L}_i$ ne se recoupe pas et que
$$\mathcal{L}_i\cap \mathcal{L}_j=\{z^*\},$$
pour tout $i\neq j$, ($i,j=1,...,m+1$). En faisant des coupures le
long de ces chemins, on obtient une surface
$$\sigma=\left(\overline{\mathbb{C}}\backslash\{z_1,...,z_m\}\right)\backslash\bigcup_i\mathcal{L}_i,$$
hom\'{e}omorphe \`{a} un disque, donc simplement connexe.

\begin{prop}
Sur la surface $\sigma$, le prolongement analytique des $w_k(z)$
ne d\'{e}pend pas du chemin.
\end{prop}
\emph{D\'{e}monstration}: Il suffit d'utiliser le th\'{e}or\`{e}me
de monodromie et la construction de la surface $\sigma$. $\square$

Consid\'{e}rons une coupure le long de $\mathcal{L}_j$ et la
solution $w_k(z)$. On tourne autour de $z_j$ pour atteindre
l'autre coupure. Ceci revient \`{a} transformer la solution
$w_k(z)$ en une solution $w_{l_k}(z)$. Donc pour chaque $j$
fix\'{e}, on a une permutation $\pi_j$ qui transforme $k$ en
$l_k$. Prenons $n$ copies $\sigma_1,...,\sigma_n$ de $\sigma$ et
identifions le bord $B_j$ de $\sigma_i$ avec le bord $A_j$ de
$\sigma_{k=\pi_j(i)}$. On identifie ainsi tous les bords deux
\`{a} deux et on obtient une surface compacte. Nous allons montrer
que cette surface, not\'{e}e $X$, est connexe.

\begin{thm}
La surface de Riemann $X$ obtenue est connexe.
\end{thm}
\emph{D\'{e}monstration}: Supposons que $X$ n'est pas connexe. Il
est donc possible de num\'{e}roter les copies
$\sigma_1,...,\sigma_n$ de $\sigma$ de sorte que les $k$
premi\`{e}res $\sigma_1,...,\sigma_k$, $k<n$, soient reli\'{e}es
entre elles, formant ainsi une des composantes connexes de $X$. Du
point de vue des permutations, cel\`{a} signifie que si $\pi_j$
envoit les indices $i=1,...,n$ sur les indices $j_1,...,j_n$,
alors $j_1,...,j_k$ est une permutation des indices $1,...,k$ et
$j_{k+1},...,j_n$ est une permutation des indices $k+1,...,n$. On
consid\`{e}re
\begin{eqnarray}
P(w,z)&=&a_0(z)\prod_{i=1}^n(w-w_i(z)),\nonumber\\
G(w,z)&=& (w-w_1(z))...(w-w_k(z)),\nonumber\\
&=&w^k-E_1w^{k-1}+E_2w^{k-2}-...+(-1)^kE_k,\nonumber
\end{eqnarray}
o\`{u}
\begin{eqnarray}
E_1&=&\sum_{i=1}^kw_i(z),\nonumber\\
E_2&=& \sum_{\overset{i,j=1}{i< j}}^kw_iw_j,\nonumber\\
&\vdots&\nonumber\\
E_k&=&\sum_{\overset{i_1,...,i_j\in \{1,...,k\}}{i_1<...<
i_j}}^kw_{i_1}...w_{i_j},\nonumber
\end{eqnarray}
sont des polyn\^{o}mes sym\'{e}triques en $w_1,...,w_k$. Ce sont
des fonctions rationnelles. En effet, la permutation $\pi_j
(j=1,...,m+1)$ transforme l'ensemble $\{w_1,...,w_k\}$ en lui
m\^{e}me. Comme $E_k$ sont des polyn\^{o}mes sym\'{e}triques en
$w_1,...,w_k$, alors cette permutation transforme aussi $E_k$ en
$E_k$. Pour passer les coupures $L_j$, on applique la permutation
$\pi_j$ qui conserve la valeur de $E_k$, donc $E_k(z)$ est
univalu\'{e}e sur $\sigma$. En outre $E_k$ est holomorphe sauf
peut-\^{e}tre aux points de branchement o\`{u} elle pourrait
\'{e}ventuellement avoir des p\^{o}les. Donc $E_k$ est
m\'{e}romophe et par cons\'{e}quent rationnelle. D\`{e}s lors
$G(w,z)$ est une fonction rationnelle en $w$ et $z$. En
multipliant $G(w,z)$ par le d\'{e}nominateur commun des $E_k$, on
obtient un polyn\^{o}me $Q(w,z)$ de degr\'{e} $k$ en $w$ et dont
les racines sont les fonctions $w_1,...,w_k$. D'apr\`{e}s la
proposition 10.3, appendice 10.2), il existe deux polyn\^{o}mes
$U$ et $V$ tels que :
$$UP+VQ=\mbox{polyn\^{o}me (r\'{e}sultante) } R(z) \mbox{ en } z.$$
Par hypoth\`{e}se $P$ est irr\'{e}ductible. Comme $\mbox{deg
}Q<\mbox{deg } P$, $P$ et $Q$ sont premiers entre eux, donc $P$ et
$Q$ n'ont pas de facteur commun et par cons\'{e}quent $R(z)\neq
0$. Mais que vaut l'expression
$$UP+VQ=R(z),\quad \forall z,$$
pour les valeurs particuli\`{e}res $w=w_j(z), 1\leq j\leq k$? On a
$$UP+VQ\mid_{w=w_j(z)}=U.0+V.0=0=R(z).$$
Donc $R(z)=0, \forall z,$ ce qui est contradictoire et le
th\'{e}or\`{e}me est d\'{e}montr\'{e}. $\square$

\section{Formes diff\'{e}rentielles}

\subsection{G\'{e}n\'{e}ralités}

Soient $n,k\in \mathbb{N}$, $M$ une vari\'{e}t\'{e}
diff\'{e}rentiable de dimension $n$ et $U\subset M$ un ouvert.
Soit
$$\omega=\sum_{1\leq
i_1<...<i_k\leq n}f_{i_1,...,i_k}dx_{i_1}\wedge...\wedge
dx_{i_k},$$ une $k$-forme diff\'{e}rentielle sur $U$. Les
$f_{i_1,...,i_k}$ ($1\leq i_1<...<i_k\leq n$) sont des fonctions
de $U$ dans $\mathbb{R}$ (ou $\mathbb{C}$), de classe
$\mathcal{C}^\infty$. On d\'{e}finit le produit ext\'{e}rieur de
$\omega$ et $\lambda$ (une l-forme diff\'{e}rentielle dans $U$) en
posant
$$\omega\wedge \lambda=
\sum_{1\leq i_1,...,i_k,j{_1},...,j{_l}\leq n}
f_{i_1,...,i_k}g_{j_1,...,j{_l}}dx_{i_1}\wedge...\wedge
dx_{i_k}\wedge dx_{j{_1}}\wedge...\wedge dx_{j{_l}}.$$ C'est une
$(k+l)$-forme diff\'{e}rentielle dans $U$. On v\'{e}rifie
ais\'{e}ment que :
\begin{eqnarray}
k+l>n&\Longrightarrow&\omega\wedge \lambda=0,\nonumber\\
(\omega\wedge \lambda)\wedge \eta&=&\omega\wedge (\lambda\wedge
\eta),\nonumber\\
(\omega+ \eta)\wedge \lambda&=&(\omega\wedge\lambda)+(\eta\wedge
\lambda),\nonumber\\
\omega\wedge \lambda&=&(-1)^{kl}(\lambda\wedge \omega).\nonumber
\end{eqnarray}
On d\'{e}finit la diff\'{e}rentielle ext\'{e}rieure de $\omega$ en
posant
$$d\omega=
\sum_{1\leq i_1,...,i_k\leq n} df_{i_1,...,i_k}\wedge
dx_{i_1}\wedge...\wedge dx_{i_k}.$$ C'est une $(k+1)$-forme
diff\'{e}rentielle dans $U$. Dans le cas d'une $1$-forme
$\omega=\sum_{i=1}^n f_idx,$ on a $$d\omega=\sum_{1\leq i,k\leq n}
\left(\frac{\partial f_j}{\partial x_i}-\frac{\partial
f_i}{\partial x_j}\right)dx_i\wedge dx_j.$$ On v\'{e}rifie les
formules suivantes :
\begin{eqnarray}
d(a\omega+b\lambda)&=&ad\omega+bd\lambda,\quad (a,b\in \mathbb{R})\nonumber\\
d(\omega\wedge
\lambda)&=&(d\omega\wedge\lambda)+(-1)^k(\omega\wedge
d\lambda),\quad k=\mbox{deg } \omega\nonumber\\
d(d\omega)&=&0.\nonumber
\end{eqnarray}
On dit que $\omega$ est ferm\'{e}e (ou un cocycle) si
$$d\omega=0.$$ En particulier, une $1$-forme $\omega=\sum_{i=1}^n
f_idx$ est ferm\'{e}e si et seulement si
$$\frac{\partial f_i}{\partial x_j}=\frac{\partial f_j}{\partial
x_i},\quad \forall 1\leq i,j\leq n.$$ On dit que $\omega$ est
exacte (ou cohomologue à 0) s'il existe une $(k-1)$-forme
diff\'{e}rentielle $\lambda$ dans $U$ telle que :
$$\omega=d\lambda.$$ En particulier, une $1$-forme diff\'{e}rentielle
$\omega=\sum_{i=1}^n f_idx$ est exacte s'il existe une application
$h:U\rightarrow \mathbb{R}$ (de classe $\mathcal{C}^1$) telle que
: $$f_i=\frac{\partial h}{\partial x_i}.$$ Toute forme
diff\'{e}rentielle exacte est ferm\'{e}e. La r\'{e}ciproque est
fausse en g\'{e}n\'{e}ral et elle est vraie en degr\'{e} $\geq 1$
si l'ouvert $U$ est \'{e}toil\'{e} (lemme de Poincar\'{e}). Soient
$\Omega^k(M)$ l'espace vectoriel r\'{e}el des $k$-formes
diff\'{e}rentielles de classe $\mathcal{C}^\infty$ sur $M$ et
$d:\Omega^k(M)\longrightarrow \Omega^{k+1}(M)$, la
diff\'{e}rentielle ext\'{e}rieure.
\begin{defn}
On appelle groupe de cohomologie de la vari\'{e}t\'{e} $M$,
l'espace vectoriel r\'{e}el
\begin{eqnarray}
H^k(M,\mathbb{R})&=&\frac{\ker \left[d:\Omega^k(M)\longrightarrow
\Omega^{k+1}(M)\right]}{\mbox{Im}
\left[d:\Omega^{k-1}(M)\longrightarrow
\Omega^k(M)\right]},\nonumber\\
&=&\frac{\{k-\mbox{formes diff\'{e}rentielles ferm\'{e}es sur}
M\}}{\{k-\mbox{formes diff\'{e}rentielles exactes sur}
M\}}.\nonumber
\end{eqnarray}
C'est le groupe de cohomologie de De Rham que l'on d\'{e}signe
aussi par $H^k_{DR}(M,\mathbb{R})$.
\end{defn}
Un \'{e}l\'{e}ment de ce groupe est une classe d'\'{e}quivalence
de formes ferm\'{e}es diff\'{e}rent l'une de l'autre par une
diff\'{e}rentielle :
$$\omega_1\sim \omega_2\quad \mbox{si}\quad
\omega_1-\omega_2=d\lambda.$$
\begin{Exmp}
$H^0(M,\mathbb{R})$ est un espace vectoriel de dimension finie
\'{e}gal au nombre de composantes connexes de $M$. En effet, nous
n'avons ici que des 0-formes diff\'{e}rentielles, i.e., des
fonctions $f(x)$ sur $M$. Il n'y a pas de formes
diff\'{e}rentielles exactes. D\`{e}s lors, $H^0(M,\mathbb{R})=\{f:
f \mbox{ est ferm\'{e}e}\}$. Comme $df(x)=0$, alors dans toute
carte $(U,x_1,...,x_n)$ de $M$, on a $$\frac{\partial f}{\partial
x_1}=\frac{\partial f}{\partial x_2}=...=\frac{\partial
f}{\partial x_n}=0.$$ Par cons\'{e}quent, $f(x)=$ constante
localement, i.e., $f(x)=$ constante sur chaque composante connexe
de $M$. Donc le nombre de composantes connexes de $M$ est la
dimension en question.
\end{Exmp}

Le lemme de Poincar\'{e} peut-\^{e}tre vu comme un
th\'{e}or\`{e}me d'annulation de la cohomologie :
$H^k(U,\mathbb{R})=0$, $k\geq 1$ o\`{u} $U\subset \mathbb{R}^n$
est un ouvert \'{e}toil\'{e}.

Soit $I=[0,1]$, $I^k=I\times ...\times I$ (k-fois). Un
$k$-simplexe (de classe $\mathcal{C}^r$, $r\geq 1$) dans $U$ est
une application $\varphi : I^k\longrightarrow U$, qui est de
classe $\mathcal{C}^r$ sur $I^k$. Comme $I^k$ est le volume
form\'{e} par $k$ vecteurs ind\'{e}pendants dans $U$,
l'application $\varphi$ est donc une d\'{e}formation de $I^k$. Par
exemple, un $0$-simplexe est un point. Un $1$-simplexe dans
$\mathbb{R}^3$ est une courbe dans $\mathbb{R}^3$. Un $2$-simplexe
dans $\mathbb{R}^3$ est une surface dans $\mathbb{R}^3$
hom\'{e}omorphe \`{a} $I^2$ (un triangle). Un $3$-simplexe dans
$\mathbb{R}^3$ est une boule, un cube, un volume dans
$\mathbb{R}^3$ hom\'{e}omorphe \`{a} $I^3$ (un t\'{e}tra\`{e}dre).

Rappelons que pour int\'{e}grer une $k$-forme diff\'{e}rentielle
sur une vari\'{e}t\'{e} $M$, il faut que celle-ci soit orientable.
Cel\`{a} signifie que $M$ v\'{e}rifie l'une des conditions
\'{e}quivalentes suivantes : $(i)$ $M$ est munie d'une forme
volume, i.e., une $k$-forme diff\'{e}rentielle qui ne s'annule
nulle part. $(ii)$ Il existe sur $M$ un atlas tel que le jacobien
de tout changement de carte soit $>0$. Par exemple, $\mathbb{R}^n$
est orient\'{e}e par la forme volume $dx_1\wedge...\wedge dx_n$.
Le cercle $S^1$ est orient\'{e} par $d\theta$. Le tore
$T^2=S^1\times S^1$ est orient\'{e} par la forme volume
$d\theta\wedge d\varphi$. Toutes les vari\'{e}t\'{e}s holomorphes
sont orientables. La sph\`{e}re $S^2$, l'espace projectif
$\mathbb{RP}^n$ ($n$ paire), la bande de M\"{o}bius ne sont pas
orientables.

Soient $\omega$ une $k$-forme diff\'{e}rentielle $U\subset M$
(orientable) et $\varphi$ un $k$-simplexe dans $U$ de classe
$\mathcal{C}^1$. L'int\'{e}grale de $\omega$ sur $\varphi$ est
d\'{e}finie par
$$\int_\varphi \omega=\int_{I^k}\omega (\varphi).$$

On suppose que $k\geq 2$ et soient $\pi$ une permutation de
$(1,...,k)$,
$$\varphi : I^k\longrightarrow U, (u_1,...,u_k)\longmapsto \varphi
(u_1,...,u_k),$$ un $k$-simplexe dans $U$, et
$$\varphi_\pi : I^k\longrightarrow U, (u_1,...,u_k)\longmapsto
\varphi_\pi (u_1,...,u_k)=\varphi (u_{\pi(1)},...,u_{\pi(1)}),$$
un $k$-simplexe dans $U$. Pour toute $k$-forme diff\'{e}rentielle
$\omega$ dans $M$, on a
$$\int_{\varphi_\pi}\omega=\mbox{sign }\pi \int_\varphi \omega,$$
o\`{u} $\mbox{sign }\pi$ d\'{e}signe la signature de la
permutation $\pi$. Cette relation montre que l'ensemble des
$\varphi_\pi$ lorsque $\pi$ parcourt les permutations de
$(1,2,...,k)$ peut \^{e}tre divis\'{e} en deux classes: La
premi\`{e}re correspond au cas o\`{u} $\pi$ est paire ($\mbox{sign
}\pi=1$); $\varphi_\pi$ et $\varphi$ ont m\^{e}me orientation et
on posera dans ce cas $\varphi_\pi=\varphi$. La seconde correspond
au cas o\`{u} $\pi$ est impaire ($\mbox{sign }\pi=-1$);
$\varphi_\pi$ et $\varphi$ ont des orientations oppos\'{e}es et on
posera dans ce cas $\varphi_\pi=\varphi_-$. Pour cette seconde
classe, $\pi$ est une transposition, i.e., deux \'{e}l\'{e}ments
seulement sont permut\'{e}s. Par exemple, si $k=2$ et $n=3$, on a
$\varphi_-(u_1,u_2)=\varphi(u_2,u_1)$. Nous avons suppos\'{e} que
$k\geq2$. Dans le cas où $k=1$, alors $\varphi_-$ s'obtient par la
formule $\varphi_-(u)=\varphi(1-u)$.

Un $k$-complexe (de classe $\mathcal{C}^r$, $r\geq 1$) dans $U$
est une famille finie $\Phi=(\varphi_1,...,\varphi_m)$ de
$k$-simplexes $\varphi_j$, $1\leq j\leq m$, dans $U$ (de classe
$\mathcal{C}^r$, $r\geq 1$). G\'{e}om\'{e}triquement, un complexe
peut se visualiser comme une r\'{e}union de surfaces, celles
d\'{e}finies par les simplexes le composant.

Si $\Phi=(\varphi_1,...,\varphi_m)$ est un $k$-complexe (de classe
$\mathcal{C}^1$) dans $U$ et si $\omega$ est une $k$-forme
diff\'{e}rentielle dans $U$, alors l'int\'{e}grale de $\omega$ sur
$\Phi$ est d\'{e}finie par
$$\int_\Phi\omega=\sum_{j=1}^m\int_{\varphi_j}\omega.$$
En particulier, un $k$-complexe $\Phi=(\varphi_1,...,\varphi_m)$
dans $U$ tel que : $$\varphi_j(1)=\varphi_{j+1}(0),\quad 1\leq
j\leq m-1,$$ est un chemin dans $U$ par morceaux. Si en outre
$\varphi_m(1)=\varphi_1(0)$, alors $\Phi$ est un cycle ou chemin
ferm\'{e}.

Soient $\omega$ une $1$-forme diff\'{e}rentielle exacte dans $U$,
$\Phi=(\varphi_1,...,\varphi_m)$ et $\Psi=(\psi_1,...,\psi_m)$
deux chemins dans $U$. Si $\omega=df$, alors
$$\int_\Phi=f\left(\varphi_m(1)\right)-f\left(\varphi_1(0)\right).$$
Si $\varphi_1(0)=\psi_1(0)$ et $\varphi_m(1)=\psi_m(1)$, alors
$$\int_\Phi \omega=\int_\Psi \omega.$$ Enfin si $\Phi$ est un
chemin ferm\'{e}, alors $$\int_\Phi \omega=0.$$

Le bord d'un $k$-simplexe $\varphi : I^k\longrightarrow U, k\geq
2$, est le $(k-1)$-complexe, not\'{e} $\partial \varphi$,
d\'{e}fini par
$$\partial \varphi=\left(\varphi_{\mbox{sign}(-1)^{j+\alpha}}^{j, \alpha} : 1\leq j\leq k, \alpha=0, 1\right),$$
o\`{u}
\begin{eqnarray}
\varphi_{+}^{j, \alpha}\equiv \varphi^{j, \alpha} :
I^{k-1}&\longrightarrow& U,\nonumber\\
(u_1,...,u_{k-1})&\longmapsto& \varphi^{j,
\alpha}(u_1,...,u_{k-1})=\varphi
(u_1,...,u_{j-1},\alpha,u_j,...,u_{k-1}),\nonumber
\end{eqnarray}
et $\varphi_{+}^{j, \alpha}=\left(\varphi^{j, \alpha}\right)_-$.

Lors de la d\'{e}termination du bord d'un simplexe $\varphi$, les
notations suivantes peuvent \^{e}tre utiles pour pr\'{e}ciser
l'image de $\varphi$ ainsi que l'orientation. Soient
$x_0,x_1,...,x_k \in U$, $k+1$ points et
$$\varphi : I^k\longrightarrow U, u\longmapsto
\varphi(u)=x_0+\sum_{j=1}^ku_j(x_j-x_0).$$ L'application $\varphi$
d\'{e}finit un $k$-simplexe dans $U$ et $\mbox{Im }\varphi$ est un
parall\'{e}lipip\`{e}de (à $k$ dimensions) construit sur les $k$
segments joignant $x_0$ à $x_j$, $1\leq j\leq k$. Ces
$k$-simplexes sont appel\'{e}s $k$-parall\'{e}lotopes orient\'{e}s
et seront not\'{e}s : $[x_0, x_1,...,x_k]$. On a $x_0=\varphi(0)$
et $x_j=\varphi(e_j)$, $1\leq j\leq k$, o\`{u} $e_j$ d\'{e}signe
le $j^{\mbox{ième}}$ vecteur de base de $\mathbb{R}^k$. Pour
$k=1$, on note $[x_0, x_1]$ et $x_0
\underset{orientation}{\longrightarrow } x_1$.
\begin{Exmp}
Soit
$$\varphi:I^2\longrightarrow \mathbb{R}^2, (u_1,u_2)\longmapsto
(u_1,u_2),$$ le $2$-simplexe identit\'{e} (injection canonique).
On a
\begin{eqnarray}
\partial \varphi&=&\left(\varphi_{\mbox{sign}(-1)^{j+\alpha}}^{j, \alpha} : 1\leq j\leq 2, \alpha=0, 1\right),\nonumber\\
&=&\left(\varphi_-^{1,0}, \varphi^{1,1}, \varphi^{2,0},
\varphi_-^{2,1}\right),\nonumber
\end{eqnarray}
o\`{u}
\begin{eqnarray}
\varphi_-^{1,0}(u)&=&\left(\varphi^{1,0}\right)_-(u),\nonumber\\
&=&\varphi^{1,0}(1-u),\nonumber\\
&=&\varphi(0,1-u),\nonumber\\
&=&(0,1-u),\nonumber\\
&=&[(0,1),(0,0)](u),\nonumber\\
&=&[e_2,0](u),\nonumber
\end{eqnarray}
\begin{eqnarray}
\varphi^{1,1}(u)&=&\varphi(1,u),\nonumber\\
&=&(1,u),\nonumber\\
&=&[(1,0),(1,1)](u),\nonumber\\
&=&[e_1,e_1+e_2](u),\nonumber
\end{eqnarray}
\begin{eqnarray}
\varphi^{2,0}(u)&=&\varphi(u,0),\nonumber\\
&=&(u,0),\nonumber\\
&=&[(0,0),(1,0)](u),\nonumber\\
&=&[0,e_1](u),\nonumber
\end{eqnarray}
\begin{eqnarray}
\varphi_-^{2,1}(u)&=&\left(\varphi^{2,1}\right)_-(u),\nonumber\\
&=&\varphi^{2,1}(1-u),\nonumber\\
&=&(1-u,1),\nonumber\\
&=&[(1,1),(0,1)](u),\nonumber\\
&=&[e_1+e_2,e_2](u).\nonumber
\end{eqnarray}
Comme $I^2=[0,1]\times [0,1]$, alors $\varphi (I^2)$ est le
carr\'{e} construit sur les segments joignant $0$ à $e_1$, $e_1$
\`{a} $e_1+e_2$, $e_1+e_2$ à $e_2$ et $e_2$ à $0$. On obtient
$$(\partial \varphi)(I^2)=\bigcup_{\underset{\alpha =0,1}{1\leq j\leq 2}}\left(\varphi_{\mbox{sign}(-1)^{j+\alpha}}^{j,
\alpha}(I)\right)=\mbox{fr }(\varphi(I^2)).$$
\end{Exmp}
Le bord d'un $k$-complexe $\Phi=(\varphi_1,...,\varphi_m)$ est le
$(k-1)$-complexe d\'{e}fini par
\begin{eqnarray}
\partial \Phi&=&(\partial \varphi_1,...,\partial \varphi_m),\nonumber\\
&=&\left(\varphi_{\mbox{sign}(-1)^{l+\alpha}}^{j, l, \alpha} :
1\leq j\leq m, 1\leq l\leq k, \alpha=0, 1\right).\nonumber
\end{eqnarray}

Soient $M_1$, $M_2$ des vari\'{e}t\'{e}s diff\'{e}rentiables de
dimension $m_1$, $m_2$ respectivement et $U_1\subset M_1$,
$U_2\subset M_2$ des ouverts. Pour toute application
diff\'{e}rentiable $g : U_1\longrightarrow U_2$, et toute
$k$-forme diff\'{e}rentielle dans $U_2$, on peut d\'{e}finir une
$k$-forme diff\'{e}rentielle dans $U_1$ (appel\'{e}e le pull-back
par $g$ ou image inverse ou encore transpos\'{e}e de $\omega$ par
$g$) en posant
$$g^*\omega =\sum_{1\leq i_1,...,i_k\leq m_2}\left(f_{i_1,...,i_k}\circ
g\right)dg_{i_1}\wedge...\wedge dg_{i_k},$$ o\`{u}
$$dg_{i_l}=\sum_{j=1}^{m_1}\frac{\partial g_{i_l}}{\partial
y_j}dy_j,$$ sont des $1$-formes dans $U_1$. Notons que $g^*$ est
un op\'{e}rateur lin\'{e}aire de l'espace des $k$-formes sur $N$
dans l'espace des $k$-formes sur $U_1$ (l'ast\'{e}risque indique
que $g^*$ op\`{e}re dans le sens inverse de $g$). Soit $\omega$
est une $k$-forme diff\'{e}rentielle dans $U_2$ et $\lambda$ une
$l$-forme diff\'{e}rentielle dans $U_2$. Alors si $k=l$,
$$g^*(\omega+\lambda)=g^*\omega+g^*\lambda,$$
et
$$g^*(\omega \wedge\lambda)=g^*\omega\wedge g^*\lambda.$$
Si $h : U_2\rightarrow \mathbb{R}$ est une application continue,
$$g^*(h\omega)=(h\circ g)g^*\omega.$$
Si $\omega$ est de classe $\mathcal{C}^1$ dan $U_2$ et $g$ de
classe $\mathcal{C}^2$ dans $U_1$,
$$g^*(d\omega)=d(g^*\omega).$$
Si $M_3$ est une autre vari\'{e}t\'{e} diff\'{e}rentiable,
$U_3\subset M_3$ un ouvert et $h : U_3 \rightarrow U_1$ une
application de classe $\mathcal{C}^1$, alors
$$(g \circ h)^*\omega=h^*(g^*\omega).$$
On d\'{e}duit de ces propri\'{e}t\'{e}s que si $g :
U_1\longrightarrow U_2$ est une application diff\'{e}rentiable,
alors il y a des applications lin\'{e}aires induites
$$g^* : H^k(U_2,\mathbb{R})\longrightarrow H^k(U_1,\mathbb{R}),$$
telles que :
$$g^*([\omega]\wedge[\lambda])=g^*[\omega]\wedge g^*[\lambda].$$
En outre, si $g$ est un diff\'{e}omorphisme local, alors $g^*$ est
un isomorphisme d'alg\`{e}bres (les groupes de cohomologie donnent
donc des invariants diff\'{e}rentiables).

On montre que si $\varphi$ est un $k$-simplexe dans $M$ de classe
$\mathcal{C}^1$, $\omega$ une $k$-forme diff\'{e}rentielle dans
$M$, $\tau^k : I^k\rightarrow \mathbb{R}^k, u\mapsto \tau^k(u)=u$
un $k$-simplexe identit\'{e} dans $I^k\subset \mathbb{R}^k$
(injection canonique de classe $\mathcal{C}^\infty$), alors
$$\int_\varphi \omega=\int_{\tau^k}\varphi^*\omega.$$
On en d\'{e}duit que si $\omega$ est une $(k-1)$-forme
diff\'{e}rentielle de classe $\mathcal{C}^1$ dans $M$ (orientable)
et $\Phi=(\varphi_1,...,\varphi_m)$ un $k$-complexe dans $M$ de
classe $\mathcal{C}^2$, alors on a la formule de Stokes-Cartan
$$\int_{\partial \Phi}\omega=\int_\Phi d\omega.$$

Soit $A$ un anneau unitaire. On notera $C_k$ le $A$-module libre
engendr\'{e} par tous les $k$-simplexes dans un complexe $\Phi$.
Un \'{e}l\'{e}ment de $C_k$ est appel\'{e} une $k$-cha\^{i}ne dans
le complexe $\Phi$. C'est une somme finie formelle de la forme
$$c_k=\sum_k\alpha_k\sigma_k,$$
o\`{u} $\sigma_k$ est un $k$-simplexe et $\alpha_k\in A$. Le bord
$\partial c_k$ d'une $k$-cha\^{i}ne est d\'{e}finie par
$$\partial c_k=\sum_k\alpha_k\partial\sigma_k.$$
Un cycle est une cha\^{i}ne $c_k$ telle que : $\partial c_k=0$. Un
bord est une cha\^{i}ne $c_k$ telle qu'il existe une cha\^{i}ne
$c_{k+1}$ avec $\partial c_{k+1}=c_k$. Par analogie avec les
formes diff\'{e}rentielles, on peut dire qu'un cycle est une
cha\^{i}ne ferm\'{e}e et qu'un bord est une cha\^{i}ne exacte. On
montre que pour toute $k$-cha\^{i}ne $c_k$, le bord $\partial c_k$
est une $(k-1)$-cha\^{i}ne et que $\partial (\partial c_k)=0$. En
outre, les $k$-cha\^{i}nes forment un groupe ab\'{e}lien en
introduisant une loi d'addition sur les cha\^{i}nes comme suit :
si $c_k=\sum_k\alpha_k\sigma_k$, $c'_k=\sum_k\alpha'_k\sigma_k$,
alors
$$c_k+c'_k=\sum_k(\alpha_k+\alpha'_k)\sigma_k.$$
Les cycles forment un groupe not\'{e}
$$Z_k=\ker (\partial :C_k\longrightarrow C_{k-1}).$$
Deux cycles $c'_k$ et $c''_k$ sont \'{e}quivalents ou homologues
si et seulement si $c'_k-c''_k=\partial c_{k+1}$. De m\^{e}me, les
bords forment aussi un groupe not\'{e}
$$B_k=\mbox{Im} (\partial :C_{k+1}\longrightarrow C_k).$$
On pose $B_0=0$. On a les inclusions : $B_k\subset Z_k\subset
C_k$.
\begin{defn}
On appelle groupe d'Homologie, not\'{e} $H_k$, le quotient
$Z_k/B_k$.
\end{defn}
Il est clair qu'un bord est un cycle. Mais tout cycle n'est pas un
bord.
\begin{Exmp}
Soient $a$ et $b$ deux cycles ind\'{e}pendants dans $H_1(T)$. Ces
cycles forment une base d'homologie du tore $T$. On a
$$H_1(T)\overset{(1)}{=}Z_1(T)/B_1(T)\overset{(2)}{=}\mathbb{Z}/2\mathbb{Z}.$$
Les deux groupes $(1)$ et $(2)$ ont m\^{e}me structures; ils sont
engendr\'{e}s par deux \'{e}l\'{e}ments $a$ et $b$.
\end{Exmp}

\subsection{Formes diff\'{e}rentielles sur les surfaces de Riemann
(Diff\'{e}rentielles ab\'{e}liennes)}

Une forme diff\'{e}rentielle\footnote{On omettra de pr\'{e}ciser
qu'il s'agit d'une $1$-forme diff\'{e}rentielle.} sur une surface
de Riemann $X$ s'\'{e}crit
$$\omega =f(\tau) d\tau,$$
o\`{u} $\tau$ est le param\`{e}tre local et $f$ une fonction
complexe de $\tau$.
\begin{defn}
On dit que $\omega $ est une diff\'{e}rentielle ab\'{e}lienne si
$f(\tau)$ est une fonction m\'{e}romorphe sur $X$, holomorphe si
$f(\tau)$ est une fonction holomorphe sur $X$, ayant un p\^{o}le
d'ordre $k$ ou un z\'{e}ro d'ordre $k$ en un point $p$ si
$f(\tau)$ a un p\^{o}le d'ordre $k$ ou un z\'{e}ro d'ordre $k$ en
ce point.
\end{defn}
\begin{defn}
On dit qu'une diff\'{e}rentielle ab\'{e}lienne est de
$1^{\grave{e}re}$esp\`{e}ce si elle est holomorphe sur $X$, de
$2^{\grave{e}me}$esp\`{e}ce si son
r\'{e}sidu\footnote{$\mbox{R\'{e}s} _p \omega=\mbox{Rés} _p
f=a_{-1}$ o\`{u} $p\in X$ et $f(\tau)=\sum_{k=N}^\infty a_kz^k$.
On a aussi, $\mbox{Rés} _p \omega=\int_\gamma \omega$ o\`{u}
$\gamma$ est une courbe ferm\'{e}e sur $X$.} est nul en chaque
point de $X$ et enfin de $3^{\grave{e}me}$esp\`{e}ce si son
r\'{e}sidu est non nul en au moins un point de $X$.
\end{defn}

Soit $( a_{1},\ldots ,a_{g},b_{1},\ldots ,b_{g}) $ une base de
cycles dans le groupe d'homologie $H_{1}( X,\mathbb{Z}) $ de telle
fa\c{c}on que les produits d'intersection de cycles deux \`{a}
deux s'\'{e}crivent : $$\left( a_{j},a_{j}\right) =\left(
b_{j},b_{j}\right) =0,\quad \left( a_{j},b_{k}\right) =\delta
_{jk},\quad 1\leq j,k\leq g.$$ On dira que $( a_{1},\ldots
,a_{g},b_{1},\ldots ,b_{g}) $ est une base symplectique du groupe
d'homologie $H_{1}( X,\mathbb{Z})$. On d\'{e}signe par $X^*$ la
repr\'{e}sentation normale de la surface de Riemann $X$ de genre
$g$, i.e., un polyg\^{o}ne \`{a} $4g$ c\^{o}t\'{e}s identifi\'{e}s
deux \`{a} deux selon le symbole
$a_1b_1a_1^{-1}b_1^{-1}...a_gb_ga_g^{-1}b_g^{-1}$ et peut \^{e}tre
d\'{e}finit \`{a} partir d'une triangulation de la surface $X$.
Les $4g$ c\^{o}t\'{e}s sont nomm\'{e}s dans l'ordre o\`{u} ils se
pr\'{e}sentent sur le bord orient\'{e} du polyg\^{o}ne, $a_1^{-1}$
devant \^{e}tre identifi\'{e} \`{a} $a_1$ apr\`{e}s avoir
renvers\'{e} son orientation, etc...
\begin{prop}
Soit $X$ une surface de Riemann compacte de genre $g$. Soit
$\omega$ une diff\'{e}rentielle ab\'{e}lienne sur $X$. Alors
$$\sum_{p\in X}\mbox{Rés}_p \omega=0.$$
\end{prop}
\emph{D\'{e}monstration}: Soit $X^*$ la repr\'{e}sentation normale
de $X$. On sait que $X^*$ est un polyg\^{o}ne à $4g$ c\^{o}t\'{e}s
identifi\'{e}s deux \`{a} deux. Si l'on parcourt le bord $\partial
X^*$ de ce polyg\^{o}ne, on constate que chaque c\^{o}t\'{e} est
parcouru deux foix, l'un dans le sens de son orientation et
l'autre dans le sens oppos\'{e}. Donc $\int_{\partial
X^*}\omega=0$. Or, d'apr\`{e}s le th\'{e}or\`{e}me des
r\'{e}sidus, on a
$$\int_{\partial X^*}\omega=2\pi i\sum_{p\in
X}\mbox{R\'{e}s}_p\omega,$$ car tous les points de $X^*$ sont des
points de $X$. D'o\`{u}, $\sum_{p\in X}\mbox{R\'{e}s}_p\omega=0$.
$\square$

Soient $X$ une surface de Riemann compacte, $\omega$ une forme
diff\'{e}rentielle ferm\'{e}e, $\triangle$ une cha\^{i}ne sur $X$
et $\partial \triangle$ son bord. D'apr\`{e}s le th\'{e}or\`{e}me
de Stokes-Cartan, on a
$$\int_{\partial \triangle}=\int\int_\triangle d\omega=0.$$
D\`{e}s lors, si $\gamma$ est un chemin contenu dans $X$ et si
$[\gamma]\in H_1(X,\mathbb{Z})$ est la classe d'homologie
contenant $\gamma$, alors l'application
$$\int_{[\gamma]} : H_1(X,\mathbb{Z})\longrightarrow
\mathbb{C},\quad [\gamma]\longmapsto
\int_{[\gamma]}\omega=\int_{\gamma}\omega,$$ est bien d\'{e}finie.
En particulier, si $\omega$ est une forme diff\'{e}rentielle
holomorphe alors l'int\'{e}grale ci-dessus est bien d\'{e}finie
puisque toute forme diff\'{e}rentielle holomorphe $\omega$ est
ferm\'{e}e.

\begin{defn}
On appelle p\'{e}riodes de $\omega $ suivant les cycles
$a_{j},b_{j},$ les nombres $\int_{a_{j}}\omega $ et
$\int_{b_{j}}\omega$.
\end{defn}

Soit $\omega$ une forme diff\'{e}rentielle holomorphe ou de
$2^{\grave{e}me}$esp\`{e}ce sur $X$. Soit
$$\gamma=\sum_{j=1}^g\alpha_ia_i+\sum_{j=1}^g\beta_jb_j,$$
un cycle sur $X$. On d\'{e}duit du th\'{e}or\`{e}me des
r\'{e}sidus, que
$$\int_\gamma\omega=\sum_{j=1}^g\alpha_j\int_{a_j}\omega+\sum_{j=1}^g\beta_j\int_{b_j}\omega.$$

Soit $p_0\in X$, un point fix\'{e}. Pour tout point $p\in X$,
l'int\'{e}grale $\int_{p_0}^p\omega$ est bien d\'{e}finie si et
seulement si toutes les p\'{e}riodes de $\omega$ sont nulles. La
condition n\'{e}cessaire est \'{e}vidente car si
$\int_{p_0}^p\omega$ est bien d\'{e}finie, alors l'int\'{e}grale
de $\omega$ sur tout chemin ferm\'{e} est nulle. Concernant la
condition suffisante, on va utiliser un raisonnement par
l'absurde. Supposons que $\int_{p_0}^p\omega$ n'est pas bien
d\'{e}finie, d\`{e}s lors deux chemins $\gamma_1$ et $\gamma_2$
conduisent à des r\'{e}sultats diff\'{e}rents. On peut donc
trouver un point $p$ appartenant à l'int\'{e}rieur du chemin
ferm\'{e} $\gamma\equiv \gamma_1 \cup \gamma_2$ tel que :
$\mbox{R\'{e}s}_p\omega\neq 0$. Consid\'{e}rons maintenant deux
cycles homologues, un de chaque c\^{o}t\'{e} du point $p$, de
telle mani\`{e}re \`{a} ce que la bande de surface ainsi
d\'{e}termin\'{e}e soit assez fine pour ne contenir qu'un seul
point à r\'{e}sidu non nul. Soit $q$ un point du premier cycle,
$r$ un autre point du second cycle et tra\c{c}ons un chemin qui
relie les deux points $a$ et $r$ mais qui ne passe pas par le
point $p$. En partant du point $q$, on parcourt le premier cycle
et on revient au point $q$; on d\'{e}signe ce chemin parcouru par
$\mathcal{C}_1$. Ensuite, on emprunte le chemin $\mathcal{C}_2$
qui m\`{e}ne du point $q$ au point $r$. Apr\`{e}s, on parcourt le
second cycle et on notera par $\mathcal{C}_3$ le chemin parcouru.
Enfin, on reprend le chemin inverse $\mathcal{C}_4$ qui m\`{e}ne
du point $r$ au point $q$. En d\'{e}signant par
$\mathcal{C}\equiv\mathcal{C}_1\cup \mathcal{C}_2\cup
\mathcal{C}_3\cup \mathcal{C}_4$ le chemin ferm\'{e} de ces
parcourts, on voit que $\int_{\mathcal{C}}\omega\neq 0$ puisque
$\mbox{R\'{e}s}_p\omega\neq 0$. D'un autre c\^{o}t\'{e}, on a
$$\int_\mathcal{C}\omega=\int_{\mathcal{C}_1}\omega+\int_{\mathcal{C}_2}\omega
+\int_{\mathcal{C}_3}\omega+\int_{\mathcal{C}_4}\omega.$$ Par
hypoth\`{e}se les p\'{e}riodes de $\omega$ sont nulles, i.e.,
$$\int_{\mathcal{C}_1}\omega=\int_{\mathcal{C}_3}\omega=0.$$
Or
$$\int_{\mathcal{C}_2}\omega=-\int_{\mathcal{C}_4}\omega,$$
donc $\int_{\mathcal{C}}\omega=0$, ce qui est absurde.

\begin{thm}
L'ensemble\footnote{On \'{e}crira indiff\'{e}ramment $\Omega^1(X)$
ou $H^0(X,\Omega^1)$.}
$$\Omega^1(X)=\{\omega : \omega \mbox{diff\'{e}rentielle holomorphe sur $X$}\},$$
des diff\'{e}rentielles holomorphes sur $X$ est de dimension $g$.
\end{thm}
\emph{D\'{e}monstration}: Montrons tout d'abord que $\dim
\Omega^1(X)\leq g$. Pour cel\`{a}, nous allons montrer qu'il
n'existe pas de suite libre de $g+1$ formes diff\'{e}rentielles
$\omega_0, \omega_1,...,\omega_g$ sur $X$. Donc d\'{e}terminons
des constantes non toutes nulles $c_0,...,c_g$ telles que :
$$\omega\equiv\sum_{k=0}^gc_k\omega_k=0.$$
Soit $\omega_k=u_k+iv_k$. Comme $\omega_k$ est holomophe, alors sa
partie r\'{e}elle $u_k=\mbox{Re }\omega_k$ et sa partie imaginaire
$v_k=\mbox{Im } \omega_k$ sont harmoniques. Consid\'{e}rons la
diff\'{e}rentielle r\'{e}elle
$$\theta\equiv\sum_{k=0}^g(\xi_ku_k+\eta_kv_k),\quad (\xi_k ,\eta_k\in
\mathbb{R}),$$ et d\'{e}terminons $\xi_k ,\eta_k$ de sorte que
toutes les p\'{e}riodes de $\theta$ soient nulles, i.e.,
\begin{eqnarray}
&&\sum_{k=0}^g(\xi_k\int_{a_1}u_k+\eta_k\int_{a_1}v_k)=0,\nonumber\\
&&\quad\vdots\nonumber\\
&&\sum_{k=0}^g(\xi_k\int_{a_g}u_k+\eta_k\int_{a_g}v_k)=0,\nonumber\\
&&\sum_{k=0}^g(\xi_k\int_{b_1}u_k+\eta_k\int_{b_1}v_k)=0,\nonumber\\
&&\quad\vdots\nonumber\\
&&\sum_{k=0}^g(\xi_k\int_{b_g}u_k+\eta_k\int_{b_g}v_k)=0.\nonumber
\end{eqnarray}
Ce sont $2g$ \'{e}quations à $2g+2$ inconnues, il existe donc une
solution $\xi_k ,\eta_k$ non triviale. Toutes les p\'{e}riodes de
$\theta$ sont nulles et d'apr\`{e}s ce qui pr\'{e}c\`{e}de, la
fonction $\varphi(p)\equiv\int_0^p\theta$ est bien d\'{e}finie.
Cette fonction \'{e}tant born\'{e}e et harmonique, alors
d'apr\`{e}s le principe du maximum, $\varphi$=
constante\footnote{En effet, la fonction $\varphi$ n'a ni maximum
ni minimum car sinon elle devrait l'atteindre sur le bord de la
surface de Riemann $X$. Or $X$ n'a pas de bord et comme $\varphi$
est born\'{e}e, alors $\varphi$= constante.}. Notons que
$\varphi(0)=0$, donc $\varphi\equiv 0$ et par cons\'{e}quent
$\theta=0$. Montrons maintenant que :
$$\sum_{k=0}^gc_k\omega_k=0.$$ En effet, posons
$c_k=\xi_k-i\eta_k$, d'o\`{u}
\begin{eqnarray}
\int_0^p\omega&=&\int_0^p\sum_{k=0}^g(\xi_k-i\eta_k)(u_k+iv_k),\nonumber\\
&=&\int_0^p\sum_{k=0}^g(\xi_ku_k+\eta_kv_k)+i\int_0^p\sum_{k=0}^g(\xi_kv_k-\eta_ku_k),\nonumber\\
&=&\int_0^p\mbox{Re }\omega+i\int_0^p\mbox{Im } \omega.\nonumber
\end{eqnarray}
Or $$\sum_{k=0}^g(\xi_ku_k+\eta_kv_k)=\theta=0,$$ donc
$\int_0^p\mbox{Re }\omega=0$. D'apr\`{e}s les conditions de
Cauchy-Riemann, on a aussi $\int_0^p\mbox{Im }\omega=0$. Donc
$\int_0^p\omega=0$, $\forall p$ et par cons\'{e}quent $\omega=0$.
Il reste \`{a} prouver que : $\dim \Omega^1(X)\geq g$. On se
contente ici de donner les \'{e}tapes essentielles.
Consid\'{e}rons l'espace d'Hilbert $L^2(X)$ des formes
diff\'{e}rentielles mesurables sur $X$ et $C^k(X)$ l'ensemble des
formes diff\'{e}rentielles $\zeta(\tau)d\tau$ avec $\zeta \in
\mathcal{C}^k$ et $\tau$ un param\`{e}tre local. Soit $\omega \in
C^1(X)\cap L^2(X)$. On montre que la forme $\omega$ se
d\'{e}compose de mani\`{e}re unique comme suit
$$\omega=\omega_h+df,$$
o\`{u} $\omega_h$ est une diff\'{e}rentielle harmonique et $df \in
L^2(X)$. Si $H(X)$ est l'ensemble des formes diff\'{e}rentielles
harmoniques sur $X$ alors $\dim H\geq 2g$ et $$H(X)=\Omega^1\oplus
\overline{\Omega}^1.$$ D\`{e}s lors, $\dim H(X)=2\dim\Omega^1$ et
par cons\'{e}quent $\dim\Omega^1(X)\geq g$, ce qui ach\`{e}ve la
d\'{e}monstration. $\square$

On d\'{e}duit de l'\'{e}quation $$\frac{\partial F}{\partial
z}dz+\frac{\partial F}{\partial w}dw=0,$$ sur la surface de
Riemann $X : F(w,z)=0$, que les diff\'{e}rentielles rationnelles
s'\'{e}crivent sur $X$ :
$$\omega=\frac{P(w,z)}{\frac{\partial F}{\partial w}}dz=-\frac{P(w,z)}{\frac{\partial F}{\partial z}}dw,$$
o\`{u} $P(w,z)$ une fonction rationnelle arbitraire. En fait, la
formule du r\'{e}sidu de Poincar\'{e} affirme que toutes les
diff\'{e}rentielles holomorphes sur la surface de Riemann $X$
s'\'{e}crivent sous la forme ci-dessus avec pour $P(w,z)$ un
polyn\^{o}me de degr\'{e} $\leq n-3$.

\begin{Exmp}
La diff\'{e}rentielle
$$\omega =\frac{dz}{\sqrt{z\left( z-1\right) \left( z-\lambda \right)
}},$$ est l'unique diff\'{e}rentielle holomorphe sur la courbe
elliptique $X$ d'\'{e}quation :
$$F(w,z) =w^{2}-z(z-1)(z-\lambda)=0,\text{ }\lambda \neq
0,1.$$ En effet, les points de branchements de $X$ sont : $0, 1,
\lambda, \infty$. Comme $g(X)=1$, alors d'apr\`{e}s le
th\'{e}or\`{e}me 3.7, il existe une seule diff\'{e}rentielle
holomorphe sur $X$. La formule du r\'{e}sidu de Poincar\'{e}
s'\'{e}crit dans ce cas
$$\omega=\frac{z^{k}w^{j}}{2\sqrt{z(z-1)(z-\lambda)}}dz.$$
Prenons $j=0$, d'o\`{u}
$$\omega=\frac{z^{k}}{2\sqrt{z(z-1)(z-\lambda)}}dz.$$
Au voisinage du point $z=0$, on choisit comme param\`{e}tre local
$t=\sqrt{z}$. D'o\`{u}
\begin{eqnarray}
\omega&=&\frac{t^{2k}}{\sqrt{(t^2-1)(t^2-\lambda)}}dt,\nonumber\\
&=&at^{2k}(1+o(t^2))dt,\quad a\equiv\mbox{constante},\nonumber
\end{eqnarray}
ce qui montre que $\omega$ est holomorphe pour $k\geq 0$. Au
voisinage de $z=1$ et de $z=\lambda$, on utilise la m\^{e}me
m\'{e}thode et on obtient le m\^{e}me conclusion. Au voisinage de
$z=\infty$, on choisit comme param\`{e}tre local
$t=\frac{1}{\sqrt{z}}$ et on obtient
\begin{eqnarray}
\omega&=&-\frac{dt}{\sqrt{t^{2k}(1-t^2)(1-\lambda t^2)}}dt,\nonumber\\
&=&bt^{-2k}(1+o(t^2))dt,\quad b\equiv\mbox{constante},\nonumber
\end{eqnarray}
ce qui montre que $\omega$ est holomorphe pour $k\leq 0$. Donc au
voisinage de $0, 1, \lambda, \infty$, la forme diff\'{e}rentielle
$\omega$ est holomorphe si $k=0$, i.e., si
$$\omega =\frac{dz}{\sqrt{z\left( z-1\right) \left( z-\lambda \right)
}}.$$ Elle l'est aussi en dehors de ces points car si $c\in
\mathbb{C}$, $z_0\neq 0, 1, \lambda, \infty$, on choisit comme
param\`{e}tre local $t=z-z_0$, d'o\`{u}
\begin{eqnarray}
\omega&=&\frac{dt}{2\sqrt{(t+z_0)(t+z_0-1)(t+z_0-\lambda)}}dt,\nonumber\\
&=&c(1+o(t^3))dt,\quad c\equiv\mbox{constante},\nonumber
\end{eqnarray}
est holomorphe. En conclusion,
$$\omega =\frac{dz}{\sqrt{z\left( z-1\right) \left( z-\lambda \right)
}},$$ est l'unique diff\'{e}rentielle holomorphe sur la courbe
elliptique $X$.

\end{Exmp}
\begin{Exmp}
Les diff\'{e}rentielles
$$\omega _{k}=\frac{z^{k-1}dz}{\sqrt{\prod_{j=1}^{2g+1}\left( z-z_{j}\right) }},\text{ }k=1,2,\ldots
,g,$$ forment une base de diff\'{e}rentielles holomorphes sur la
courbe hyperelliptique $X$ de genre $g$ associ\'{e}e \`{a}
l'\'{e}quation
$$F(w,z) =w^{2}-\prod_{j=1}^{2g+1}\left( z-z_{j}\right)
=0.$$  En effet, les points de branchements de $X$ sont :
$z_1,...,z_{2g+1}, \infty$. D'apr\`{e}s le th\'{e}or\`{e}me 3.7,
le nombre de diff\'{e}rentielles holomorphes sur $X$ est \'{e}gal
au genre $g$ de $X$. Il suffit d'utiliser un raisonnement
similaire \`{a} celui de l'exemple pr\'{e}c\'{e}dent. Au voisinage
de $z=z_i$, $1\leq i \leq 2g+1$, on choisit comme param\`{e}tre
local $t=\sqrt{z-z_i}$ et on obtient
$$\frac{2\left( z_{i}+t^{2}\right)
^{k-1}tdt}{\sqrt{\prod_{j=1}^{2g+1}\left( t^{2}+z_{i}-z_{j}\right)
}},$$ qui sont holomorphes pour $k\geq 1$. De m\^{e}me, au
voisinage de $z=\infty$, on choisit comme param\`{e}tre local
$t=\frac{1}{\sqrt{z}}$ et on obtient
$$\omega=-\frac{2t^{2g+1}dt}{t^{2\left( k-g\right) }\sqrt{\prod_{j=1}^{2g+1}\left( 1-z_{j}t^{2}\right)
}},$$ qui sont holomorphes pour $k\leq g$. En dehors des points
$z_1,...,z_{2g+1}, \infty$ les diff\'{e}rentielles en questions
sont \'{e}videmment holomorphes.
\end{Exmp}

\section{Relations bilin\'{e}aires de Riemann}

\begin{thm}
Soient $\omega $\ et $\omega ^{\prime }$\ et deux
diff\'{e}rentielles holomorphes sur $X$. Alors,
$$\sum_{k=1}^{g}\left( \int_{a_{k}}\omega \int_{b_{k}}\omega ^{\prime }
-\int_{b_{k}}\omega \int_{a_{k}}\omega ^{\prime }\right) =0.$$ Si
en outre $\omega $\ est non nulle, alors
$$i\sum_{k=1}^{g}\left( \int_{a_{k}}\omega
\int_{b_{k}}\overline{\omega }-\int_{b_{k}}\omega \int_{a_{k}}
\overline{\omega }\right) \text{ }>  0.$$ (Ces deux expressions
s'appellent relations bilin\'{e}aires de Riemann).
\end{thm}
\emph{D\'{e}monstration}: Soit $X^*$ la repr\'{e}sentation normale
de $X$. Posons
$$f(p)=\int_{p_0}^p\omega,\qquad g(p)=\int_{p_0}^p\omega',$$
o\`{u} $p_0$ est un point fix\'{e} n'appartenant ni à $a_k$, ni
\`{a} $b_k$. Les fonctions $f$ et $g$ sont bien d\'{e}finies sur
$X^*$. Si $p\in a_k$, alors il est identifi\'{e} à $p^*\in
a_k^{-1}$. D'o\`{u}
\begin{equation}\label{eqn:euler}
f(p^*)=\int_{p_0}^{p^*}\omega=f(p)+\int_{b_k}\omega.
\end{equation}
Si $p\in b_k$, alors il est identifi\'{e} à $p^*\in b_k^{-1}$.
D'o\`{u}
\begin{equation}\label{eqn:euler}
f(p^*)=\int_{p_0}^{p^*}\omega=f(p)-\int_{a_k}\omega.
\end{equation}
Puisque la forme $\omega$ est exacte sur $X^*$, i.e., $\omega=df$,
alors
$$\omega\wedge \omega'=df\wedge \omega'=d(f\omega').$$
En appliquant la formule de Stokes-Cartan, on obtient
$$\int_X\omega\wedge \omega'=\int_{\partial X^*}f\omega'.$$
En tenant compte de (4.1) et (4.2), on obtient
\begin{eqnarray}
&&\int_{\partial
X^*}f\omega'\nonumber\\
&=&\sum_{k=1}^g\left[\int_{\partial a_k}f\omega' +\int_{\partial
b_k}f\omega'+\int_{\partial
a_k^{-1}}\left(f+\int_{b_k}\omega\right)\omega' +\int_{\partial
b_k^{-1}}\left(f-\int_{a_k}\omega\right)\omega'\right].\nonumber
\end{eqnarray}
Or
$$\int_{\partial
a_k^{-1}}\left(f+\int_{b_k}\omega\right)\omega'=-\int_{\partial
a_k}f\omega'-\int_{\partial
a_k}\left(\int_{b_k}\omega\right)\omega',$$
$$\int_{\partial
b_k^{-1}}\left(f-\int_{a_k}\omega\right)\omega'=-\int_{\partial
b_k}f\omega'+\int_{\partial
b_k}\left(\int_{a_k}\omega\right)\omega',$$ d'o\`{u}
$$\int_{\partial X^*}f\omega'=\sum_{k=1}^{g}\left( \int_{a_{k}}\omega \int_{b_{k}}\omega ^{\prime }
-\int_{b_{k}}\omega \int_{a_{k}}\omega ^{\prime }\right),$$ i.e.,
\begin{equation}\label{eqn:euler}
\int_X\omega\wedge \omega'=\sum_{k=1}^{g}\left( \int_{a_{k}}\omega
\int_{b_{k}}\omega ^{\prime } -\int_{b_{k}}\omega
\int_{a_{k}}\omega ^{\prime }\right).
\end{equation}
Si $\omega$, $\omega'$ sont deux formes diff\'{e}rentielles
holomorphes sur $X$ avec $\omega=f(z)dz$ et $\omega'=g(z)dz$,
alors
$$\omega\wedge \omega'=fgdz\wedge dz=0,$$
et par cons\'{e}quent la relation (4.3) implique que :
$$\sum_{k=1}^{g}\left( \int_{a_{k}}\omega \int_{b_{k}}\omega ^{\prime }
-\int_{b_{k}}\omega \int_{a_{k}}\omega ^{\prime }\right) =0.$$ Si
$\omega$ est une forme diff\'{e}rentielle holomorphe sur $X$ avec
$\omega\neq 0$ et $\omega=f(z)dz$, $z=x+iy$ alors
$$\omega\wedge \overline{\omega}=-2i|f|^2dx\wedge dy.$$
D\`{e}s lors
\begin{eqnarray}
i\sum_{k=1}^{g}\left( \int_{a_{k}}\omega
\int_{b_{k}}\overline{\omega }-\int_{b_{k}}\omega \int_{a_{k}}
\overline{\omega }\right)&=&i\int_X\omega\wedge
\overline{\omega},\mbox{ d'apr\`{e}s (4.3) avec }
\omega'=\overline{\omega},\nonumber\\
&=&\int_X|f|^2dx\wedge dy > 0, \quad \omega\neq 0,\nonumber
\end{eqnarray}
et le th\'{e}or\`{e}me est d\'{e}montr\'{e}. $\square$

Soit $\left( \omega _{1},\ldots ,\omega _{g}\right) $ une base de
diff\'{e}rentielles holomorphes sur la surface de Riemann $X$ de
genre $g$.
\begin{defn}
La matrice des p\'{e}riodes de $X$ est d\'{e}finie par $\Omega
=\left( E,F\right) ,$ o\`{u}
$${E}=\left(\begin{array}{ccc}
\int_{a_{1}}\omega _{1}&\cdots &\int_{a_{g}}\omega _{1}\\
\vdots &\ddots & \vdots \\
\int_{a_{1}}\omega _{g}&\cdots &\int_{a_{g}}\omega _{g}
\end{array}\right),\qquad
{F}=\left(\begin{array}{ccc}
\int_{b_{1}}\omega _{1}&\cdots &\int_{b_{g}}\omega _{1}\\
\vdots &\ddots & \vdots \\
\int_{b_{1}}\omega _{g}&\cdots &\int_{b_{g}}\omega _{g}
\end{array}\right).
$$
\end{defn}

\begin{prop}
Les 2g vecteurs $$\left(\begin{array}{c}
\int_{a_{1}}\omega _{1}\\
\vdots \\
\int_{a_{1}}\omega _{g}
\end{array}\right),...,\left(\begin{array}{c}
\int_{a_{g}}\omega _{1}\\
\vdots \\
\int_{a_{g}}\omega _{g}\end{array}\right),\left(\begin{array}{c}
\int_{b_{1}}\omega _{1}\\
\vdots \\
\int_{b_{1}}\omega
_{g}\end{array}\right),...,\left(\begin{array}{c}
\int_{b_{g}}\omega _{1}\\
\vdots \\
\int_{b_{g}}\omega _{g}
\end{array}\right),$$
sont $\mathbb{R}$-lin\'{e}airement ind\'{e}pendants.
\end{prop}
\emph{D\'{e}monstration}: En effet, proc\`{e}dons par l'absurde en
supposant que ces vecteurs soient $\mathbb{R}$-d\'{e}pendants.
Autrement dit, soit $k_1,...,k_g, m_1,...m_g$ des scalaires, tous
non nuls, tels que :
\begin{equation}\label{eqn:euler}
\sum_{j=1}^g\left[k_j\left(\begin{array}{c}
\int_{a_{j}}\omega _{1}\\
\vdots \\
\int_{a_{j}}\omega _{g}
\end{array}\right)+m_j\left(\begin{array}{c}
\int_{b_{j}}\omega _{1}\\
\vdots \\
\int_{b_{j}}\omega _{g}
\end{array}\right)\right]=0.
\end{equation}
On a aussi
\begin{equation}\label{eqn:euler}
\sum_{j=1}^g\left[k_j\left(\begin{array}{c}
\int_{a_{j}}\overline{\omega} _{1}\\
\vdots \\
\int_{a_{j}}\overline{\omega} _{g}
\end{array}\right)+m_j\left(\begin{array}{c}
\int_{b_{j}}\overline{\omega} _{1}\\
\vdots \\
\int_{b_{j}}\overline{\omega} _{g}
\end{array}\right)\right]=0,
\end{equation}
o\`{u} $\overline{\omega}_j$, $1\leq j\leq g$, d\'{e}signe le
conjugu\'{e} complexe de $\omega_j$. Consid\'{e}rons la matrice
$$
\Omega^*=\left(\begin{array}{cccccc}
\int_{a_{1}}\omega _{1}&\cdots &\int_{a_{g}}\omega _{1}&\int_{b_{1}}\omega _{1}&\cdots &\int_{b_{g}}\omega _{1}\\
\vdots &\ddots & \vdots&\vdots &\ddots & \vdots \\
\int_{a_{1}}\omega _{g}&\cdots &\int_{a_{g}}\omega
_{g}&\int_{b_{1}}\omega _{g}&\cdots &\int_{b_{g}}\omega _{g}\\
\int_{a_{1}}\overline{\omega} _{1}&\cdots &\int_{a_{g}}\overline{\omega} _{1}&\int_{b_{1}}\overline{\omega} _{1}&\cdots &\int_{b_{g}}\overline{\omega} _{1}\\
\vdots &\ddots & \vdots&\vdots &\ddots & \vdots \\
\int_{a_{1}}\overline{\omega} _{g}&\cdots
&\int_{a_{g}}\overline{\omega} _{g}&\int_{b_{1}}\overline{\omega}
_{g}&\cdots &\int_{b_{g}}\overline{\omega} _{g}
\end{array}\right).
$$
D'apr\`{e}s (4.4) et (4.5), on d\'{e}duit que $\mbox{rang
}\Omega^*<2g$ et d\`{e}s lors
\begin{eqnarray}
\int_{a_k}\sum_{j=1}^g\left(\alpha_j\omega_j+\beta_j\overline{\omega}_j\right)&=&0,\quad 1\leq k\leq g\nonumber\\
\int_{b_k}\sum_{j=1}^g\left(\alpha_j\omega_j+\beta_j\overline{\omega}_j\right)&=&0,\quad
1\leq k\leq g\nonumber
\end{eqnarray}
o\`{u} $\alpha_1,...,\alpha_g,\beta_1,...,\beta_g$ sont des
nombres complexes non tous nuls. En notant $$\omega\equiv
\sum_{j=1}^g\alpha_j\omega_j,$$ et $$\eta\equiv
\sum_{j=1}^g\overline{\beta}_j\omega_j,$$ les expressions
ci-dessus peuvent encore s'\'{e}crire sous la forme
\begin{equation}\label{eqn:euler}
\left\{\begin{array}{rl}
\int_{a_k}(\omega+\overline{\eta})=0,&\quad 1\leq k\leq g\\
\int_{b_k}(\omega+\overline{\eta})=0,&\quad 1\leq k\leq g
\end{array}\right.
\end{equation}
Montrons que la forme diff\'{e}rentielle $\omega+\overline{\eta}$
est exacte. Soient $\gamma$ un chemin d'extr\'{e}mit\'{e}s $p_0$
et $p$ et contenu dans $X$; le point $p_0$ \'{e}tant fix\'{e}
tandis que $p$ est arbitraire. On pose $$h(p)=\int_\gamma
(\omega+\overline{\eta}).$$ Cette int\'{e}grale ne d\'{e}pend pas
des extr\'{e}mit\'{e}s $p_0$ et $p$ de $\gamma$. En effet,
d\'{e}signons par $\gamma_1$ et $\gamma_2$ deux chemins joignant
les points $p_0$ et $p$. Notons que $\gamma=\gamma_1\cup \gamma_2$
est un chemin ferm\'{e} dans $X$ et on peut \'{e}crire
$$\gamma=\sum_{j=1}^g(m_ja_j+n_jb_j)+\partial \triangle,\quad
(m_j, n_j \in \mathbb{Z}),$$ o\`{u} $\triangle$ est une cha\^{i}ne
dans $X$. D'apr\`{e}s les \'{e}quations (4.6) et le
th\'{e}or\`{e}me de Stokes-Cartan, on a
$$\int_{\partial
\triangle}(\omega+\overline{\eta})=\int\int_{\triangle}
d(\omega+\overline{\eta})=\int\int_{\triangle} 0=0.$$ D\`{e}s
lors,
$$\int_{\gamma_1}(\omega+\overline{\eta})=\int_{\gamma_2}(\omega+\overline{\eta}),$$
ce qui montre que $h(p)=\int_{\gamma}(\omega+\overline{\eta})$ ne
d\'{e}pend pas du choix de $\gamma$ et elle est bien d\'{e}finie.
Par cons\'{e}quent, $\omega+\overline{\eta}=dh$, i.e.,
$\omega+\overline{\eta}$ est exacte. Nous allons maintenant
montrer que $\omega=\eta=0$, ce qui contredit l'ind\'{e}pendance
de $\omega_1,...,\omega_g$. Posons $\omega=f(z)dz$ et
$\eta=g(z)dz$. On a
$$\omega \wedge\eta=fgdz\wedge dz=0,$$
et
$$\eta\wedge \overline{\eta}=|g(z)|^2dz\wedge
d\overline{z}=-2i|g(z)|^2dx\wedge dy,$$ où $z=x+iy$. Pour montrer
que $\eta=0$, on va utiliser un raisonnement par l'absurde.
Supposons donc que $\eta\neq 0$, d'o\`{u}
$$\frac{i}{2}\int\int_X\eta\wedge \overline{\eta}=\int\int_X
|g(z)|^2dx\wedge dy >0,\quad \eta\neq 0.$$ On a
\begin{eqnarray}
\eta\wedge \overline{\eta}&=&\eta\wedge \omega+\eta\wedge
\overline{\eta},\nonumber\\
&=&\eta\wedge (\omega+\overline{\eta},\nonumber\\
&=&\eta\wedge dh,\nonumber\\
&=&-dh\wedge \eta,\nonumber\\
&=&-d(h\eta)-h\wedge d\eta,\nonumber\\
&=&-d(h\eta)-h\wedge dg\wedge dz,\nonumber\\
&=&-d(h\eta)-h\wedge (\frac{\partial g}{\partial z}dz+
\frac{\partial g}{\partial \overline{z}}d\overline{z})\wedge dz,\nonumber\\
&=&-d(h\eta),\nonumber
\end{eqnarray}
car la forme $\eta$ \'{e}tant holomorphe, la fonction $g$ est
aussi holomorphe et d'apr\`{e}s les conditions de Cauchy-Riemann,
on a $\frac{\partial g}{\partial \overline{z}}=0$. D'o\`{u}
$$\int\int_X \eta \wedge \overline{\eta}=-\int\int_X d(h\eta)=0,$$
ce qui est absurde. Donc $\eta=0$. En utilisant un raisonnement
similaire, on montre qu'on a aussi $\omega=0$ et la
d\'{e}monstration est compl\`{e}te. $\square$

\begin{thm}
Les relations bilin\'{e}aires de Riemann  sont \'{e}quivalentes
respectivement aux relations suivantes :
$$\Omega Q\Omega ^{\top }=0,\qquad i\Omega
Q\overline{\Omega }^{\top }\text{ } > 0,$$ o\`{u}
${Q}=\left(\begin{array}{cc}
0 &I\\
-I &0
\end{array}\right)$ est une matrice d'ordre $2g$ (dite matrice d'intersection) et $I$ la matrice
unit\'{e} d'ordre $g$.
\end{thm}
\emph{D\'{e}monstration}: Soient $\omega$ et $\omega'$ deux
diff\'{e}rentielles holomorphes,
$$\Phi=\left(\int_{a_1}\omega ... \int_{a_g}\omega \int_{b_1}\omega ... \int_{b_g}\omega \right),$$
la matrice des p\'{e}riodes de $\omega$ et
$$\Phi'=\left(\int_{a_1}\omega' ... \int_{a_g}\omega' \int_{b_1}\omega' ... \int_{b_g}\omega' \right),$$
celle de $\omega'$. Evidemment, la premi\`{e}re relation
bilin\'{e}aire de Riemann s'\'{e}crit
$$\Phi Q\Phi'^\top=0.$$
Or
\begin{eqnarray}
\omega&=&\sum_{j=1}^g\alpha_j\omega_j,\nonumber\\
\omega'&=&\sum_{j=1}^g\alpha'_j\omega_j,\nonumber
\end{eqnarray}
donc
\begin{eqnarray}
\Phi&=&\left(\sum_{j=1}^g\alpha_j\int_{a_1}\omega_j ...
\sum_{j=1}^g\alpha_j\int_{a_g}\omega_j
\sum_{j=1}^g\alpha_j\int_{b_1}\omega_j ... \sum_{j=1}^g\alpha_j\int_{b_g}\omega_j \right),\nonumber\\
&=&\alpha \Omega,\nonumber
\end{eqnarray}
o\`{u} $\alpha\equiv (\alpha_1 ... \alpha_g)$ et
\begin{eqnarray}
\Phi'&=&\left(\sum_{j=1}^g\alpha'_j\int_{a_1}\omega_j ...
\sum_{j=1}^g\alpha'_j\int_{a_g}\omega_j
\sum_{j=1}^g\alpha'_j\int_{b_1}\omega_j ... \sum_{j=1}^g\alpha'_j\int_{b_g}\omega_j \right),\nonumber\\
&=&\alpha' \Omega,\nonumber
\end{eqnarray}
o\`{u} $\alpha'\equiv (\alpha'_1 ... \alpha'_g)$. Si
$$\Phi Q \Phi'^\top=0,$$ i.e.,
$$\alpha\Omega Q\Omega^\top\alpha^\top=0,$$ quel que soit $ \alpha ,\alpha'$,
alors $\Omega Q \Omega^\top=0$. R\'{e}ciproquement, si $\Omega Q
\Omega^\top=0$, alors $\alpha\Omega Q\Omega^\top\alpha^\top=0$.
Autrement dit, la premi\`{e}re relation bilin\'{e}aire de Riemann
est \'{e}quivalente \`{a}
$$\Omega Q \Omega^\top=0.$$
Soit
$$\overline{\Phi}=\left(\int_{a_1}\overline{\omega} ... \int_{a_g}\overline{\omega}
\int_{b_1}\overline{\omega} ... \int_{b_g}\overline{\omega}
\right),$$ la matrice des p\'{e}riodes de $\overline{\omega}$. On
montre que la seconde relation bilin\'{e}aire de Riemann
s'\'{e}crit sous la forme
$$i\Phi Q \overline{\Phi}^\top>0,$$
ou encore
$$i\alpha \Omega Q \overline{\Omega}^\top \overline{\alpha}^\top>0.$$
Cette relation est vraie pour tout $\alpha$, donc elle est
\'{e}quivalente \`{a} celle-ci
$$i \Omega Q \overline{\Omega}^\top >0,$$
ce qui ach\`{e}ve la d\'{e}monstration. $\square$

\begin{prop}
La matrice $i\Omega Q\overline{\Omega }^{\top }$ est hermitienne
et elle est d\'{e}finie positive. On a
\begin{eqnarray}
\Omega Q \Omega^\top&=&EF^\top-FE^\top,\nonumber\\
\Omega Q
\overline{\Omega}^\top&=&E\overline{F}^\top-F\overline{E}^\top.\nonumber
\end{eqnarray}
En outre, la matrice $E$ est inversible.
\end{prop}
\emph{D\'{e}monstration}: On sait que $i \Omega Q
\overline{\Omega}^\top >0$. Par ailleurs, on a
\begin{eqnarray}
\left(i \Omega Q \overline{\Omega}^\top \right)^\top
&=&Q \overline{\Omega} Q^\top \Omega^\top ,\nonumber\\
&=&-Q \overline{\Omega} Q \Omega^\top ,\nonumber\\
&=&\overline{\left(Q \Omega Q \overline{\Omega}^\top
\right)},\nonumber
\end{eqnarray}
ce qui montre que $i \Omega Q \overline{\Omega}^\top$ est une
matrice hermitienne. On a
\begin{eqnarray}
\Omega Q \Omega^\top&=&(E\quad F)\left(\begin{array}{cc} O&I\\
-I&O
\end{array}\right)
\left(\begin{array}{c} E^\top\\
F^\top
\end{array}\right),\nonumber\\
&=&(-F\quad E)\left(\begin{array}{c} E^\top\\
F^\top
\end{array}\right),\nonumber\\
&=&-FE^\top+EF^\top.\nonumber
\end{eqnarray}
De m\^{e}me, on a
\begin{eqnarray}
\Omega Q \overline{\Omega}^\top&=&(E\quad F)\left(\begin{array}{cc} O&I\\
-I&O
\end{array}\right)
\left(\begin{array}{c} \overline{E}^\top\\
\overline{F}^\top
\end{array}\right),\nonumber\\
&=&(-F\quad E)\left(\begin{array}{c} \overline{E}^\top\\
\overline{F}^\top
\end{array}\right),\nonumber\\
&=&-F\overline{E}^\top+E\overline{F}^\top.\nonumber
\end{eqnarray}
Soit $s$ la solution de l'\'{e}quation : $sE=0$. D'apr\`{e}s ce
qui pr\'{e}c\'{e}de, on a
\begin{eqnarray}
s\left(i \Omega Q \overline{\Omega}^\top \right)^\top
\overline{s}^\top
&=&i s\left(E\overline{F}^\top-F\overline{E}^\top\right)\overline{s}^\top ,\nonumber\\
&=&i \left((sE)(\overline{sF})^\top-(sF)(\overline{sE})^\top\right) ,\nonumber\\
&=&0,\nonumber
\end{eqnarray}
puisque $sE=0$. D\`{e}s lors l'\'{e}quation $sE=0$ admet la seule
solution triviale $s=0$ et par cons\'{e}quent $\det E\neq 0$. Donc
$E$ est inversible et ach\`{e}ve la d\'{e}monstration. $\square$

Soit $(\omega_1,...,\omega_g)$ une base de diff\'{e}rentielles
holomorphes sur une surface de Riemann compacte $X$ de genre $g$
et soit $(a_1,...,a_g,b_1,...,b_g)$ une base symplectique du
groupe d'homologie $H_1(X,\mathbb{Z})$. Soit $\Omega=(E, F)$ la
matrice des p\'{e}riodes de $X$ associ\'{e}e \`{a} ces bases.
Notons que les p\'{e}riodes $\int_{a_j}\omega_j$ et
$\int_{b_j}\omega_j$, d\'{e}pendent non seulement de la structure
complexe de la surface $X$ mais d\'{e}pendent aussi de ces bases.
Nous allons voir qu'on n'est pas oblig\'{e} de maintenir ce
dernier choix. On montrera ci-dessous qu'il existe une unique base
(dite normalis\'{e}e) $(\eta_1,...,\eta_g)$ de $H^0(X,\Omega^1)$
telle que :
$$(\eta_1,...,\eta_g)^\top=E^{-1}(\omega_1,...,\omega_g)^\top,$$
ou ce qui revient au m\^{e}me, une unique base
$(\eta_1,...,\eta_g)$ telle que :
$$
\int_{a_l}\eta_k=\delta_{lk}\equiv\left\{\begin{array}{rl}
1&\mbox{si } k=l\\
0&\mbox{si } k\neq l
\end{array}\right.
$$

\begin{prop}
Soit $(\omega_1,...,\omega_g)$ une base de diff\'{e}rentielles
holomorphes sur une surface de Riemann de genre $g$. Alors, il
existe une nouvelle base $\left( \eta _{1},\ldots ,\eta
_{g}\right) :$
$$\eta _{k}=\sum_{j=1}^{g}c_{kj}\omega _{j},\quad1\leq k\leq g,$$
telle que :
$$
\int_{a_l}\eta_k=\delta_{lk},\quad 1\leq k,l\leq g
$$
c'est-\`{a}-dire de telle sorte que la matrice des p\'{e}riodes
associ\'{e}e \`{a} cette base soit de la forme $\left( I,\text{
}Z\right) $ o\`{u} $I$ est la matrice unit\'{e} et $Z=E^{-1}F.$ La
matrice $Z$ est sym\'{e}trique et \`{a} partie imaginaire
d\'{e}finie positive.
\end{prop}
\emph{D\'{e}monstration}: On a
$$\int_{b_l}\eta_k\sum_{j=1}^gc_{kj}\int_{b_l}\omega_j,$$
ce qui s'exprime en langage matriciel par $Z=CF$ o\`{u}
$C=(c_{kj})_{1\leq k,j\leq g}$ est la matrice des $c_{kj}$. Or
$\int_{a_l}\eta_k=\delta_{lk}$, donc
$$\sum_{j=1}^gc_{kj}\int_{a_l}\omega_j=\delta_{lk},$$
i.e., $CE=I$. D'apr\`{e}s la proposition pr\'{e}c\'{e}dente, la
matrice $E$ est inversible, d'o\`{u} $C=E^{-1}$ et par
cons\'{e}quent $Z=E^{-1}F$. Montrons tout d'abord que :
$Z=Z^\top$. En effet, on a
$$Z^\top-Z=F^\top(E^{-1})^\top-E^{-1}F,$$
d'o\`{u}
\begin{eqnarray}
E\left(Z^\top-Z\right)E^\top&=&E\left(F^\top(E^{-1})^\top-E^{-1}F\right)E^\top,\nonumber\\
&=&EF^\top-FE^\top,\nonumber\\
&=&\Omega Q \Omega^\top,\quad\mbox{proposition 4.5},\nonumber\\
&=&0,\nonumber
\end{eqnarray}
en vertu de la premi\`{e}re relation bilin\'{e}aire de Riemann.
Donc $Z^\top-Z=0$. Montrons maintenant que : $\mbox{Im} Z>0$. En
effet, d'apr\`{e}s la proposition 4.5, la matrice $i\Omega
J\overline{\Omega }^{\top }$ est d\'{e}finie positive. D'o\`{u}
\begin{eqnarray}
0&<&E^{-1}\left(i\Omega
Q\overline{\Omega }^{\top }\right)\left(\overline{E}^\top\right)^{-1},\nonumber\\
&=&iE^{-1}\left(E\overline{F}^\top-F\overline{E}^\top\right)\left(\overline{E}^\top\right)^{-1},
\quad\mbox{proposition 4.5},\nonumber\\
&=&i\left(\overline{F}^\top \left(\overline{E}^\top\right)^{-1}-E^{-1}F\right),\nonumber\\
&=&i\left(\overline{Z}^\top-Z\right),\nonumber\\
&=&i\left(\overline{Z}-Z\right),\quad\mbox{car Z est sym\'{e}trique},\nonumber\\
&=&\mbox{Im Z},\nonumber
\end{eqnarray}
la d\'{e}monstration s'ach\`{e}ve. $\square$

\section{Diviseurs}

Soient $p$ un point de $X$, $\tau_p:X\longrightarrow
\overline{\mathbb{C}}$ un param\`{e}tre local en $p$ (ou une
uniformisante locale en $p$, i.e., une carte locale en $p$
appliquant $p$ sur $0$) et $f$ une fonction m\'{e}romorphe au
voisinage de $p$. L'ordre de $f$ en $p$, est l'unique entier $n$
tel que : $f=\tau_p^n.g,$ o\`{u} $g$ est holomorphe ne s'annulant
pas en $p$. Dans le cas o\`{u} $f=0$, on choisit par convention
$n=+\infty$. L'entier $n$ d\'{e}pend de $p$ et de $f$ et on le
note $\mbox{ord}_p(f)$. On a
$$\mbox{ord}_p(f_1+f_2)\geq \inf (\mbox{ord}_p(f_1)
, \mbox{ord}_p (f_2)),$$ et
$$\mbox{ord}_p(f_1f_2)=\mbox{ord}_p(f_1)+\mbox{ord}_p(f_2).$$

\begin{defn}
Un diviseur sur une surface de Riemann $X$ est une combinaison
formelle du type\footnote{Il s'agit d'une notation utile qui
d\'{e}signe une collection finie de points sans ordre y compris
des entiers relatifs li\'{e}s \`{a} chaque point. Lorsqu'on
\'{e}crit par exemple : $2p_1+3p_2+p_3+...+p_n$, il faut bien
comprendre que cette \'{e}criture est purement formelle et qu'elle
n'a rien \`{a} voir avec $2$ coordonn\'{e}es $p_1$, $3$
coordonn\'{e}es $p_2$, etc...}
$$\mathcal{D}=\sum_{p\in X}n_p.p=\sum_j n_{j}p_{j},\quad n_j\in
\mathbb{Z},$$ avec $(p_{j})$ une famille localement finie de
points de $X$.
\end{defn}
La somme ci-dessus \'{e}tant finie (puisque $X$ est compacte), on
peut donc d\'{e}finir le support d'un diviseur comme \'{e}tant
l'ensemble fini de points $p_j$ pour lesquels le coefficient $n_j$
est non nul. Le diviseur $\mathcal{D}$ est fini si son support est
fini et ce sera toujours le cas si $X$ est une surface de Riemann
compacte. L'ensemble des diviseurs sur $X$ est un groupe
ab\'{e}lien not\'{e} $\mbox{Div} (X)$. L'addition des diviseurs
est d\'{e}finie par l'addition des coefficients.
\begin{defn}
Le degr\'{e}\footnote{ou ordre de $\mathcal{D}$, not\'{e}
$o(\mathcal{D})$.} du diviseur $\mathcal{D}$ est un entier
not\'{e} $\mbox{ deg } \mathcal{D}$ et est d\'{e}fini par $$\mbox{
deg } \mathcal{D}=\sum_jn_j.$$
\end{defn}
L'application $$\mbox{deg} : \mbox{ Div } (X)\longrightarrow
\mathbb{Z},\quad \mathcal{D}\longmapsto\mbox{ deg } \mathcal{D},$$
est un homomorphisme de groupe. Le noyau de cet homomorphisme est
l'ensemble des diviseurs de degr\'{e} $0$, not\'{e}
$\mbox{Div}^\circ (X)$, et forme un sous-groupe de $\mbox{Div}
(X)$.

Soit $f\neq 0$, une fonction m\'{e}romorphe sur $X$. A cette
fonction $f$, on fait correspondre un diviseur not\'{e} $(f)$ en
prenant pour $p_j$ les z\'{e}ros et les p\^{o}les de $f$ et pour
$n_j$ l'ordre de $p_j$ avec un signe n\'{e}gatif pour les
p\^{o}les. De mani\`{e}re plus pr\'{e}cise, on pose
$$(f)=\sum_{p\in X}\mbox{ord}_p(f) . p ,$$ o\`{u} les $\mbox{ord}_pf$
sont nuls sauf un nombre fini d'entre eux. En d\'{e}signant par
$\alpha_1,...,\alpha_l$ les z\'{e}ros de $f$ de multiplicit\'{e}
$n_1,...,n_l$ respectivement et par $\beta_1,...,\beta_m$ les
p\^{o}les de $f$ de multiplicit\'{e} $p_1,...,p_m$ respectivement,
on obtient
\begin{eqnarray}
(f)&=&\sum_{j=1}^l n_{j}\alpha_{j}-\sum_{j=1}^m p_{j}\beta_{j},\nonumber\\
&=&(f)_0-(f)_\infty,\nonumber
\end{eqnarray}
o\`{u} $(f)_0=(\text{diviseur des z\'{e}ros de }f)$ et
$(f)_\infty=(\text{diviseur des p\^{o}les de}f)$.
G\'{e}om\'{e}triquement, cel\`{a} signifie que $(f)_0$ correspond
\`{a} l'intersection de $X$ avec la courbe $f=0$ et $(f)_\infty$
\`{a} l'intersection de $X$ avec $\frac{1}{f}=0$.
\begin{Exmp}
Les diviseurs des fonctions
$$f_1(z)=z(z-1), \quad f_2(z)=\frac{1}{z}, \quad f_3(z)=z^2,$$
$$ f_4(z)=\frac{z}{(z-a)(z-b)}, \quad
f_5(z)=\frac{z^2}{(z-a)(z-b)}, (a,b\in \mathbb{C}),$$ sont
\'{e}videmment \'{e}gaux \`{a}
$$\left(f_1\right)=\{0\}-2\{\infty\}, \quad
\left(f_2\right)=\{0\}+\{\infty\}, \quad
\left(f_3\right)=2\{0\}-\{\infty\},$$
$$\left(f_4\right)=-\{a\}-\{b\}+\{0\}+\{\infty\}, \quad
\left(f_5\right)=-\{a\}-\{b\}+2\{0\}.$$
\end{Exmp}
On a
\begin{eqnarray}
(fg)&=&(f)+(g),\nonumber\\
(f^{-1})&=&-(f),\nonumber\\
(f)&=&(g)\Longrightarrow \frac{f}{g}=\mbox{constante}.\nonumber
\end{eqnarray}

Tout diviseur d'une fonction m\'{e}romorphe est dit diviseur
principal. L'ensemble des diviseurs principaux forme un
sous-groupe de $\mbox{Div}^\circ (X)$. Sur toute surface de
Riemann compacte, une fonction m\'{e}romorphe $f\neq 0$ a le
m\^{e}me nombre de z\'{e}ros que des p\^{o}les, donc $\mbox{ deg
}(f)=0$. Autrement dit, tout diviseur principal a le degr\'{e} $0$
mais en g\'{e}n\'{e}ral, les diviseurs de degr\'{e} $0$ ne sont
pas tous principaux.

Signalons deux notions qui seront \'{e}tudi\'{e}es dans les
sections suivantes : Le groupe de Picard, not\'{e} $\mbox{Pic}
(X)$, est le groupe des diviseurs quotient\'{e} par les diviseurs
principaux. La vari\'{e}t\'{e} jacobienne, not\'{e}e $\mbox{Jac}
(X)$, s'identifie au quotient du groupe des diviseurs de degr\'{e}
$0$ par les diviseurs principaux. Nous verrons que $\mbox{Jac}
(X)\simeq \mbox{Pic}^\circ(X)$ o\`{u} ce dernier d\'{e}signe le
sous-groupe de $\mbox{Pic} (X)$ form\'{e} par les classes des
diviseurs de degr\'{e} $0$.

On dit qu'un diviseur $\mathcal{D}$ est positif (ou effectif) et
on note $\mathcal{D}\geq 0$, si les entiers $n_{j}$ qui
interviennent dans la somme sont positifs. Plus
g\'{e}n\'{e}ralement, on d\'{e}finit la relation d'ordre partiel
$\geq$ sur les diviseurs par $\mathcal{D}_1 \geq \mathcal{D}_2$ si
et seulement si $\mathcal{D}_1-\mathcal{D}_2$ est positif. Deux
diviseurs $\mathcal{D}_1$ et $\mathcal{D}_2$ sont dits
lin\'{e}airement \'{e}quivalents (et on note $\mathcal{D}_1\sim
\mathcal{D}_2$) si $\mathcal{D}_1-\mathcal{D}_2$ est principal,
i.e., si $\mathcal{D}_1-\mathcal{D}_2=(f)$ o\`{u} $f$ est une
fonction m\'{e}romorphe.

Si $\mathcal{D}=\sum_{p\in X}n_p . p$ est un diviseur, on notera
$\mathcal{L}(\mathcal{D})$ l'ensemble des fonctions
m\'{e}romorphes $f$ telles que: $\mbox{ord}_p(f)+n_p\geq 0$, pour
tout $p\in X$. Autrement dit,
$$\mathcal{L}(\mathcal{D})=\{f\mbox{ m\'{e}romorphe sur }X:(f)+\mathcal{D}\geq 0\},$$
i.e., l'espace vectoriel des fonctions de $X$ dont le diviseur est
plus grand que $-\mathcal{D}$. Si $(f)+\mathcal{D}$ n'est $\geq 0$
pour aucun $f,$ on posera $\mathcal{L}(\mathcal{D})=0$. Par
exemple, si le diviseur $\mathcal{D}$ est positif alors
$\mathcal{L}(\mathcal{D})$ est l'ensemble des fonctions
holomorphes en dehors de $\mathcal{D}$ et ayant au plus des
p\^{o}les simples le long de $\mathcal{D}$.

\begin{prop}
Si $\mathcal{D}_1\sim \mathcal{D}_2$, alors
$\mathcal{L}(\mathcal{D}_1)$ est isomorphe \`{a}
$\mathcal{L}(\mathcal{D}_2)$ et en outre $\mbox{ deg }
\mathcal{D}_1=\mbox{ deg } \mathcal{D}_2$.
\end{prop}
\emph{D\'{e}monstration}: Par hypoth\`{e}se $\mathcal{D}_1\sim
\mathcal{D}_2$, donc par d\'{e}finition
$\mathcal{D}_1-\mathcal{D}_2=(f)$ o\`{u} $f$ est une fonction
m\'{e}romorphe. Par ailleurs, pour tout $g\in
\mathcal{L}(\mathcal{D}_1)$, on a
$$(g)+\mathcal{D}_1\geq 0,$$ et d\`{e}s lors
\begin{eqnarray}
(fg)+\mathcal{D}_2&=&(f)+(g)+\mathcal{D}_2,\nonumber\\
&=&\mathcal{D}_1-\mathcal{D}_2+(g)+\mathcal{D}_2,\nonumber\\
&=&(g)+\mathcal{D}_1\geq 0.\nonumber
\end{eqnarray}
L'application $$\mathcal{L}(\mathcal{D}_1)\longrightarrow
\mathcal{L}(\mathcal{D}_2),\quad g\longmapsto fg,$$ est
lin\'{e}aire et admet comme r\'{e}ciproque
$$\mathcal{L}(\mathcal{D}_2)\longrightarrow
\mathcal{L}(\mathcal{D}_1), \quad g\longmapsto \frac{g}{f}.$$
D'o\`{u} $\mathcal{L}(\mathcal{D}_1)\cong
\mathcal{L}(\mathcal{D}_2)$. En ce qui concerne la derni\`{e}re
assertion, on a
$$\mbox{ deg }(f)=\mbox{ deg }\mathcal{D}_1-\mbox{ deg } \mathcal{D}_2,$$ et
le r\'{e}sultat d\'{e}coule du fait que tout diviseur principal a
le degr\'{e} $0$. $\square$

\begin{Exmp}
Consid\'{e}rons le cas o\`{u} $X=\mathbb{P}^1(\mathbb{C})$ est la
droite projective complexe et soit $\mathcal{D}=\{0\}+\{1\}$ un
diviseur sur $X$. Les fonctions qui appartiennent \`{a} l'espace
$\mathcal{L}(\mathcal{D})$ sont les constantes et les fonctions :
$\frac{1}{z}$, $\frac{1}{z-1}$, $\frac{1}{z(z-1)}$.
\end{Exmp}

On peut associer \`{a} chaque forme diff\'{e}rentielle $\omega$ un
diviseur not\'{e} $(\omega)$. Si $\omega=fd\tau_p$ avec $f$ une
fonction m\'{e}romorphe sur $X$ et $\tau_p$ un param\`{e}tre local
en $p\in X$, on d\'{e}finit l'ordre de $\omega$ en $p$ par
$\mbox{ord}_p(\omega)=\mbox{ord}_0(f)$ et le diviseur $(\omega)$
de $\omega$ par $$(\omega)=\sum_{p\in X}\mbox{ord}_p(\omega).p.$$

Si $\mathcal{D}=\sum_{p\in X}n_p . p$ est un diviseur, on
d\'{e}finit de fa\c{c}on analogue \`{a}
$\mathcal{L}(\mathcal{D})$, l'espace lin\'{e}aire
$\mathcal{I}(\mathcal{D})$ comme \'{e}tant l'ensemble des formes
diff\'{e}rentielles m\'{e}romorphes $\omega$ sur $X$ telles que :
$\mbox{ord}_p\omega+n_p\geq 0$, pour tout $p\in X$. Autrement dit,
$$\mathcal{I}(\mathcal{D})=\{\omega\text{ m\'{e}romorphe sur }X:(\omega)+\mathcal{D}\geq 0\},$$
i.e., l'ensemble des formes diff\'{e}rentielles m\'{e}romorphes
$\omega$ sur $X$ telles que : $(\omega)+\mathcal{D}\geq 0$. Si
$\mathcal{D}_1\sim \mathcal{D}_2$, alors
$\mathcal{I}(\mathcal{D}_1)$ est isomorphe \`{a}
$\mathcal{I}(\mathcal{D}_2)$, d'o\`{u} $\dim
\mathcal{I}(\mathcal{D}_1)=\dim \mathcal{I}(\mathcal{D}_2)$. Dans
le cas o\`{u} le diviseur $\mathcal{D}$ est n\'{e}gatif, alors
$(\omega)+\mathcal{D}\geq 0$ est l'ensemble des formes
diff\'{e}rentielles qui n'ont pas de p\^{o}les et qui ont des
z\'{e}ros au moins aux points de $\mathcal{D}$. Notons enfin qu'en
vertu du th\'{e}or\`{e}me des r\'{e}sidus, on a $\sum_{p\in
X}\mbox{R\'{e}s }(\omega)=0$, o\`{u} $\mbox{R\'{e}s }(\omega)$ est
le r\'{e}sidu en $p$ de $\omega$, i.e., le coefficient de
$\frac{1}{\tau_p}$ dans le d\'{e}veloppement de $f$ en s\'{e}rie
de Laurent.

\begin{defn}
On appelle diviseur canonique sur $X$ et l'on d\'{e}signe par $K$,
le diviseur $(\omega)$ d'une $1-$forme m\'{e}romorphe $\omega\neq
0$ sur $X$.
\end{defn}

\begin{prop}
Soit $\mathcal{D}$ un diviseur sur une surface de Riemann compacte
$X$ et $K$ un diviseur canonique sur $X$. Alors, l'application
\begin{equation}\label{eqn:euler}
\psi:\mathcal{L}(K-\mathcal{D})\longrightarrow
\mathcal{I}(-\mathcal{D}),\quad f\longmapsto f\omega,
\end{equation}
est un isomorphisme.
\end{prop}
\emph{D\'{e}monstration}: En effet, on a
\begin{eqnarray}
(f\omega)=(f)+(\omega)&=&(f)+K,\nonumber\\
&\geq&-(K-\mathcal{D})+K,\nonumber\\
&=&-(-\mathcal{D}),\nonumber
\end{eqnarray}
ce qui montre que l'application $\psi$ est bien d\'{e}finie. Cette
derni\`{e}re est injective, i.e., l'\'{e}quation $f\omega=g\omega$
entraine $f=g$. Montrons que $\psi$ est surjective. Soit $\eta \in
\mathcal{I}(-\mathcal{D})$, d'o\`{u} il existe une fonction
m\'{e}romorphe $h$ sur $X$ telle que : $h\omega=\eta$. D\`{e}s
lors,
\begin{eqnarray}
(h)+K&=&(h)+(\omega),\nonumber\\
&=&(h\omega),\nonumber\\
&=&(\eta)\geq-(-\mathcal{D}),\nonumber
\end{eqnarray}
d'o\`{u} $(h)\geq -(K-\mathcal{D})$, $h\in
\mathcal{L}(K-\mathcal{D})$ et par cons\'{e}quent $\psi$ est
surjective. $\square$

\section{Le th\'{e}or\`{e}me de Riemann-Roch}

Nous allons maintenant \'{e}tudier un des th\'{e}or\`{e}mes les
plus importants de la th\'{e}orie des surfaces de Riemann
compactes : le th\'{e}or\`{e}me de Riemann-Roch. Il permet de
d\'{e}finir le genre d'une surface de Riemann qui est un invariant
fondamental. Il s'agit d'un th\'{e}or\`{e}me d'existence efficace
qui permet, entre autres, de d\'{e}terminer le nombre de fonctions
m\'{e}romorphes lin\'{e}airement ind\'{e}pendantes ayant certaines
restrictions sur leurs p\^{o}les. A cause de l'importance de ce
th\'{e}or\`{e}me, nous allons donner une preuve
\'{e}l\'{e}mentaire constructive bien qu'un peu technique et nous
mentionnons aussi quelques cons\'{e}quences de ce
th\'{e}or\`{e}me.

\begin{thm}
(de Riemann-Roch) : Si $X$ une surface de Riemann compacte et
$\mathcal{D}$ est un diviseur sur $X$, alors
\begin{equation}\label{eqn:euler}
\dim \mathcal{L}(\mathcal{D})-\dim
\mathcal{L}(K-\mathcal{D})=\mbox{ deg } \mathcal{D}-g+1,
\end{equation}
o\`{u} $K$ est le diviseur canonique sur $X$ et $g$ est le genre
de $X$. Cette formule peut s'\'{e}crire sous la forme
\'{e}quivalente
\begin{equation}\label{eqn:euler}
\dim \mathcal{L}(\mathcal{D})-\dim
\mathcal{I}(-\mathcal{D})=\mbox{ deg } \mathcal{D}-g+1.
\end{equation}
\end{thm}
\emph{D\'{e}monstration}: L'\'{e}quivalence entre les formules
(6.1) et (6.2) r\'{e}sulte imm\'{e}diatement de l'isomorphisme
(5.1). La preuve du th\'{e}or\`{e}me est imm\'{e}diate dans le cas
o\`{u} $\mathcal{D}=0$ car $\mathcal{L}(0)$ est l'ensemble des
fonctions holomorphes sur $X$. Or toute fonction holomorphe sur
une surface de Riemann compacte est constante, donc
$\mathcal{L}(0)=\mathbb{C}$. En outre, on sait que la dimension de
l'espace des formes diff\'{e}rentielles holomorphes sur $X$ est le
genre $g$ de $X$, d'o\`{u} le r\'{e}sultat. La
preuve du th\'{e}or\`{e}me va se faire en plusieurs \'{e}tapes :\\
\underline{\'{E}tape 1} : Soit $\mathcal{D}$ un diviseur positif.
Autrement dit,
$$\mathcal{D}=\sum_{k=1}^mn_kp_k,\quad n_k\in \mathbb{N}^*,\quad
p_k\in X.$$  Soit $f\in \mathcal{L}(\mathcal{D})$, i.e., $f$ a au
plus un p\^{o}le d'ordre $n_k$ en $p_k$. Au voisinage de $p_k$, on
a
$$df=\left(\sum_{j=-n_k-1}^\infty c_j^k\tau^j\right)d\tau,$$
donc $df$ est m\'{e}romorphe. Plus pr\'{e}cis\'{e}ment, $df$ a un
p\^{o}le d'ordre $n_k+1$ en $p_k$. Comme $f$ est m\'{e}romorphe,
alors $df$ ne peut pas avoir de p\^{o}le simple et d\`{e}s lors
son r\'{e}sidu est nul, i.e., $c_{-1}^k=0.$ Soit $\left(
a_{1},\ldots ,a_{g},b_{1},\ldots ,b_{g}\right) $ une base de
cycles dans le groupe d'homologie $H_{1}\left( X,\mathbb{Z}\right)
$ de telle fa\c{c}on que les indices d'intersection de cycles deux
\`{a} deux s'\'{e}crivent :
$$\left( a_{i},a_{i}\right) =\left( b_{i},b_{i}\right)
=0,\quad\left( a_{i},b_{j}\right) =\delta _{ij},\quad1\leq i,j\leq
g.$$ Rappelons (analyse harmonique) que pour tout $n\geq 2$ et
pour tout $p\in X$, il existe une diff\'{e}rentielle
m\'{e}romorphe unique $\eta$ sur $X$ telle que : $\eta$ est
holomorphe sur $X\backslash p$ tandis qu'autour de $p$,
$$\eta=\left(\frac{1}{\tau^n}+\circ (\tau)\right)d\tau,$$ o\`{u} $\tau$ est
un param\`{e}tre local en $p$ choisi avec $p=0$ et en outre, on a
$\int_{a_i}\eta=0$. Posons
\begin{equation}\label{eqn:euler}
\eta=df-\sum_{k=1}^m\sum_{j=2}^{n_k}c_{j}^k\eta_k^j,
\end{equation}
o\`{u} $\eta_k^j$ sont des diff\'{e}rentielles m\'{e}romorphes
ayant un p\^{o}le d'ordre $j$ en $p_k$ et holomorphes sur
$X\backslash p_k$. D'o\`{u}
$$\int_{a_i}\eta=\int_{a_i}df-\sum_{k=1}^m\sum_{j=2}^{n_k}c_j^k\int_{a_i}\eta_k^j.$$
L'int\'{e}grale d'une diff\'{e}rentielle exacte le long d'un
chemin ferm\'{e} \'{e}tant nulle, donc $\int_{a_i}df=0$. La forme
$\eta$ \'{e}tant holomorphe, alors
$$\eta=c_1\omega_1+...+c_g\omega_g,$$ o\`{u} $(\omega_1,...,\omega_g)$
est une base de $\Omega(X)$ et d\`{e}s lors
$$\eta=c_1\int_{a_i}\omega_1+...+c_g\int_{a_i}\omega_g,\quad i=1,...,g$$
Puisque la matrice $$E=\left(\int_{a_i}\omega_j\right)_{1\leq
i,j\leq g},$$ est inversible, alors $c_1=...=c_g=0$, donc $\eta=0$
et d'apr\`{e}s (6.3), on a
$$df=\sum_{k=1}^m\sum_{j=2}^{n_k}c_j^k\eta_k^j.$$
Consid\'{e}rons l'application
$$\varphi:\mathcal{L}(\mathcal{D})\longrightarrow
V\equiv\left\{(c_j^k):\sum_{k=1}^m\sum_{j=2}^{n_k}c_j^k\int_{b_l}\eta_k^j=0\right\},\quad
f\longmapsto c_j^k.$$ Notons que
\begin{eqnarray}
\mbox{Ker } \varphi&=&\left\{f : \mbox{m\'{e}romorphe sur X et
n'ayant pas de
p\^{o}le}\right\},\nonumber\\
&=&\{f : f\mbox{ est une constante }\},\nonumber
\end{eqnarray}
d'o\`{u} $\dim \mbox{ Ker } \varphi=1$ et par cons\'{e}quent,
$$\dim \mathcal{L}(\mathcal{D})=\dim V+1.$$ Les espaces
$\frac{\mathcal{L}(\mathcal{D})}{\mathbb{C}}$ et $V$ sont
isomorphes et on a
\begin{eqnarray}
\dim\mathcal{L}(\mathcal{D})-1&=&\dim V,\nonumber\\
&=&\dim\left\{(c_j^k):\sum_{k=1}^m\sum_{j=2}^{n_k}c_j^k\int_{b_l}\eta_k^j=0\right\},\nonumber\\
&=&\mbox{ deg }\mathcal{D}-\mbox{ rang }\mathcal{M},\nonumber
\end{eqnarray}
o\`{u}
$$
{\mathcal{M}}=\left(\begin{array}{ccccccccc}
\int_{b_l}\eta_1^2&\int_{b_l}\eta_1^3&...&\int_{b_l}\eta_1^{n_1+1}&\int_{b_l}\eta_2^2&...&
\int_{b_l}\eta_2^{n_2+1}&...&\int_{b_l}\eta_m^{n_m+1}\\
\int_{b_2}\eta_1^2&\int_{b_2}\eta_1^3&...&\int_{b_2}\eta_1^{n_1+1}&\int_{b_2}\eta_2^2&...&
\int_{b_2}\eta_2^{n_2+1}&...&\int_{b_2}\eta_m^{n_m+1}\\
\vdots&\vdots&...&\vdots&\vdots&...&\vdots&...&\vdots\\
\int_{b_g}\eta_1^2&\int_{b_g}\eta_1^3&...&\int_{b_g}\eta_1^{n_1+1}&\int_{b_g}\eta_2^2&...&
\int_{b_g}\eta_2^{n_2+1}&...&\int_{b_g}\eta_m^{n_m+1}
\end{array}\right),
$$
est la matrice dont le nombre de lignes est $g$ et le nombre de
colonnes est $\mbox{deg } \mathcal{D}$. Notons que
\begin{eqnarray}
\mbox{rang } \mathcal{M}&=&\mbox{ Nombre de colonnes }-
\mbox{Nombre de relations entre ces colonnes },\nonumber\\
&=&\mbox{deg } \mathcal{D}-\dim V,\nonumber\\
&=&\mbox{deg } \mathcal{D}-\dim \mathcal{L}(\mathcal{D})+1.
\end{eqnarray}
Calculons maintenant le rang de $\mathcal{M}$ d'une autre
fa\c{c}on. Soit $(\omega_1,...,\omega_g)$ une base orthonorm\'{e}e
de $\Omega^1(X)$. Au voisinage de $p_k$, la forme $\omega_s$ admet
un d\'{e}veloppement en s\'{e}rie de Taylor,
$$\omega_s=\left(\sum_{j=0}^\infty \alpha_{sj}^k\tau^j\right)d\tau.$$
Posons $\varphi_s\equiv \int_0^z\omega_s$ et soit $X^*$ la
repr\'{e}sentation normale de la surface de Riemann $X$. Notons
que si $\tau\in a_j$, alors il est identifi\'{e} \`{a} $\tau^*\in
a_j^{-1}$, d'o\`{u}
$$\varphi_s(\tau^*)=\varphi_s(\tau)+\int_{b_j}\omega_s.$$ De m\^{e}me,
si $\tau\in b_j$, alors il est identifi\'{e} à $\tau^*\in
b_j^{-1}$ et
$$\varphi_s(\tau^*)=\varphi_s(\tau)+\int_{a_j}\omega_s.$$ On a
\begin{eqnarray}
&&\int_{\partial X^*}\varphi_s\eta_k^n\nonumber\\
&=&\sum_{j=1}^g\left(\int_{a_j}\varphi_s\eta_k^n+\int_{b_j}\varphi_s\eta_k^n
+\int_{a_j^{-1}}\left(\varphi_s+\int_{b_j}\omega_s\right)\eta_k^n
+\int_{b_j^{-1}}\left(\varphi_s-\int_{a_j}\omega_s\right)\eta_k^n\right),\nonumber\\
&=&\sum_{j=1}^g\left(-\int_{b_j}\omega_s\int_{a_j}\eta_k^n
+\int_{a_j}\omega_s\int_{b_j}\eta_k^n\right),\nonumber\\
&=&\sum_{j=1}^g\left(-\omega_s(b_j)\eta_k^n(a_j)+\omega_s(a_j)\eta_k^n(b_j)\right),\nonumber\\
&=&\sum_{j=1}^g\omega_s(a_j)\eta_k^n(b_j),\nonumber\\
&=&\eta_k^n(b_s).
\end{eqnarray}
Or
\begin{eqnarray}
\int_{\partial X^*}\varphi_s\eta_k^n&=& 2\pi i\sum_k
\mbox{ R\'{e}s }_{p_k}\left(\varphi_s\eta_k^n\right),\nonumber\\
&=&2\pi i \frac{\alpha_{s,n-2}^k}{n-1},\nonumber
\end{eqnarray}
donc d'apr\`{e}s (6.5), la matrice $\mathcal{M}$ a comme
coefficient
$$\int_{b_s}\eta_k^n=\eta_k^n(b_s)=2\pi i \frac{\alpha_{s,n-2}^k}{n-1}.$$
D\`{e}s lors $$\det \mathcal{M}=C\det \mathcal{N},$$ o\`{u} $$
C\equiv(2\pi i)\left(\pi i\right)...\left(\frac{2\pi
i}{n_1}\right)(2\pi i)...\left(\frac{2\pi
i}{n_2}\right)...\left(\frac{2\pi i}{n_m}\right),$$ est une
constante et
$$
{\mathcal{N}}= \left(\begin{array}{ccccccccc}
\alpha_{1,0}^1&\alpha_{1,1}^1&...&\alpha_{1,n_1-1}^1&\alpha_{1,0}^2&...
&\alpha_{1,n_2-1}^2&...&\alpha_{1,n_m-2}^m\\
\alpha_{2,0}^1&\alpha_{2,1}^1&...&\alpha_{2,n_1-1}^1&\alpha_{2,0}^2&...
&\alpha_{2,n_2-1}^2&...&\alpha_{2,n_m-2}^m\\
\vdots&\vdots&...&\vdots&\vdots&...&\vdots&...&\vdots\\
\alpha_{g,0}^1&\alpha_{g,1}^1&...&\alpha_{g,n_1-1}^1&\alpha_{g,0}^2&...
&\alpha_{g,n_2-1}^2&...&\alpha_{g,n_m-2}^m\\
\end{array}\right).
$$
Calculons maintenant la dimension de l'espace
$\mathcal{L}(K-\mathcal{D})$ ou ce qui revient au m\^{e}me de
l'espace $\mathcal{I}(-\mathcal{D})$, i.e., celui des formes
diff\'{e}rentielles m\'{e}romorphes $\omega$ qui s'annulent $n_k$
fois au point $p_k$. On a
$$\omega=\sum_{s=1}^gX_s\omega_s=
\sum_{s=1}^gX_s\left(\alpha_{s,0}^k+\alpha_{s,1}^k\tau+\alpha_{s,2}^k\tau^2+\cdots\right)d\tau.$$
Pour que $\omega$ s'annule $n_k$ fois au point $p_k$, il faut que
les $n_k$ premiers termes dans l'expression ci-dessus soient
nulles. D\`{e}s lors,
$$(X_1,...,X_g).\mathcal{N}=0,$$
tandis que la dimension de $\mathcal{L}(K-\mathcal{D})$ coincide
avec celle de l'ensemble de $(X_1,...,X_g)$ tel que :
$\omega=\sum_{s=1}^gX_s\omega_s$ s'annule $n_k$ fois au point
$p_k$, i.e.,
\begin{eqnarray}
\dim \mathcal{L}(K-\mathcal{D})&=&g-\mbox{rang }
\mathcal{N},\nonumber\\
&=&g-\mbox{rang }\mathcal{M}.\nonumber
\end{eqnarray}
D'o\`{u}
$$\mbox{rang } \mathcal{M}=g-\dim \mathcal{L}(K-\mathcal{D}),$$
et en tenant compte de (6.4), on obtient finalement
$$\dim \mathcal{L}(\mathcal{D})-\dim \mathcal{L}(K-\mathcal{D})=\mbox{deg }
\mathcal{D}-g+1.$$ \underline{\'{E}tape 2} : La preuve donn\'{e}e
dans l'\'{e}tape 1 est valable pour tout diviseur lin\'{e}airement
\'{e}quivalent \`{a} un diviseur positif \'{e}tant donn\'{e} que
$\dim \mathcal{L}(\mathcal{D}),$ $\dim \mathcal{L}(K-\mathcal{D})$
(ou $\dim \mathcal{I}(-\mathcal{D})$)
et $\mbox{ deg } \mathcal{D}$ ne seront pas affect\'{e}s.\\
\underline{\'{E}tape 3} : Soit $f$ une fonction m\'{e}romorphe,
$\mathcal{D}$ un diviseur positif et posons
$$\mathcal{D}'=(f)+\mathcal{D}_0,$$ autrement dit, $\mathcal{D}'$ et
$\mathcal{D}_0$ sont lin\'{e}airement \'{e}quivalents. Nous avons
les assertions suivantes :
\begin{eqnarray}
&\textbf{(i)}&\dim \mathcal{L}(\mathcal{D}')=\dim
\mathcal{L}(\mathcal{D}_0).\nonumber\\
&\textbf{(ii)}&\dim \mathcal{L}(K-\mathcal{D'})=\dim
\mathcal{L}(K-\mathcal{D}_0).\nonumber\\
&\textbf{(iii)}&\mbox{ deg } \mathcal{D}'=\mbox{deg }
\mathcal{D}_0.\nonumber
\end{eqnarray}
qui d\'{e}coulent imm\'{e}diatement de la proposition 5.3.
Envisageons maintenant les diff\'{e}rents cas possibles :\\
$1^{\mbox{re}} cas$ : $\dim\mathcal{L}(\mathcal{D})>0$. Soit
$f_0\in \mathcal{L}(\mathcal{D})$, d'o\`{u} $(f_0)+\mathcal{D}>0,$
et
$$\dim\mathcal{L}((f_0)+\mathcal{D})-\dim\mathcal{L}\left(K-(f_0)-\mathcal{D}\right)
=\mbox{deg } ((f_0)+\mathcal{D})-g+1,$$ i.e.,
$$\dim\mathcal{L}(\mathcal{D})-\dim\mathcal{L}(K-\mathcal{D})
=\mbox{deg } \mathcal{D}-g+1.$$ $2^{\mbox{\`{e}me}} cas$ :
$\dim\mathcal{L}(\mathcal{D})=0$ et
$\dim\mathcal{L}(K-\mathcal{D})\neq 0$. En appliquant la formule
ci-dessus à $K-\mathcal{D}$, on obtient
\begin{equation}\label{eqn:euler}
\dim\mathcal{L}(K-\mathcal{D})-\dim\mathcal{L}(\mathcal{D})
=\mbox{deg} (K-\mathcal{D})-g+1.
\end{equation}
Pour la suite, on aura besoin du r\'{e}sultat int\'{e}ressant
suivant : Pour tout diviseur canonique $K$ sur une surface de
Riemann compacte $X$, on a
\begin{equation}\label{eqn:euler}
\mbox{deg } K=2g-2.
\end{equation}
o\`{u} $g$ est le genre de $X$. En effet, en posant
$\mathcal{D}=K$ dans la formule (6.1), on obtient $$\dim
\mathcal{L}(\mathcal{D})-\dim \mathcal{L}(0)=\mbox{deg
}\mathcal{D}-g+1.$$ Or $\mathcal{L}(0)=\mathbb{C}$, donc $\dim
\mathcal{L}(0)=1$ et on a
$$\mbox{ deg } K=g+\dim\mathcal{L}(K)-2.$$ Par ailleurs, en posant
$\mathcal{D}=0$ dans la formule (6.1), on obtient $$\dim
\mathcal{L}(0)-\dim\mathcal{L}(K)=\mbox{deg } 0-g+1,$$ d'o\`{u} $
\dim \mathcal{L}(K)=g$ et par cons\'{e}quent $\mbox{deg } K=2g-2$.
Ceci ach\`{e}ve la preuve du r\'{e}sultat annonc\'{e}. Pour
terminer la preuve du $2^{\mbox{ème}} cas$, on utilise ce
r\'{e}sultat et la formule (6.6), on obtient
$$
\dim \mathcal{L}(K-\mathcal{D})-\dim
\mathcal{L}(\mathcal{D})=-\mbox{deg }\mathcal{D}+g-1.
$$
$3^{\mbox{ème}} cas$ :
$\dim\mathcal{L}(\mathcal{D})=\dim\mathcal{L}(K-\mathcal{D})=0$.
Pour ce cas, on doit montrer que : $\mbox{deg }\mathcal{D}=g-1$.
Pour cel\`{a}, consid\'{e}rons deux diviseurs positifs
$\mathcal{D}_1$ et $\mathcal{D}_2$ n'ayant aucun point en commun
et posons $\mathcal{D}\equiv \mathcal{D}_1-\mathcal{D}_2$. On a
$$\mbox{ deg }\mathcal{D}=\mbox{deg }\mathcal{D}_1-\mbox{deg }
\mathcal{D}_2,$$ et
\begin{eqnarray}
\dim\mathcal{L}(\mathcal{D})&\geq& \mbox{deg }
\mathcal{D}_1-g+1,\nonumber\\
&=&\mbox{deg } \mathcal{D}+\mbox{ deg }
\mathcal{D}_2-g+1,\nonumber
\end{eqnarray}
i.e., $$\mbox{deg }\mathcal{D}_2-\dim
\mathcal{L}(\mathcal{D}_1)\leq \mbox{deg }\mathcal{D}+g-1.$$ Or
$$\mbox{ deg }\mathcal{D}_2-\dim\mathcal{L}(\mathcal{D}_1)\geq 0,$$
car sinon il existe une fonction $f\in \mathcal{L}(\mathcal{D}_1)$
qui s'annule en tout point de $\mathcal{D}_2$, donc $\mbox{deg
}\mathcal{D}\leq g-1$. En appliquant le m\^{e}me raisonnement
\^{a} $K-\mathcal{D}$, on obtient $\mbox{ deg }(K-\mathcal{D})\leq
g-1$. Comme $\mbox{deg }K=2g-2$ (voir (6.7)), alors $\mbox{deg
}\mathcal{D}\geq g-1.$ Finalement, $\mbox{deg }\mathcal{D}=g-1$,
ce qui ach\`{e}ve la d\'{e}monstration du th\'{e}or\`{e}me.
$\square$

\begin{rem}
La formule (6.1) peut s'\'{e}crire sous la forme suivante :
$$\dim H^0(X,\mathcal{O}_\mathcal{D})-\dim H^1(X,\mathcal{O}_\mathcal{D})=\mbox{ deg }
\mathcal{D}-g+1.$$ En introduisant la caract\'{e}ristique
d'Euler-Poincar\'{e} :
$$
\chi(\mathcal{D})\equiv\dim \mathcal{L}(\mathcal{D})-\dim
\mathcal{I}(-\mathcal{D})=\dim H^0(X,\mathcal{O}_\mathcal{D})-\dim
H^1(X,\mathcal{O}_\mathcal{D}),
$$
pour un diviseur $\mathcal{D}$ sur une surface de Riemann $X$ de
genre $g$, le th\'{e}or\`{e}me de Riemann-Roch s'\'{e}crit
$$\chi(\mathcal{D})=\mbox{deg }\mathcal{D}-g+1.$$
\end{rem}

\begin{rem}
Le th\'{e}or\`{e}me de Riemann-Roch se d\'{e}montre rapidement si
on utilise des th\'{e}ories encore plus pouss\'{e}es : la
dualit\'{e} de Kodaira-Serre et autres techniques (voir par
exemple [9] ou [18]). En effet, notons tout d'abord que
$K-\mathcal{D}$ est le diviseur correspondant au fibr\'{e}
$K\otimes\mathcal{L}(\mathcal{D})^*$ o\`{u} $\mathcal{L}^*$ est le
dual du fibr\'{e} $\mathcal{L}$. Les dimensions des espaces
$H^0(X, K\otimes\mathcal{L}(\mathcal{D})^*)$ et $H^1(X,
\mathcal{L}(\mathcal{D}))$ sont \'{e}gales puisque $H^0(X,
K\otimes\mathcal{L}(\mathcal{D})^*)$ est dual de $H^1(X,
\mathcal{L}(\mathcal{D}))$ en vertu de la dualit\'{e} de
Kodaira-Serre.
\end{rem}

\begin{rem}
On sait que toute fonction holomorphe sur une surface de Riemann
compacte $X$ est constante. Une question se pose : Que se passe
t-il dans le cas des fonctions m\'{e}romorphes? La r\'{e}ponse
d\'{e}coule du th\'{e}or\`{e}me de Riemann-Roch. Plus
pr\'{e}cis\'{e}ment, si $p$ est un point quelconque de $X$, on
peut trouver une fonction m\'{e}romorphe non constante, holomorphe
sur $X\backslash\{p\}$ et ayant un p\^{o}le d'ordre inf\'{e}rieur
ou \'{e}gal \`{a} $g+1$ en $p$. De m\^{e}me, on montre qu'il
existe sur $X$ des formes diff\'{e}rentielles holomorphes non
nulles, qui s'annulent en au moins un point.
\end{rem}

\begin{rem}
Soient $(\omega_1,...,\omega_g)$ une base de $\Omega^1(X)$. Soit
$(U,\tau)$ une carte locale en $p\in X$ avec $\tau(p)=0$. Il
existe des fonctions $f_j$ holomorphes sur $U$ telles que :
$\omega_j=f_j(\tau)d\tau$. Le wronskien de $\omega_1,...,\omega_g$
est d\'{e}fini par le d\'{e}terminant
$$W_\tau(\omega_1,...,\omega_g)\equiv W(f_1,...,f_g)=\det
\left(f^{(k-1)}_j\right)_{1\leq j,k\leq g}.$$ On dit que $p$ est
un point de Weierstrass si $W_\tau(\omega_1,...,\omega_g)$
s'annule. Dans le cas o\`{u} $p$ est un point de Weierstrass alors
on peut trouver une fonction m\'{e}romorphe sur $X$ ayant un
p\^{o}le unique d'ordre inf\'{e}rieur ou \'{e}gal au genre $g$ au
point $p$. Une autre application du th\'{e}or\`{e}me de
Riemann-Roch, permet de montrer l'existence d'une suite de $g$
entiers : $1=n_1<n_2<...<2g,$ $g\geq 1$, pour lesquels il n'existe
aucune fonction holomorphe sur $X\setminus {p}$, $p\in X$, et
ayant un p\^{o}le en $p$ d'ordre exactement $n_j$. On montre que
$p$ est un point de Weierstrass si et seulement si la suite des
$n_j$ est distincte de $\{1,2,...,g\}$.
\end{rem}

\section{La formule de Riemann-Hurwitz}

Nous donnons une preuve analytique de l'importante formule de
Riemann-Hurwitz. Elle exprime le genre d'une surface de Riemann à
l'aide du nombre de ses points de ramifications et du nombre de
ses feuillets. Nous montrons que cette formule fournit un moyen
efficace pour d\'{e}terminer le genre d'une surface de Riemann
donn\'{e}e. Quelques exemples int\'{e}ressants seront
\'{e}tudi\'{e}s.

Soient $X$ et $Y$ deux surfaces de Riemann compactes connexes et
soit $f$ une application holomorphe non constante de $X$ dans $Y$.
Notons que $f$ est un rev\^{e}tement, i.e., un morphisme surjectif
fini. Pour tout point $p\in X$, il existe une carte $\varphi$
(resp. $\psi$) de $X$ (resp. $Y$) centr\'{e}e en $p$ (resp.
$f(p)$) telles que : $f_{\psi\circ\varphi}(\tau)=\tau^n$, o\`{u}
$n$ est un entier strictement positif. L'entier $n-1$ s'appelle
indice de ramification de $f$ au point $p$ et on le note $V_f(p)$.
Lorsque $V_f(p)$ est strictement positif, alors on dit que $p$ est
un point de ramification (ou de branchement) de $f$. Une condition
n\'{e}cessaire et suffisante pour que $p$ soit un point de
ramification de $f$ est que le rang de $f$ en $p$ soit nul.
L'image $J$ des points de ramifications de $f$ ainsi que que son
image r\'{e}ciproque $I$ sont ferm\'{e}s et discrets. La
restriction de $f$ à $X\setminus I$ est un rev\^{e}tement de
$Y\setminus J$ dont le nombre de feuillets est le degr\'{e} de
l'application $f$ et on a
$$m\equiv \sum_{p\in f^{-1}(q)}\left(V_f(p)+1\right),\quad \forall q\in Y.$$

\begin{thm}
(Formule de Riemann-Hurwitz). Soient $X$ et $Y$ deux surfaces de
Riemann compactes de genre $g(X)$ et $g(Y)$ respectivement. Soit
$f$ une application holomorphe non constante de $X$ dans $Y$.
Alors
$$g(X)=m\left(g(Y)-1\right)+1+\frac{V}{2},$$
o\`{u} $m$ est le degr\'{e} de $f$ et $V$ est la somme des indices
de ramification de $f$ aux diff\'{e}rents points de $X$.
\end{thm}
\emph{D\'{e}monstration}: Soit $f:X\longrightarrow Y$, une
application holomorphe non constante de degr\'{e} $m$. Soit
$\omega$ (resp. $\eta$) une forme diff\'{e}rentielle
m\'{e}romorphe non nulle sur $Y$ (resp. $X$). Soit $\tau$ (resp.
$\upsilon$) un param\`{e}tre local sur $X$ (resp. $Y$) et
supposons que : $\upsilon=f(\tau)$. En d\'{e}signant par
$\omega=h(v)dv$, la forme diff\'{e}rentielle m\'{e}romorphe sur
$Y$, alors forme diff\'{e}rentielle $\eta$ sur $X$ s'\'{e}crit en
terme de $\tau$ sous la forme, $\eta=h(f(\tau))f'(\tau)d\tau$.
Nous allons voir que cette derni\`{e}re est aussi m\'{e}romorphe.
Notons que si on remplace $\tau$ par $\tau_1$, avec
$\tau=w(\tau_1)$, alors en terme de $\tau_1$ l'application $f$
s'\'{e}crit $\upsilon=(f\circ w)(\tau_1)$, et donc nous attribuons
\`{a} $\tau_1$ l'expression
$h(f(w(\tau_1)))f'(w(\tau_1))w'(\tau_1)d\tau_1$ ce qui montre que
$\eta$ est une forme diff\'{e}rentielle m\'{e}romorphe. On peut
supposer que $\tau$ s'annule en $p\in X$ et que $\upsilon$
s'annule en $f(p)$. D\`{e}s lors, $\upsilon=\tau^{V_f(p)+1}$
o\`{u} $V_f(p)$ est l'indice de ramification de $f$ au point $p$.
Par cons\'{e}quent,
$$\mbox{ord}_p\eta=(V_f(p)+1)\mbox{ ord}_{f(p)}\omega+V_f(p),$$ et
$$\sum_{p\in X}\mbox{ ord}_p\eta=\sum_{p\in X}(V_f(p)+1)\mbox{ ord}_{f(p)}\omega+V,$$
o\`{u} $V=\sum_{p\in X}V_f(p)$. D'apr\`{e}s la formule 6.7, on a
$$\sum_{p\in X}\mbox{ ord}_p\eta=2g(X)-2,$$
et
\begin{eqnarray}
\sum_{p\in X}(V_f(p)+1)\mbox{ ord}_{f(p)}\omega&=&\sum_{p\in X,
V_f(p)=0}\mbox{ ord}_{f(p)}\omega,\nonumber\\
&=&\sum_{q\in Y}m . \mbox{ ord}_q\omega,\nonumber\\
&=&m(2g(Y)-2).\nonumber
\end{eqnarray}
Par cons\'{e}quent,
$$2g(X)-2=m(2g(Y)-2)+V,$$
ce qui ach\`{e}ve la preuve du th\'{e}or\`{e}me. $\square$

Une des cons\'{e}quences les plus int\'{e}ressantes de la formule
de Riemann-Hurwitz est de donner un moyen efficace de calculer le
genre d'une surface de Riemann donn\'{e}e.

\begin{Exmp}
Un cas particulier important est repr\'{e}sent\'{e} par les
courbes hyperelliptiques $X$ de genre $g(X)$ d'\'{e}quations
$$w^{2}=p_{n}(z)=\prod_{j=1}^n(z-z_j),$$
o\`{u} $p_{n}(z)$ est un polyn\^{o}me sans racines multiples,
i.e., tous les $z_j$ sont distincts. Notons que
$$f:X\longrightarrow Y=\mathbb{P}^1(\mathbb{C})=\mathbb{C}\cup \{\infty\},$$
est un rev\^{e}tement double ramifi\'{e} le long des points $z_j$.
Chaque $z_j$ est ramifi\'{e} d'indice $1$ et en outre le point
\`{a} l'infini $\infty$ est ramifi\'{e} si et seulement si $n$ est
impair. D'apr\`{e}s la formule de Riemann-Hurwitz, on a
\begin{eqnarray}
g(X)&=&m(g(Y)-1)+1+\frac{V}{2},\nonumber\\
&=&2(0-1)+1+\frac{1}{2}\sum_{p\in X}V_f(p),\nonumber\\
&=&E\left(\frac{n-1}{2}\right),\nonumber
\end{eqnarray}
o\`{u} $E\left(\frac{n-1}{2}\right)$ d\'{e}signe la partie
enti\`{e}re de $\left(\frac{n-1}{2}\right)$. Les courbes
hyperelliptiques de genre $g$ sont associ\'{e}es aux \'{e}quations
de la forme : $w^2=p_{2g+1}(z),$ ou $w^2=p_{2g+2}(z),$ (selon que
le point à l'infini $\infty$ est un point de branchement ou non)
avec $p_{2g+1}(z)$ et $p_{2g+2}(z)$ des polyn\^{o}mes sans racines
multiples.
\end{Exmp}

\begin{Exmp}
D\'{e}terminons le genre $g$ de la surface de Riemann $X$
associ\'{e}e \`{a} l'\'{e}quation :
$$F(w,z)=w^3+p_2(z)w^2+p_4(z)w+p_6(z)=0,$$
o\`{u} $p_j(z)$ d\'{e}signe un polyn\^{o}me de degr\'{e} $j$. On
proc\`{e}de comme suit : on a
\begin{eqnarray}
F(w,z)&=&w^3+az^2w^2+bz^4w+cz^6+\text{termes d'ordre
inf\'{e}rieur},\nonumber\\
&=&\prod_{j=1}^3\left(z+\alpha_jz^2\right)+\text{termes d'ordre
inf\'{e}rieur}.\nonumber
\end{eqnarray}
Consid\'{e}rons $F$ comme un rev\^{e}tement par rapport \`{a} $z $
et cherchons ce qui ce passe quand $z \nearrow \infty .$ On a
\begin{eqnarray}
(w)_\infty&=&-2P-2Q-2R,\nonumber\\
(z)_\infty&=&-P-Q-R.\nonumber
\end{eqnarray}
Posons $t=\frac{1}{z}$, d'o\`{u}
$$F(w,z)= \frac{1}{t^{6}}(t^6z^3+at^4z^2+bt^2z+c)+\cdots.$$
Ceci sugg\`{e}re le changement de cartes suivant :
$$(w,z)\longmapsto \left(\zeta=t^2w,t=\frac{1}{z }\right).$$ On a
\begin{eqnarray}
\frac{\partial F}{\partial w}&=&3w ^2+2p_2(z)w+p_4(z),\nonumber\\
&=&3w^2+2az^2w+bz^4+\cdots,\nonumber\\
&=&\frac{3\zeta^2}{t^4}+\frac{2a\zeta}{t^4}+\frac{b}{t^4}+\cdots\nonumber
\end{eqnarray}
La fonction $\frac{\partial F}{\partial w}$ \'{e}tant
m\'{e}romorphe sur la surface de Riemann $X$, alors Le nombre de
z\'{e}ros de cette fonction coincide avec celui de ses p\^{o}les.
Comme
\begin{eqnarray}
\left(\frac{\partial F}{\partial
w}\right)_{P}&=&-4P,\nonumber\\
\left(\frac{\partial F}{\partial
w}\right)_{Q}&=&-4Q,\nonumber\\
\left(\frac{\partial F}{\partial
w}\right)_{R}&=&-4R,\nonumber\\
\left(\frac{\partial F}{\partial
w}\right)_{\infty}&=&-4(P+Q+R),\nonumber
\end{eqnarray}
alors le nombre de z\'{e}ros de $\frac{\partial F}{\partial w }$
dans la partie affine $X\setminus \{P,Q,R\}$ est \'{e}gal à $8$,
et d'apr\`{e}s la formule de Riemann-Hurwitz, on a $g(X) =4.$
\end{Exmp}

\begin{Exmp}
Calculons le genre de la surface de Riemann $X$ associ\'{e}e au
polyn\^{o}me :
$$w^4=z^4-1.$$
Ici, on a quatre feuillets. Les points de ramifications \`{a}
distance finie sont $1,-1,i$ et $-i$. On note que $z=\infty$ n'est
pas un point de ramification. L'indice de ramification \'{e}tant
\'{e}gal \`{a} $12$, alors d'apr\`{e}s la formule de
Riemann-Hurwitz, le genre de la surface de Riemann en question est
\'{e}gal à $3$.
\end{Exmp}

\begin{Exmp}
Consid\'{e}rons la courbe de Fermat $X$ associ\'{e}e \`{a}
l'\'{e}quation :
$$w^n+z^n=1,\quad n\geq 2.$$
Ici on a un rev\^{e}tement de degr\'{e} $n$. Chaque racine
$n^{\mbox{\`{e}me}}$ de l'unit\'{e} est ramifi\'{e} d'indice $n-1$
tandis que le point \`{a} l'infini $\infty$ n'est pas un point de
ramification et par cons\'{e}quent $$g(X)=\frac{(n-1)(n-2)}{2}.$$
L'\'{e}quation de Fermat : $$U^n+V^n=W^n,$$  (avec
$w=\frac{U}{Z}$, $z=\frac{V}{Z}$) \'{e}tant de genre $\geq 1$ pour
$n\geq 3$, elle n'admet donc qu'un nombre fini de solutions. Ce
fut une des pistes utilis\'{e}es par A. Wiles pour prouver le
grand th\'{e}or\`{e}me de Fermat : pour $n\geq 3$ cette
\'{e}quation n'a pas de solution non triviale.
\end{Exmp}

\section{Le théorème d'Abel}

\begin{thm}
Soient $p_1,...,p_m, q_1,...q_m$ des points de $X$.
Alors les deux conditions suivantes sont \'{e}quivalentes:\\
(i) Il existe une fonction m\'{e}romorphe $f$ telle que:
$$(f) =\sum_{j=1}^{m}q_{j}-\sum_{j=1}^{m}p_{j}.$$
(ii) Il existe un chemin ferm\'{e} $\gamma $\ tel que: $$\forall
\omega \in \Omega^1(X) ,
\quad\sum_{j=1}^{m}\int_{p_{j}}^{q_{j}}\omega =\int_{\gamma
}\omega.$$
\end{thm}
\emph{D\'{e}monstration}: Montrons que : $(i)\Longrightarrow
(ii)$. Soit $f$ une fonction m\'{e}romorphe sur $X$ et soit
$\omega=d\log f\in \Omega^1(X)$. Posons $\varphi(p)\equiv
\int_{p_0}^p\omega$. Nous allons calculer $\int_{\partial
X^*}\varphi d\log f$ o\`{u} $X^*$ est la repr\'{e}sentation
normale de $X$. D'apr\`{e}s le th\'{e}or\`{e}me des r\'{e}sidus,
on a
$$\int_{\partial X^*}\varphi d\log f=2\pi \sum \mbox{Rés }(\varphi d\log
f).$$ Notons que puisque la fonction $\varphi$ est holomorphe,
elle n'a donc pas de p\^{o}les et par cons\'{e}quent pour calculer
le r\'{e}sidu de $\varphi d\log f$ sur $X^*$, il suffit de
d\'{e}terminer les p\^{o}les de $d\log f$. Par hypoth\`{e}se, la
fonction $f$ est m\'{e}romorphe. Donc on a au voisinage d'un
p\^{o}le $r$ d'ordre $m$,
$$f(z)=(z-r)^m g(z),\quad m>0,$$
et au voisinage d'un z\'{e}ro $r$ d'ordre $m$,
$$f(z)=(z-r)^m g(z),\quad m<0,$$
avec $g(z)$ une fonction holomorphe au voisinage de $r$ et telle
que : $g(r)\neq 0$. On a
$$\log f=m\log (z-r)+\log g(z),$$
et
$$d\log f=\frac{m}{z-r}+\frac{g'(z)}{g(z)}.$$
Notons que $d\log f$ a des p\^{o}les aux z\'{e}ros et aux
p\^{o}les de $f$. D\`{e}s lors,
$$\mbox{R\'{e}s }(\varphi d\log f)=\varphi(r_j).m_j,$$
o\`{u} $r_j$ est un p\^{o}le ou un z\'{e}ro de $f$ avec le signe
positif ou n\'{e}gatif suivant que $r_j$ est un z\'{e}ro ou un
p\^{o}le de $f$ tandis que les entiers $m_j$ d\'{e}signent la
multiplicit\'{e} de $r_j$. Donc
\begin{eqnarray}
\int_{\partial X^*}\varphi d\log f&=&2\pi \sum
\varphi(r_j).m_j,\nonumber\\
&=&2\pi \sum_{j=1}^m \int_{p_j}^{q_j} \omega,
\end{eqnarray}
en vertu de la d\'{e}finition de $\varphi$, $r_j$ et $m_j$.
Calculons cette int\'{e}grale d'une autre mani\`{e}re. En
raisonnant comme dans la preuve de la premi\`{e}re relation
bilin\'{e}aire de Riemann (th\'{e}or\`{e}me 4.1), on obtient
$$\int_{\partial X^*}\varphi d\log f=\sum_{m=1}^g\left(\omega(a_m)\int_{b_m}d\log f-\omega(b_m)\int_{a_m}d\log
f\right).$$ On a
\begin{eqnarray}
\int_{b_m}d\log f(z)&=&\int_{b_m}d\log |f(z)|+i\int_{b_m}d(\arg f(z)),\nonumber\\
&=&2\pi i \alpha_m,\quad \alpha_m\in \mathbb{Z},\nonumber
\end{eqnarray}
et de m\^{e}me
$$\int_{a_m}d\log f(z)=2\pi i \beta_m,\quad
\beta_m\in \mathbb{Z}.$$ D'o\`{u}
$$\int_{\partial X^*}\varphi d\log f=2\pi i\sum_{m=1}^g\left(\alpha_m\omega(a_m)-\beta_m \omega(b_m)\right).$$
Posons $\gamma=\displaystyle{\sum_{m=1}^g (\alpha_m a_m-\beta_m
b_m)}$, d'o\`{u}
$$\int_{\partial X^*}\varphi d\log f=2\pi i\int_\gamma \omega.$$
En comparant avec (8.1), on obtient
$$\int_\gamma \omega=\sum_{j=1}^m \int_{p_j}^{q_j}\omega.$$
Montrons maintenant que : $(ii)\Longrightarrow (i)$. Soient
$p_1,...,p_m, q_1,...q_m$ des points de $X$ et $\gamma$ un chemin
ferm\'{e} tel que:
$$
\forall\omega\in\Omega^1(X),\quad\sum_{j=1}^{k}\int_{p_{j}}^{q_{j}}\omega
=\int_{\gamma}\omega.
$$
Montrons qu'il existe une fonction m\'{e}romorphe telle que:
$$(f) =\sum_{j=1}^{m}q_{j}-\sum_{j=1}^{m}p_{j}\equiv \mathcal{D}.$$
Rappelons (analyse harmonique) que pour tout $p,q \in X$, il
existe une diff\'{e}rentielle m\'{e}romorphe $\eta$ sur $X$ ayant
des p\^{o}les simples en $p, q$ et telle que : $\mbox{Rés}_p \eta
=1$, $\mbox{Rés}_q \eta =-1$. On en d\'{e}duit que si
$p_1,...,p_n\in X$ et $c_1,...,c_n\in \mathbb{C}$ avec
$\sum_{j=1}^nc_j=0$, alors il existe une diff\'{e}rentielle
m\'{e}romorphe $\eta$ sur $X$ ayant des p\^{o}les simples en $p_j$
et telle que : $\mbox{R\'{e}s}{_p{_j}} \eta =c_j$ et
$\int_{a_j}\eta=0$. En effet, soient $p_1,...,p_n\in X$, $q\in X$,
$q\neq p_j$ $1\leq j\leq n$ et $c_1,...,c_n\in \mathbb{C}$ avec
$\sum_{j=1}^nc_j=0$. On peut trouver des diff\'{e}rentielles
m\'{e}romorphes $\eta_j$ sur $X$ ayant des p\^{o}les simples en
$p_j$, $q$ et telles que : $\mbox{R\'{e}s}{_p{_j}} \eta_j =1$ et
$\mbox{R\'{e}s}{_q} \eta_j =-1$. La forme diff\'{e}rentielle
$\lambda_1=\sum_{j=1}^nc_j\eta_j$ a des p\^{o}les simples en $p_j$
avec $\mbox{R\'{e}s}{_{p_j}}\lambda_1=c_j$ mais n'a pas de
p\^{o}les en $q$ avec
$\mbox{R\'{e}s}{_q}\lambda_1=(-1)\left(\sum_{j=1}^nc_j\right)=0$.
Soit $(\omega_1,...,\omega_g)$ une base de diff\'{e}rentielles
holomorphes sur $X$ et consid\'{e}rons la forme diff\'{e}rentielle
$$\lambda_2=\lambda_1+\sum_{k=1}^g\alpha_k\omega_k,$$ o\`{u}
$\alpha_1,...,\alpha_g$ sont des constantes \`{a} d\'{e}terminer.
On ajoute \`{a} $\lambda_1$ une forme diff\'{e}rentielle
holomorphe, on ne change donc rien \`{a}  ses p\^{o}les qui sont
simples ni \`{a} ses r\'{e}sidus. Il reste \`{a} montrer que :
$\int_{a_j}\lambda_2=0$ ou ce qui revient au m\^{e}me \`{a}
d\'{e}terminer les constantes $\alpha_1,...,\alpha_g$ telles que :
$$\int_{a_j}\lambda_1+\sum_{k=1}^g\alpha_k\int_{a_j}\omega_k=0.$$
Ceci revient \`{a} r\'{e}soudre le syst\`{e}me de $g$
\'{e}quations \`{a} $g$ inconnues suivant :
$$
\left(
\begin{array}{ccc}
\int_{a_1}\omega_1&\cdots&\int_{a_1}\omega_g\\
\vdots&\ddots&\vdots\\
\int_{a_g}\omega_1&\cdots&\int_{a_g}\omega_g
\end{array}
\right) \left(
\begin{array}{c}
\alpha_1\\
\vdots\\
\alpha_g
\end{array}
\right)= \left(
\begin{array}{c}
-\int_{a_1}\lambda_1\\
\vdots\\
-\int_{a_g}\lambda_1
\end{array}
\right).
$$
Notons que la matrice \`{a} gauche est la transpos\'{e}e de la
matrice $E$ intervenant dans la d\'{e}finition de la matrice des
p\'{e}riodes (d\'{e}finition 4.2). D'apr\`{e}s la proposition 4.5,
la matrice $E$ est inversible, donc sa matrice transpos\'{e}e
aussi et par cons\'{e}quent le syst\`{e}me ci-dessus admet une
solution pour laquelle $\lambda_2$ est le $\eta$ cherch\'{e}.
Revenons maintenant au diviseur $$\mathcal{D}\equiv(f)
=\sum_{j=1}^{m}q_{j}-\sum_{j=1}^{m}p_{j},$$ et notons que l'on
peut l'\'{e}crire sous la forme
$$\mathcal{D}=\sum_{j=1}^{n}c_{j}p_j, \quad n<m, \quad c_j\in \mathbb{Z};$$
il suffit de regrouper les $p_j$ et $q_j$ qui sont les m\^{e}mes.
La somme des coefficients n'a pas chang\'{e}; elle valait
$m.1+m(-1)=0$, donc $\sum_{j=1}^{n}c_{j}=0$. D'apr\`{e}s ce qui
pr\'{e}c\`{e}de, il existe une diff\'{e}rentielle m\'{e}romorphe
$\eta$ sur $X$ ayant des p\^{o}les simples aux points $p_j$ et
telle que : $\mbox{Rés}{_p{_j}} \eta =c_j$ et $\int_{a_j}=0$. Soit
$(\omega_1,...,\omega_g)$ une base orthonorm\'{e}e de
$\Omega^1(X)$ et posons $\varphi_j(p)=\int_{p_0}^p\omega_j$. En
raisonnant comme dans la preuve de la premi\`{e}re relation
bilin\'{e}aire de Riemann (th\'{e}or\`{e}me 4.1), on obtient
\begin{eqnarray}
\int_{\partial X^*}\varphi_j
\eta&=&\sum_{m=1}^g\left(\omega_j(a_m)\eta{b_m}-\omega_j(b_m)\eta{a_m}\right),\nonumber\\
&=&\sum_{m=1}^g\left(\int_{a_m}\omega_j\int_{b_m}\eta-\int_{b_m}\omega_j\int_{a_m}\eta\right),\nonumber\\
&=&\int_{b_j}\eta,\nonumber\\
&=&\eta(b_j),
\end{eqnarray}
car $\int_{a_m}\omega_j=\delta_{mj}$ et $\int_{a_m}\eta=0$. D'un
autre c\^{o}t\'{e}, on a
$$\int_{\partial X^*}\varphi_j
\eta=2\pi i\sum_{k=1}^n\mbox{R\'{e}s}_{p_k}(\varphi_k\eta).$$ En
utilisant un raisonnement similaire \`{a} celui fait
pr\'{e}c\'{e}demment pour montrer que $(i)\Rightarrow (ii)$, on
obtient
$$\int_{\partial X^*}\varphi_j
\eta=2\pi i\sum_{k=1}^nc_k\varphi_j(p_k).$$ En tenant compte des
d\'{e}finitions de $\varphi_j$, $p_k$ et $c_k$, on obtient
\begin{eqnarray}
\int_{\partial X^*}\varphi_j \eta&=&2\pi
i\sum_{k=1}^n c_k\int_{p_0}^p\omega_j,\nonumber\\
&=&2\pi i\sum_{k=1}^m \int_{q_k}^{p_k}\omega_j,\nonumber\\
&=&2\pi i\int_\gamma \omega_j.\nonumber
\end{eqnarray}
En comparant cette derni\`{e}re expression avec celle obtenue dans
(8.2), on obtient
$$\eta(b_j)=\int_{b_j}\eta=2\pi i\int_\gamma\omega_j.$$
Comme $\gamma$ est un chemin ferm\'{e}, il peut s'\'{e}crire sous
la forme
$$\gamma=\sum_{j=1}^g m_ja_j+\sum_{j=1}^g m_{g+j}b_j.$$
D\`{e}s lors,
\begin{eqnarray}
\eta(b_k)&=&\int_{b_k}\eta,\nonumber\\
&=&2\pi i\sum_{j=1}^g\left(m_j\int_{a_j}\omega_k+m_{g+j}\int_{b_j}\omega_k\right),\nonumber\\
&=&2\pi i\left(m_k+\sum_{j=1}^g m_{g+j}\int_{b_j}\omega_k\right),\nonumber\\
&=&2\pi i\left(m_k+\sum_{j=1}^g m_{g+j}\int_{b_k}\omega_j\right),
\mbox{ car } Z \mbox{ est sym\'{e}trique }, \nonumber\\
&=&2\pi
i\left(m_k+\sum_{j=1}^gm_{g+j}\omega_j(b_k)\right).\nonumber
\end{eqnarray}
Posons
$$\theta\equiv \eta-2\pi i\sum_{j=1}^gm_{g+j}\omega_j(b_k).$$
La forme diff\'{e}rentielle $\sum_{j=1}^gm_{g+j}\omega_j(b_k)$
\'{e}tant holomorphe, on en d\'{e}duit que $\theta$ (comme $\eta$)
est une diff\'{e}rentielle m\'{e}romorphe ayant des p\^{o}les
simples en $p_j$ et dont les r\'{e}sidus sont
$\mbox{R\'{e}s}_{p_j}\theta=c_j\in \mathbb{Z}$. En outre, on a
\begin{eqnarray}
\int_{b_k}\theta&=&\theta(b_k),\nonumber\\
&=&\eta(b_k)-2\pi i\sum_{j=1}^gm_{g+j}\omega_j(b_k),\nonumber\\
&=&2\pi i\left(m_k+\sum_{j=1}^gm_{g+j}\omega_j(b_k)\right)-2\pi i\sum_{j=1}^gm_{g+j}\omega_j(b_k),\nonumber\\
&=&2\pi i m_k,\nonumber
\end{eqnarray}
et
\begin{eqnarray}
\int_{a_k}\theta&=&\theta(a_k),\nonumber\\
&=&\eta(a_k)-2\pi i\sum_{j=1}^gm_{g+j}\omega_j(a_k),\nonumber\\
&=&2\pi i m_k,\nonumber
\end{eqnarray}
car $\eta(a_k)=\int_{a_k}\eta=0$ et
$\omega_j(a_k)=\int_{a_k}\omega_j=\delta_{kj}=0$. Donc
l'int\'{e}grale de $\theta$ le long de tout chemin ferm\'{e} est
d\'{e}finie \`{a} un multiple entier de $2\pi i$ pr\`{e}s. La
fonction que l'on cherche à d\'{e}terminer est
$$f(p)=e^{\int_{p_0}^p\theta}.$$
En effet, cette fonction est bien d\'{e}finie et nous allons voir
qu'elle est m\'{e}romorphe et que
$$(f) =\sum_{j=1}^{m}q_{j}-\sum_{j=1}^{m}p_{j}.$$
En effet, au voisinage de $p_j$, on a
$$\theta=\left(\frac{c_j}{t}+g(t)\right)dt,$$
o\`{u} g(t) est une fonction holomorphe et
\begin{eqnarray}
f(\varepsilon)&=&e^{\int_{t}^\varepsilon \left(\frac{c_j}{t}+g(t)\right)dt},\nonumber\\
&=&e^{c_j\log\varepsilon-c_j\log t+\int_{t}^\varepsilon g(t)dt},\nonumber\\
&=&\varepsilon^{c_j}G(t),\nonumber
\end{eqnarray}
o\`{u} $G(t)$ est une fonction holomorphe. On en d\'{e}duit que
suivant le signe de $c_j$, $p_j$ est un z\'{e}ro ou un p\^{o}le de
$f$ d'ordre $|c_j|$. Finalement, on a
\begin{eqnarray}
(f)&=&\sum_{j=1}^{n}c_jp_j,\nonumber\\
&=&\sum_{j=1}^{m}q_{j}-\sum_{j=1}^{m}p_{j},\nonumber
\end{eqnarray}
et le th\'{e}or\`{e}me est d\'{e}montr\'{e}. $\square$

On d\'{e}signe par $L_\Omega$ ou tout simplement $L$ le r\'{e}seau
dans $\mathbb{C}^g$ d\'{e}fini par le $\mathbb{Z}$-module
$$L=\mathbb{Z}^g\oplus \Omega \mathbb{Z}^g,$$
ou encore
$$L=\left\{\sum_{j=1}^g\left(k_j\int_{a_j}\left(\begin{array}{c}
\omega_1\\
\vdots\\
\omega_g
\end{array}\right)+m_j\int_{b_j}\left(\begin{array}{c}
\omega_1\\
\vdots\\
\omega_g
\end{array}\right)\right) : k_j, m_j\in
\mathbb{Z}\right\},$$ i.e., le sous-groupe de $\mathbb{C}^g$
engendr\'{e} par les vecteurs colonnes de la matrice des
p\'{e}riodes de $\Omega$.

\begin{defn}
L'espace quotient $\mathbb{C}^g/L$, s'appelle vari\'{e}t\'{e}
jacobienne de $X$ et on le d\'{e}signe par $\mbox{Jac} (X)$.
\end{defn}

La vari\'{e}t\'{e} jacobienne de $X$, est un tore complexe de
dimension $g$. En effet, en utilisant la suite exponentielle
exacte de faisceaux
$$0\longrightarrow \mathbb{Z}\overset{i}{\longrightarrow}\mathcal{O}_{X}
\overset{\exp }{\longrightarrow
}\mathcal{O}_{X}^{*}\longrightarrow 0,$$ o\`{u} $\mathcal{O}_{X}$
est le faisceau des fonctions holomorphes sur $\mathcal{C}$,
$\mathcal{O}_{X}^{*}$ est le faisceau des fonctions holomorphes ne
s'annulant pas sur $X$, $i$ est l'inclusion triviale et $\exp $
est l'application exponentielle $\exp f=e^{2\pi \sqrt{-1}f}$ ainsi
que la dualit\'{e} (voir par exemple [9] ou [18]), on montre que
\begin{eqnarray}
Jac(X)&=&H^{1}\left( X,\mathcal{O}_{X}\right) /H^{1}\left(
X,\mathbb{Z}\right) ,\nonumber\\
&\simeq &H^{1}\left(X,\Omega _{X}^{1}\right)
/H^{1}\left( X,\mathbb{Z}\right) ,\nonumber\\
&\simeq &H^{0}(X,\Omega _{X}^{1}) ^{*}/H_{1}(X,\mathbb{Z}),
\nonumber\\
&\simeq &\mathbb{C}^{g}/\mathbb{Z}^{2g}.\nonumber
\end{eqnarray}

\begin{rem}
Soit $\mathcal{D}=\sum_{j=1}^mn_jq_j\in \mbox{Div }(X)$, $p\in X$,
fix\'{e} et soit $(\omega_1,...,\omega_g)$ une base de
diff\'{e}rentielles holomorphes sur $X$. L'application
$$\varphi : \mbox{Div }(X)\longrightarrow \mbox{Jac}(X),\quad
\mathcal{D}\longmapsto
\left(\sum_{j=1}^{m}n_j\int_{p}^{q_{j}}\omega_1,...,\sum_{j=1}^{m}n_j\int_{p}^{q_{j}}\omega_g\right),$$
est dite "application d'Abel-Jacobi". Dans le cas particulier
o\`{u}
$$\mathcal{D}=\mathcal{D}_1-\mathcal{D}_2=\sum_{j=1}^mq_j-\sum_{j=1}^mp_j,$$
la condition (i) signifie que $\mathcal{D}\in \mbox{Div}^0(X)$ ou
encore $\mathcal{D}_1$ est \'{e}quivalente à $\mathcal{D}_2$. La
condition (ii) peut s'\'{e}crire sous une forme condens\'{e}e,
$$\forall\omega \in \Omega^1(X) ,\quad\int_{\mathcal{D}_1}^{\mathcal{D}_2}\omega =\int_{\gamma}\omega.$$
Notons que la condition (ii) peut encore s'\'{e}crire sous la
forme
$$\varphi(\mathcal{D})\equiv
\left(\sum_{j=1}^{m}\int_{p_{j}}^{q_{j}}\omega_1,...,\sum_{j=1}^{m}\int_{p_{j}}^{q_{j}}\omega_g\right)
\equiv 0 \mbox{ mod. } L,$$ avec $\varphi$ l'application
d\'{e}finie par $\varphi : \mbox{Div}^\circ (X)\longrightarrow
\mbox{Jac}(X)$.
\end{rem}

\section{Le probl\`{e}me d'inversion de Jacobi}

Avant tout nous avons besoin d'un r\'{e}sultat qui d\'{e}coule du
th\'{e}or\`{e}me de Riemann-Roch.

\begin{prop}
Soit $\mathcal{D}$ un diviseur positif sur une surface de Riemann
compacte $X$ de genre $g>0$. Alors\\
a) $\dim \mathcal{I}(-\mathcal{D})\geq g-\mbox{deg }\mathcal{D}$.\\
b) Pour tout $p\in X$, on a $\dim \mathcal{I}(-p)\leq g-1$.
Autrement dit, il existe une diff\'{e}rentielle $\omega$
holomorphe sur $X$ telle que : $\omega(p)\neq 0$.
\end{prop}
\emph{D\'{e}monstration}: $a)$ Rappelons que les seules fonctions
holomorphes sur une surface de Riemann compacte $X$ sont les
constantes. Donc
$$\mathcal{L}(0)=\{\mbox{fonctions constantes sur } X\}\simeq
\mathbb{C},$$ et $\dim \mathcal{L}(\mathcal{D})\geq 1$.
D'apr\`{e}s le th\'{e}or\`{e}me de Riemann-Roch, on a
$$\dim \mathcal{I}(-\mathcal{D})+\mbox{deg }\mathcal{D}-g+1\geq 1,$$
ce qui implique que
$$\dim \mathcal{I}(-\mathcal{D})\geq g-\mbox{deg }\mathcal{D}.$$
$b)$ Proc\'{e}dons par l'absurde, i.e., supposons que : $\forall
\omega \in \Omega^1(X)$, $\omega(p)=0$. Donc $\dim
\mathcal{I}(-p)=g$ et en vertu du th\'{e}or\`{e}me de
Riemann-Roch, on a
$$\dim \mathcal{L}(\mathcal{D})=\dim \mathcal{I}(-p)+\mbox{deg
}p-g+1=2.$$ D'o\`{u} $\mathcal{L}(p)=\{1 ,f\}$ o\`{u} $f$ est une
fonction m\'{e}romorphe non constante (ayant au plus un p\^{o}le
simple en $p$) telle que : $(f)\geq -p$. La fonction $f :
X\longrightarrow \mathbb{P}^1(\mathbb{C})$ n'a pas de points de
branchements et la surface $X$ peut \^{e}tre vue comme \'{e}tant
un rev\^{e}tement non ramifi\'{e} de $\mathbb{P}^1(\mathbb{C})$,
$\mbox{deg }f=1$, ce qui implique que $X\simeq
\mathbb{P}^1(\mathbb{C})$. Ceci est absurde puisque $g(X)> 0$ par
hypoth\`{e}se et $g(\mathbb{P}^1(\mathbb{C}))=0$, ce qui
ach\`{e}ve la d\'{e}monstration. $\square$

Soit $X$ une surface de Riemann compacte de genre $g>0$. On note
$\mbox{Sym}^dX$ l'ensemble de tous les diviseurs positifs
$$\mathcal{D}=\sum_{j=1}^d p_j,$$
de degr\'{e} $d$ sur $X$. On dit que $\mbox{Sym}^dX$ est le
$d^{\mbox{ième}}$ produit sym\'{e}trique de $X$. On montre
ais\'{e}ment que $\mbox{Sym}^dX$ peut-\^{e}tre muni d'une
structure de vari\'{e}t\'{e} complexe. Soit
$$X^d=\underset{k-fois}{\underbrace{X\times ...\times X}},$$
le produit direct de $X$; c'est une vari\'{e}t\'{e} complexe.
Consid\'{e}rons le groupe sym\'{e}trique $\Sigma_d$ des
permutations de $\{1,...,d\}$. D\`{e}s lors, pour tout $\sigma\in
\Sigma_d$, on d\'{e}finit l'application $\sigma :
X^d\longrightarrow X^d$, en posant
$$\sigma(p_1,...,p_d)=(p_{\sigma_1},...,p_{\sigma_d}),$$
o\`{u} $\sigma_1,...,\sigma_d$ d\'{e}signent les fonctions
sym\'{e}triques \'{e}l\'{e}mentaires. Cette application est
biholomorphe, d'o\`{u} $\Sigma_d\subset \mbox{Aut }(X^d)$, i.e.,
$\Sigma_d$ peut-\^{e}tre vu comme \'{e}tant un sous groupe du
groupe des automorphismes de $X^d$. Notons que $\mbox{Sym}^dX$
h\'{e}rite de $X^d$ d'une structure d'espace topologique. D\`{e}s
lors, l'espace quotient $\mbox{Sym}^dX=X^d/\Sigma_d$ est un espace
s\'{e}par\'{e} compact. La projection $\pi : X^d\longrightarrow
\mbox{Sym}^dX$ munit $\mbox{Sym}^dX$ d'une structure de
vari\'{e}t\'{e} complexe. En effet, soit $p_j \in X$,
$\mathcal{D}=\sum p_j\in \mbox{Sym}^dX$, $p_j\neq p_k$. Autour de
chaque $p_j$, on choisit un syst\`{e}me de coordonn\'{e}es locales
$(U_j, z_j)$ dans $X$. On suppose que pour $p_j\neq p_k$, on a
$U_j\cap U_k\neq 0$ et que pour $p_j=p_k$, on a $z_j=z_k$ dans
$U_j=U_k$. L'application
$$\sum q_j \longmapsto \left(\sigma_1(z_j(q_j)),...,\sigma_d(z_j(q_j))\right),$$
d\'{e}termine, d'apr\`{e}s le th\'{e}or\`{e}me fondamental
d'alg\`{e}bre, une carte locale sur $\pi(U_1\times...\times
U_d)\subset \mbox{Sym}^dX$. En dehors des points de branchements
de la surface, l'application $\pi$ est un rev\^{e}temnt et on peut
prendre $\left(z_1(p_1)),...,z_d(p_d)\right)$ comme
coordonn\'{e}es autour de $\mathcal{D}\in \mbox{Sym}^dX$. Autour
d'un point $d.p$, l'ensemble $$(z_1+...+z_d,...,z_1...z_d),$$
forme un syst\`{e}me de coordonn\'{e}es locales.

\begin{thm}
Soit
$$\varphi_g : \mbox{Sym}^gX\longrightarrow \mbox{Jac}(X),\quad
\mathcal{D}\longmapsto
\varphi_g(\mathcal{D})=\left(\int_0^\mathcal{D}\omega_1,...,\int_0^\mathcal{D}\omega_g\right),$$
l'application d'Abel-Jacobi restreinte \`{a} l'espace
$\mbox{Sym}^gX$ o\`{u} $(\omega_1,...,\omega_g)$ est une base
normalis\'{e}e de
$\Omega^1(X)$. Alors\\
a) L'application $\varphi_g$ est bien d\'{e}finie.\\
b) L'application $\varphi_g$ est injective.\\
c) L'application $\varphi_g$ est surjective. Si $\mathcal{D}_{1}$
est un diviseur positif de degr\'{e} $g$, alors, pour tout $\left(
s_{1},\ldots ,s_{g}\right) \in \mathbb{C}^{g}$, il existe un
diviseur $\mathcal{D}_{2}$\ positif de degr\'{e} $g$\ tel que:
$$\int_{\mathcal{D}_{1}}^{\mathcal{D}_{2}}\omega =s_{k}.$$
L'application $\varphi : \mbox{Div }^0(X)\longrightarrow
\mbox{Jac}(X)$, est surjective. (Probl\`{e}me d'inversion de
Jacobi).
\end{thm}
\emph{D\'{e}monstration}: $a)$ Montrons que l'application
$\varphi_g$ est bien d\'{e}finie. Autrement dit, montrons que deux
\'{e}l\'{e}ments \'{e}quivalents dans $\mbox{Sym}^gX$, sont
envoy\'{e}s sur deux \'{e}l\'{e}ments \'{e}quivalents dans
$\mathbb{C}^{g}/L$. Soient donc $\mathcal{D}_1$ et $\mathcal{D}_2$
deux diviseurs \'{e}quivalents dans $\mbox{Sym}^gX$ et $\gamma $
un chemin ferm\'{e} sur $X$. D'apr\`{e}s le th\'{e}or\`{e}me
d'Abel et la remarque 8.1, on a
\begin{eqnarray}
\left(\int_0^{\mathcal{D}_2}\omega_1,...,\int_0^{\mathcal{D}_2}\omega_g\right)
-\left(\int_0^{\mathcal{D}_1}\omega_1,...,\int_0^{\mathcal{D}_1}\omega_g\right)
&=&\left(\int_{\mathcal{D}_1}^{\mathcal{D}_2}\omega_1,...,\int_{\mathcal{D}_1}^{\mathcal{D}_2}\omega_g\right),\nonumber\\
&=&\left(\int_\gamma\omega_1,...,\int_\gamma\omega_g\right),\nonumber
\end{eqnarray}
et il suffit de montrer que $\int_\gamma\omega_j\in L$, $1\leq
j\leq g$. Le chemin $\gamma$ \'{e}tant ferm\'{e}, on peut donc
l'\'{e}crire sous la forme suivante :
$$\gamma=\sum_{k=1}^g\left(\alpha_ka_k+\beta_kb_k\right),\quad \left(\alpha_k, \beta_k \in \mathbb{Z}\right),$$
o\`{u} $\left( a_{1},\ldots ,a_{g},b_{1},\ldots ,b_{g}\right) $
est une base de cycles dans le groupe d'homologie $H_{1}\left(
X,\mathbb{Z}\right)$. D\`{e}s lors
$$\int_\gamma\omega_j=\sum_{k=1}^g\left(\alpha_k\int_{a_k}\omega_j+\beta_k\int_{b_k}\omega_j\right),\quad 1\leq
j\leq g.$$ Nous avons montr\'{e} pr\'{e}c\'{e}demment que la
matrice $\Omega$ des p\'{e}riodes de $X$ peut s'\'{e}crire sous la
forme
\begin{eqnarray}
\Omega&=&(E, F),\quad\mbox{(d\'{e}finition 4.2)},\nonumber\\
&=&(I, Z),\quad Z=E^{-1}F,\quad\mbox{(proposition 4.6)},\nonumber\\
&=&\left(\begin{array}{cccccccc}
1&0&\cdots &0&\int_{b_1}\omega_1&\int_{b_1}\omega_2&\cdots &\int_{b_1}\omega_g\\
0&1&\cdots &0&\int_{b_2}\omega_1&\int_{b_2}\omega_2&\cdots &\int_{b_2}\omega_g\\
\vdots&&\ddots &&\vdots&\vdots&\ddots&\vdots\\
0&0& &1&\int_{b_g}\omega_1&\int_{b_g}\omega_2&\cdots
&\int_{b_g}\omega_g
\end{array}\right).\nonumber
\end{eqnarray}
D'o\`{u}
$$\int_\gamma\omega_j=\alpha_j+\sum_{k=1}^g\beta_k \int_{b_k}\omega_j,\quad 1\leq
j\leq g,$$ ce qui montre \footnote{Rappelons que si un r\'{e}seau
$L$ de $g$ points dans $\mathbb{R}^g$ est d\'{e}fini par les
vecteurs $c_1,...,c_g$ rapport\'{e}s \`{a} l'origine, alors dire
qu'un point $(x_1,...,x_g)\in L$, cel\`{a} signifie que
$(x_1,...,x_g)=\beta_1c_1+...+\beta_gc_g$, $(\beta_1,...,\beta_g
\in \mathbb{Z})$ ou encore en terme de composantes
$(c_{j1},...,c_{jg})$ du vecteur $c_j$, $1\leq j\leq g$, il faut
que : $x_j=\sum_{k=1}^g\beta_kc_{kj}$, $1\leq j\leq g$. Ici, nous
avons un r\'{e}seau de $2g$ points dans $\mathbb{C}^g$ dont les
$g$ premiers sont les vecteurs unit\'{e}s. Donc dire que
$\int_\gamma\omega_j\in L$, cel\`{a} est \'{e}quivaut du point de
vue des composantes \`{a}
$\int_\gamma\omega_j=\alpha_j+\sum_{k=1}^g\beta_k c_{kj}$, $1\leq
j\leq g$ et il suffit de choisir
$c_{kj}=\int_{b_k}\omega_j$.} que $\int_\gamma\omega_j\in L$.\\
$b)$ Montrons que l'application $\varphi_g$ est injective.
Autrement dit, montrons que  si deux diviseurs $\mathcal{D}_1$ et
$\mathcal{D}_2$ sont envoy\'{e}s sur des points \'{e}quivalents,
alors $\mathcal{D}_1\sim\mathcal{D}_2$. Soit
$\left(\int_0^{\mathcal{D}_1}\omega_1,...,\int_0^{\mathcal{D}_1}\omega_g\right)$
l'image de $\mathcal{D}_1$ et
$\left(\int_0^{\mathcal{D}_2}\omega_1,...,\int_0^{\mathcal{D}_2}\omega_g\right)$
celui de $\mathcal{D}_2$. Ces images \'{e}tant \'{e}quivalentes,
alors
$$\left(\int_{\mathcal{D}_1}^{\mathcal{D}_2}\omega_1,...,\int_{\mathcal{D}_1}^{\mathcal{D}_2}\omega_g\right)\in
L,$$ et d\`{e}s lors
\begin{eqnarray}
\left(\int_{\mathcal{D}_1}^{\mathcal{D}_2}\omega_j\right)&=&\alpha_j+\sum_{k=1}^g\beta_k
c_{kj},\nonumber\\
&=&\sum_{k=1}^g\alpha_k\delta_{kj}+\sum_{k=1}^g\beta_k c_{kj},\nonumber\\
&=&\sum_{k=1}^g\alpha_k\int_{a_k}\omega_j+\sum_{k=1}^g\beta_k \int_{b_k}\omega_j,\nonumber\\
&=&\int_\gamma\omega_j,\nonumber
\end{eqnarray}
o\`{u}
$$\gamma=\sum_{k=1}^g\left(\alpha_ka_k+\beta_kb_k\right),\quad \left(\alpha_k, \beta_k \in \mathbb{Z}\right),$$
est un chemin ferm\'{e} et ne d\'{e}pend pas de $j$. Par
cons\'{e}quent, pour tout $\omega$, on a
$$\left(\int_{\mathcal{D}_1}^{\mathcal{D}_2}\omega_j\right)=\int_\gamma\omega,$$
et d'apr\`{e}s le th\'{e}or\`{e}me d'Abel, $\mathcal{D}_1\sim\mathcal{D}_2$.\\
$c)$ La preuve va se faire en plusieurs \'{e}tapes :\\
\underline{\'{E}tape 1} : Soit $(\omega_1,...,\omega_g)$ une base
normalis\'{e}e de $\Omega^1(X)$ et choisissons un diviseur positif
$\mathcal{D}$ sur $X$ de degr\'{e} $g$ tel que :
$(\omega_j(p_j))\neq 0$, o\`{u} $p_j\in X$ et $1\leq j\leq g$.
Nous verrons ci-dessous que ce choix est toujours possible et on
dira que "$\mathcal{D}$ est g\'{e}n\'{e}ral". On veut montrer
qu'il existe un diviseur positif $\mathcal{D}_2$ de degr\'{e} $g$
tel que :
$$\int_{\mathcal{D}_1}^{\mathcal{D}_2}\omega_k=s_k.$$
Ceci est \'{e}quivalent \`{a} montrer l'existence du diviseur
$\mathcal{D}$ ci-dessus tel que :
$$\int_{\mathcal{D}}^{\mathcal{D}_2}\omega_k=t_k,\quad \forall (t_1,...,t_g)\in \mathbb{C}^g.$$
En effet, on a
$$s_k=\int_{\mathcal{D}_1}^{\mathcal{D}_2}\omega_k
=\int_{\mathcal{D}_1}^{\mathcal{D}}\omega_k+\int_{\mathcal{D}}^{\mathcal{D}_2}\omega_k,$$
d'o\`{u}
$$\int_{\mathcal{D}}^{\mathcal{D}_2}\omega_k=s_k-\int_{\mathcal{D}_1}^{\mathcal{D}}\omega_k=t_k.$$
\underline{\'{E}tape 2} : Montrons que les conditions suivantes
sont \'{e}quivalentes,
\begin{eqnarray}
&(i)&\mathcal{D} \mbox{ est g\'{e}n\'{e}ral}.\nonumber\\
&(ii)&\dim\mathcal{L}(\mathcal{D})=1.\nonumber\\
&(iii)&\dim \mathcal{I}(-\mathcal{D})=0.\nonumber
\end{eqnarray}
On a, $(ii)\Longleftrightarrow(iii)$. En effet, d'apr\`{e}s le
th\'{e}or\`{e}me de Riemann-Roch, on a
$$\dim\mathcal{L}(\mathcal{D})=\dim
\mathcal{I}(-\mathcal{D})+\mbox{deg }\mathcal{D}-g+1.$$ Or
$\mbox{deg }\mathcal{D}=g$, donc $\dim\mathcal{L}(\mathcal{D})=1$
si et seulement si $\dim \mathcal{I}(-\mathcal{D})=0$. Montrons
maintenant que $(i)\Longleftrightarrow(iii)$ ou ce qui revient au
m\^{e}me $\mbox{non }(i)\Leftrightarrow \mbox{non }(iii)$. En
effet, $\mbox{non }(iii)$ signifie que $\dim
\mathcal{I}(-\mathcal{D})\neq 0$, i.e., il existe une forme
diff\'{e}rentielle $\omega$ holomorphe telle que : $(\omega)\geq
\mathcal{D}$. Autrement dit, pour tous $p_1,...,p_g\in
\mathcal{D}$, on peut trouver des coefficients $c_1,...,c_g$ tels
que :
$$\omega=\sum_{k=1}^gc_k\omega_k(p_j),\quad 1\leq j\leq g,$$
o\`{u} $(\omega_1,...,\omega_g)$ est une base de $\Omega^1(X)$.
Les coefficients $c_1,...,c_g$ existent si et seulement si ce
syst\`{e}me homog\`{e}ne de $g$ \'{e}quations \`{a} $g$ inconnues
poss\`{e}de une solution non triviale. Autrement dit si et
seulement si
$$
\det \left(\begin{array}{ccc}
\omega _{1}(p_1)&\cdots &\omega _{1}(p_1)\\
\vdots &\ddots & \vdots \\
\omega _{g}(p_g)&\cdots &\omega _{g}(p_g)
\end{array}\right)=0,
$$
ou encore si et seulement si la condition $(i)$ n'est pas
satisfaite.\\
\underline{\'{E}tape 3} : Pour montrer que le diviseur
$\mathcal{D}$ est g\'{e}n\'{e}ral, il suffit donc de prouver que
l'une des conditions $(ii)$ ou $(iii)$ mentionn\'{e}e dans
l'\'{e}tape 2, est satisfaite. D'apr\`{e}s la proposition 9.1, on
a
$$\dim \mathcal{I} (-p_1)=g-1,$$
ce qui montre qu'il existe $p_2\in X$ tel que : $\mathcal{I}
(-p_1-p_2)\subset\mathcal{I} (-p_1)$ et
$$\dim \mathcal{I} (-p_1-p_2)=g-2.$$
De m\^{e}me, il existe $p_3\in X$ tel que : $\mathcal{I}
(-p_1-p_2-p_3)\subset\mathcal{I} (-p_1-p_2)$ et
$$\dim \mathcal{I} (-p_1-p_2-p_3)=g-3.$$
Et ainsi de suite, on peut trouver $p_g\in X$ tel que :
$\mathcal{I} (-p_1-p_2-...-p_g)\subset\mathcal{I}
(-p_1-p_2-...-p_{g-1})$ et
$$\dim \mathcal{I} (-p_1-p_2-...-p_g)=0,$$
i.e., $\dim \mathcal{I} (-\mathcal{D})=0$ et nous avons montr\'{e}
dans l'\'{e}tape 2 ci-dessus que ceci est \'{e}quivalent \`{a}
$\dim \mathcal{L}
(\mathcal{D})=1$ et aussi à $\mathcal{D}$ est g\'{e}n\'{e}ral.\\
\underline{\'{E}tape 4} : Il reste à prouver qu'il existe
$\mathcal{D}_2$ tel que :
$$\int_{\mathcal{D}}^{\mathcal{D}_2}\omega_k=t_k,\quad 1\leq k\leq
g.$$ Posons
\begin{eqnarray}
\mathcal{D}&=&\sum_{j=1}^gp_{0j},\nonumber\\
\mathcal{D}_2&=&\sum_{j=1}^gp_{j},\nonumber
\end{eqnarray}
et consid\'{e}rons la fonction
\begin{eqnarray}
f(p)&\equiv&\left(f_1(p),...,f_g(p)\right),\nonumber\\
&=&\left(\sum_{j=1}^g\int_{p_{0j}}^{p_j}\omega_1,...,\sum_{j=1}^g\int_{p_{0j}}^{p_j}\omega_g\right),\nonumber\\
&=&\left(\frac{t_1}{n},...,\frac{t_g}{n}\right),
\end{eqnarray}
o\`{u} $p=(p_1,...,p_g)$. Notons que nous avons remplac\'{e}
\footnote{pour \^{e}tre sur de travailler dans un voisinage assez
petit} $t_k$ par $\frac{t_k}{n}$ o\`{u} $n$ est un entier
suffisamment grand. D'apr\`{e}s le th\'{e}or\`{e}me des fonctions
implicites, on peut d\'{e}terminer $p$ explicitement car la
matrice jacobienne
$$
\left(\frac{\partial f_j}{\partial p_k}\right)_{1\leq j ,k g}=
\left(\begin{array}{ccc}
\omega _{1}(p_1)&\cdots &\omega _{1}(p_1)\\
\vdots &\ddots & \vdots \\
\omega _{g}(p_g)&\cdots &\omega _{g}(p_g)
\end{array}\right),$$
est inversible d'apr\`{e}s la d\'{e}finition du diviseur
$\mathcal{D}$. D'apr\`{e}s (9.1), on a
$$f_k(p)=\frac{t_k}{n},\quad 1\leq k\leq g$$
et
\begin{eqnarray}
f_k(p)&=&\sum_{j=1}^g\int_{p_{0j}}^{p_j}\omega_k,\quad 1\leq k\leq g,\nonumber\\
&=&\int_{\mathcal{D}}^{\mathcal{D}_2}\omega_k,\nonumber
\end{eqnarray}
d'o\`{u}
$$n\int_{\mathcal{D}}^{\mathcal{D}_2}\omega_k=t_k\equiv
\int_{\mathcal{D}}^{\mathcal{D}_3}\omega_k.$$ On doit donc trouver
un diviseur $\mathcal{D}_3$ tel que :
$$n\int_{\mathcal{D}}^{\mathcal{D}_2}\omega_k=
\int_{\mathcal{D}}^{\mathcal{D}_3}\omega_k.$$ Celui-ci existe
d'apr\`{e}s le th\'{e}or\`{e}me
d'addition\footnote{Th\'{e}or\`{e}me d'addition : Soient
$\mathcal{D}_{1}$\ et $\mathcal{D}_{2}$\ deux diviseurs positifs
de degr\'{e} $n,$\ $\mathcal{E}_{1}$\ un diviseur positif de
degr\'{e} $g$\ et $\omega \in \Omega^1(X)$. Alors, il existe un
diviseur $\mathcal{E}_{2}$\ positif de degr\'{e} $g$ tel que:
$\int_{\mathcal{D}_{1}}^{\mathcal{D}_{2}}\omega =
\int_{\mathcal{E}_{1}}^{\mathcal{E}_{2}}\omega.$ }.
Consid\'{e}rons enfin les diviseurs de degr\'{e} $0$ de la forme
$\mathcal{D}-pg$ o\`{u} $\mathcal{D}\in \mbox{Sym}^gX$. Ces
diviseurs forment un ensemble que l'on note
$$\mbox{Sym}^gY\equiv \mbox{Sym}^gX-pg.$$
L'application $\varphi_g$ \'{e}tant surjective sur l'espace
$\mbox{Sym}^gX$, elle est donc aussi surjective sur
$\mbox{Sym}^gY$. Par cons\'{e}quent, $\varphi$ est surjective sur
$\mbox{Div}^0(X)$ et la d\'{e}monstration s'ach\`{e}ve. $\square$

\section{Appendices}

\subsection{Courbes elliptiques et hyperelliptiques}

Nous allons dans cet appendice construire le plus intuitivement
possible la surface de Riemann dans le cas elliptique et
hyperelliptique.

Soit
$$w^{2}-P_{3}\left( z\right) =0,$$
o\`{u} $P_{3}\left( z\right) $ est un polyn\^{o}me de degr\'{e}
$3,$ ayant trois racines distinctes $e_{1},e_{2},e_{3}.$
Consid\'{e}rons
$$\mathbb{C}\longrightarrow \mathbb{C},\text{ }z\longmapsto w:w^{2}=P_{3}\left( z\right)
,$$ Il est \'{e}vident que $w$ n'est pas une fonction (uniforme).
A chaque valeur de $z$ correspond deux valeurs diff\'{e}rentes de
$w$ sauf quand $z=e_{1},$ $z=e_{2}$ et $z=e_{3}.$ En ces points,
$w$ est univalu\'{e}e: en effet, on a $w=\pm \sqrt{P_{3}\left(
e_{i}\right) }=0,$ une seule valeur. Tous les points \`{a}
l'infini dans toutes les directions seront identifi\'{e}s en un
seul point que l'on d\'{e}signe par $\infty .$ Au point $z=\infty
,$ $w$ est aussi univalu\'{e}e: en effet, posons $z=\frac{1}{t},$
d'o\`{u}
\begin{eqnarray}
w^{2}&=&P_{3}\left( \frac{1}{t}\right) ,\nonumber\\
&=&\left( \frac{1}{t}-e_{1}\right) \left( \frac{1}{t}-e_{2}\right)
\left( \frac{1}{t}-e_{3}\right) ,\nonumber\\
&=&\frac{1}{t^{3}}\left( 1-\left( e_{1}+e_{2}+e_{3}\right)
t+\left( e_{1}e_{2}+e_{1}e_{3}+e_{2}e_{3}\right)
t^{2}-e_{1}e_{2}e_{3}t^{3}\right) ,\nonumber
\end{eqnarray}
et
$$w=\pm \sqrt{P_{3}\left( \frac{1}{t}\right) }\sim \pm
\sqrt{\frac{1}{t^{3}}}.$$ Par cons\'{e}quent,
$\underset{t\rightarrow 0}{\lim }w=\pm \infty $ c'est-à-dire
$\infty$, une seule valeur.

Notre probl\`{e}me consiste \`{a} uniformiser $w$, autrement dit,
on cherche un domaine sur lequel $w$ est une fonction (uniforme).
Auparavant, \'{e}tudions le comportement de $w$ au voisinage des
racines de $P_{3}\left( z\right) =0$ c'est-à-dire
$e_{1},e_{2},e_{3}$ ainsi qu'au voisinage du point \`{a} l'infini
$\infty$. Si $z$ d\'{e}crit un circuit (c'est-à-dire un chemin
ferm\'{e}, par exemple un cercle) entourant un des points
$e_{1},e_{2},e_{3}$ et $\infty$, alors $w$ change de signe : en
effet, supposons que $z$ d\'{e}crit un cercle centr\'{e} en
$e_{1}$ et posons $z-e_{1}=re^{i\theta }$ o\`{u} $r$ est le module
de $z-e_{1}$ et $\theta$ son argument. Evidemment $r$ ne change
pas tandis que $\theta$ varie de $0$ \`{a} $2\pi$. Au voisinage de
$e_{1},$ $w=\sqrt{P_{3}\left( z\right) }$ se comporte comme
$$w=\sqrt{z-e_{1}}=r^{1/2}e^{i\theta /2}.$$ D\`{e}s lors, pour
$\theta =0,$ on a $w=r^{1/2}$ tandis que pour $\theta =2\pi ,$ on
a $w=-r^{1/2}$. Si on refait de nouveau un tour complet autour de
$z=e_{1}$, l'argument $\theta$ varie de $2\pi$ à $4\pi$ et alors
on obtient $r^{1/2}$ qui est la valeur de d\'{e}part. Pour
$z=e_{2}$ ou $z=e_{3}$, il suffit d'utiliser un raisonnement
similaire au cas pr\'{e}c\'{e}dent. En ce qui concerne le point
$\infty$, on pose comme pr\'{e}c\'{e}demment $z=\frac{1}{t}$ et on
\'{e}tudie $w^{2}=P_{3}\left( \frac{1}{t}\right)$ au voisinage de
$t=0$. On a
$$
w^{2}=\frac{1}{t^{3}}\left( 1-\left( e_{1}+e_{2}+e_{3}\right) t+
\left( e_{1}e_{2}+e_{1}e_{3}+e_{2}e_{3}\right)
t^{2}-e_{1}e_{2}e_{3}t^{3}\right),
$$
et
$$w\sim \pm \sqrt{\frac{1}{t^{3}}}=\pm t^{-3/2}.$$
Soit $t=re^{i\theta}$. Autour de $t=0$, $w=\sqrt{P_{3}\left(
\frac{1}{t}\right)}$ se comporte comme
$$w=t^{-3/2}=r^{-3/2}e^{-3i\theta /2}.$$ D\`{e}s lors, pour $\theta
=0$, on a $w=r^{-3/2}$ et pour $\theta =2\pi$, on a $w=-r^{-3/2}
$. Comme pr\'{e}c\'{e}demment, si $t$ refait de nouveau un tour
complet, $w$ reprend la valeur de d\'{e}part c'est-à-dire
$r^{-3/2}$.

Passons maintenant \`{a} la construction du domaine sur lequel $w$
serait une fonction uniforme. Cette construction se fera en
plusieurs \'{e}tapes :

\underline{${1}^{\grave{e}re}$ \'{e}tape} : Prenons deux copies ou
feuillets $\sigma _{1}$ et $\sigma _{2}$ du plan complexe
compactifi\'{e} $\mathbb{C}\cup \left\{ \infty \right\} $ ou ce
qui revient au m\^{e}me de la sph\`{e}re de Riemann puisqu'ils
sont hom\'{e}omorphes. Plaçons le feuillet $\sigma _{1}$ au dessus
de $\sigma _{2}$ et sur chacun de ces feuillets marquons les
points $e_{1},e_{2},e_{3},\infty $. Supposons que les points de
$\sigma _{1}$ seront envoy\'{e}s sur
$$w=\sqrt{\left( z-e_{1}\right) \left( z-e_{2}\right) \left( z-e_{3}\right)
},$$ et que ceux de $\sigma _{2}$ seront envoy\'{e}s sur
$$w=-\sqrt{\left( z-e_{1}\right) \left( z-e_{2}\right) \left( z-e_{3}\right)}.$$

\underline{${2}^{\grave{e}me}$ \'{e}tape} : Dans chaque feuillet,
faisons deux coupures: une le long de la courbe reliant le point
$e_{1}$ au point $e_{2}$ et l'autre le long de la courbe reliant
le point $e_{3}$ au point $\infty$. D\'{e}signons par
$A_{1},B_{1},C_{1},D_{1}$ (resp. $A_{2},B_{2},C_{2},D_{2}$) les
bords des coupures dans le feuillet $\sigma _{1}$ (resp. $\sigma
_{2}$). Rappelons que $w$ change de signe lorsque l'on tourne d'un
tour autour d'un des points $e_{1},e_{2},e_{3},\infty$. Donc en
allant de $A_{1}$ \`{a} $B_{1}$, on change le signe de $w$
c'est-\`{a}-dire on passe sur l'autre feuillet, l\`{a} o\`{u} $w$
a l'autre signe. De m\^{e}me pour les bords $A_{2}$ et
$B_{2},C_{1}$ et $D_{1},C_{2}$ et $D_{2}$. Par cons\'{e}quent, $w$
a la m\^{e}me valeur sur $A_{1}$ et $B_{2}$, sur $B_{1}$ et
$A_{2}$, sur $C_{1}$ et $D_{2}$ et enfin sur $D_{1}$ et $C_{2}$.

\underline{${3}^{\grave{e}me}$ étape} : On identifie les bords
suivants : $A_{1}$ \`{a} $B_{2}$, $B_{1}$ \`{a} $A_{2}$, $C_{1}$
\`{a} $D_{2}$ et $D_{1}$ \`{a} $C_{2}$. Apr\`{e}s recollement, on
obtient un tore \`{a} un trou.

La surface \`{a} deux feuillets obtenue s'appelle surface de
Riemann elliptique ou courbe elliptique associ\'{e}e \`{a}
l'\'{e}quation
$$w^{2}=\left( z-e_{1}\right) \left( z-e_{2}\right) \left( z-e_{3}\right)
.$$ Sur cette surface, $w$ est une fonction uniforme. Lorsqu'on
tourne autour d'un des points $e_{1},e_{2},e_{3}$, ou $\infty$, on
passe d'un feuillet \`{a} l'autre. En ces points les deux
feuillets se joignent et on les appellent points de branchement ou
de ramification de la surface.
\begin{rem}
$a)$ Si
$$w^{2}-P_{4}\left( z\right) =0,$$
o\`{u} $P_{4}\left( z\right) $ est un polyn\^{o}me de degr\'{e}
$4,$ ayant quatre racines distinctes $e_{1},e_{2},e_{3},e_{4},$
alors on obtient aussi une courbe elliptique. Les points de
branchements sont $e_{1},e_{2},e_{3}$ et $e_{4}$. Notons que si
$$w^{2}=\left( z-e_{1}\right) \left( z-e_{2}\right) \left( z-e_{3}\right) \left( z-e_{4}\right)
,$$ alors la transformation
$$\left( w,z\right) \longmapsto \left( \frac{y}{x^{2}},e_{1}+\frac{1}{x}\right)
,$$ ram\`{e}ne cette \'{e}quation \`{a} la forme
$$y^{2}=\left( 1+\left( e_{1}-e_{2}\right) x\right) \left( 1+\left( e_{1}-e_{3}\right) x\right)
\left( 1+\left( e_{1}-e_{4}\right) x\right).$$

$b)$ Signalons enfin que si
$$w^{2}-P_{n}\left( z\right) =0,$$
o\`{u} $P_{n}\left( z\right) $ est un polyn\^{o}me de degr\'{e}
$n$ sup\'{e}rieur o\`{u} \'{e}gal \`{a} $5$, ayant $n$ racines
distinctes, alors on obtient ce qu'on appelle surface de Riemann
hyperelliptique ou courbe hyperelliptique.
\end{rem}

\subsection{R\'{e}sultants et discriminants}

Soient
\begin{eqnarray}
f(x)&=&a_0x^m+a_1x^{m-1}+\cdots+a_m=a_0\prod_{k=1}^m(x-\alpha_k),\quad
a_0\neq0,\nonumber\\
g(x)&=&b_0x^n+b_1x^{n-1}+\cdots+b_n=b_0\prod_{j=1}^n(x-\beta_j),\quad
b_0\neq0,\nonumber
\end{eqnarray}
deux polyn\^{o}mes de degr\'{e} $m$ et $n$ respectivement. Ici
$(\alpha_1,...,\alpha_m)$ et $(\beta_1,...,\beta_n)$ d\'{e}signent
les racines des polyn\^{o}mes $f$ et $g$ respectivement. Le
r\'{e}sultant des polyn\^{o}mes $f$ et $g$, not\'{e}
$\mbox{R\'{e}s}(f,g)$, est le d\'{e}terminant de leur matrice de
Sylvester, i.e., le d\'{e}terminant de la matrice carr\'{e}e
d'ordre $(m+n)$ suivante :
$$
\left(\begin{array}{ccccccccc}
a_0&a_1&\ldots&\ldots&\ldots&a_m&0&\ldots&0\\
0&a_0&a_1&\ldots&\ldots&\ldots&a_m&\ddots&\vdots\\
\vdots&\ddots&\ddots&\ddots&\ddots&\ddots&\ddots&\ddots&0\\
0&\ldots&0&a_0&a_1&\ldots&\ldots&\ldots&a_m\\
b_0&b_1&\ldots&\ldots&b_n&0&\ldots&\ldots&0\\
0&b_0&b_1&\ldots&\ldots&b_n&0&\ldots&0\\
\vdots&\ddots&\ddots&\ddots&\ddots&\ddots&\ddots&\ddots&\vdots\\
\vdots&\ddots&\ddots&\ddots&\ddots&\ddots&\ddots&\ddots&0\\
0&\ldots&\ldots&0&b_0&b_1&\ldots&\ldots&b_n
\end{array}\right)
$$

\begin{prop}
Les polyn\^{o}mes $f$ et $g$ ont un facteur commun non nul si et
seulement si il existe deux polyn\^{o}mes $F$ et $G$ de degr\'{e}
strictement inf\'{e}rieur \`{a} $m$ et $n$ respectivement tels que
:
$$fG=gF.$$
\end{prop}
\emph{D\'{e}monstration}: On a
\begin{eqnarray}
f&=&Af_1^{m_1}f_2^{m_2}...f_r^{m_r},\nonumber\\
g&=&Bg_1^{n_1}g_2^{n_2}...g_s^{n_s},\nonumber
\end{eqnarray}
o\`{u} $A, B$ sont des constantes et
$f_1^{m_1},...,f_r^{m_r},g_1^{n_1},...,g_s^{n_s}$ sont des
polyn\^{o}mes irr\'{e}ductibles. Supposons que $f$ et $g$ ont un
facteur commun non nul, disons $f_1=g_1$. Consid\'{e}rons les
polyn\^{o}mes
\begin{eqnarray}
F&=&\frac{f}{f_1},\nonumber\\
G&=&\frac{g}{g_1}.\nonumber
\end{eqnarray}
D'o\`{u}
\begin{eqnarray}
\mbox{ deg } F&=&\mbox{ deg } f=m,\nonumber\\
\mbox{ deg } G&=&\mbox{ deg } g=n,\nonumber
\end{eqnarray}
et $$fG=\frac{fg}{g_1}=\frac{gf}{f_1}=gF.$$ R\'{e}ciproquement, on
a
$$fG=gF,$$
avec $\mbox{ deg } F<m$ et $\mbox{ deg } G<n$. Supposons que $f$
et $g$ n'ont pas de facteur commun. Dans ce cas, puisque
$$Af_1^{m_1}f_2^{m_2}...f_r^{m_r}.G=Bg_1^{n_1}g_2^{n_2}...g_s^{n_s}.F,$$
alors pour tout $j=1,2,...,r$, $f_j^{m_j}$ doit appara\^{i}tre
comme facteur dans $F$, i.e., $f$ doit diviser $F$ donc $ \mbox{
deg } f\leq \mbox{ deg } F$ ce qui est absurde car par
hypoth\`{e}se $\mbox{ deg } F<m$. $\square$

\begin{prop}
Les polyn\^{o}mes $f$ et $g$ ont un facteur commun non nul si et
seulement si
$$\mbox{R\'{e}s}(f,g)=0.$$
\end{prop}
\emph{D\'{e}monstration}: D'apr\`{e}s la proposition
pr\'{e}c\'{e}dente, les polyn\^{o}mes $f$ et $g$ ont un facteur
commun non nul si et seulement si il existe deux polyn\^{o}mes
\begin{eqnarray}
F(x)&=&A_0x^{m-1}+A_1x^{m-2}+...+A_{m-1},\nonumber\\
G(x)&=&B_0x^{n-1}+B_1x^{n-2}+...+B_{n-1},\nonumber
\end{eqnarray}
tels que :
$$fG=gF,$$
i.e.,
$$(a_0x^m+...+a_m)(B_0x^{n-1}+...+B_{n-1})=
(b_0x^n+...+b_n)(A_0x^{m-1}+...+A_{m-1}).$$ On identifie les
coefficients :
\begin{eqnarray}
x^{m+n-1}&:&\quad a_0B_0=b_0A_0,\nonumber\\
x^{m+n-2}&:&\quad a_0B_1+a_1B_0=b_0A_1+b_1A_0,\nonumber\\
\vdots&&\nonumber\\
x^0&:&\quad a_mB_{n-1}=b_nA_{m-1}.\nonumber
\end{eqnarray}
D'o\`{u}
\begin{eqnarray}
&& a_0B_0-b_0A_0=0,\nonumber\\
&& a_1B_0+a_0B_1-b_1A_0-b_0A_1=0,\nonumber\\
&&\quad\vdots\nonumber\\
&& a_mB_{n-1}-b_nA_{m-1}=0.\nonumber
\end{eqnarray}
On obtient un syst\`{e}me lin\'{e}aire homog\`{e}ne de $(m+n)$
\'{e}quations dont les inconnues sont
$B_0,...,B_{m-1},A_0,,A_{n-1}$. Ce syst\`{e}me admet une solution
non triviale si et seulement si
$$
\Delta\equiv \det \left(\begin{array}{cccccccc}
a_0&0&\ldots&0&-b_0&0&\ldots&0\\
a_1&a_0&\ddots&\vdots&-b_1&-b_0&\ddots&\vdots\\
\vdots&a_1&\ddots&0&\vdots&-b_1&\ddots&0\\
a_m&\vdots&\ddots&a_0&-b_n&\vdots&\ddots&-b_0\\
0&a_m&\vdots&a_1&0&-b_n&\vdots&-b_1\\
\vdots&\ddots&\ddots&\vdots&\vdots&\ddots&\ddots&\vdots\\
0&\ldots&0&a_m&0&\ldots&0&-b_n
\end{array}\right)=0.
$$
En mettant en \'{e}vidence le signe $-$ dans les $m$ derni\`{e}res
colonnes et en tenant compte du fait que le d\'{e}terminant de la
transpos\'{e}e d'une matrice est le m\^{e}me que celui de la
matrice initiale, on obtient
$$\Delta=\pm \mbox{R\'{e}s} (f,g),$$
et le r\'{e}sultat en d\'{e}coule. $\square$

\begin{prop}
Il existe deux polyn\^{o}mes $F$ et $G$ de degr\'{e} strictement
inf\'{e}rieur \`{a} $m$ et $n$ respectivement tels que :
\begin{eqnarray}
fG-gF&=&\mbox{R\'{e}s}(f,g),\nonumber\\
&=&\mbox{polyn\^{o}me en les coefficients de f et g},\nonumber\\
&=&a_0^nb_0^m\prod_{k=1}^m\prod_{j=1}^n(\alpha_k-\beta_j) ,\nonumber\\
&=&a_0^n\prod_{k=1}^mg(\alpha_k) ,\nonumber\\
&=&(-1)^{mn}b_0^m\prod_{j=1}^nf(\beta_j).\nonumber
\end{eqnarray}
\end{prop}
\emph{D\'{e}monstration}: Si $f$ et $g$ ont un facteur commun,
alors d'apr\`{e}s ce qui pr\'{e}c\`{e}de les polyn\^{o}mes $F$ et
$G$ existent et on a
$$fG-gF=0=\mbox{R\'{e}s} (f,g).$$
Si $f$ et $g$ n'ont pas de facteur commun, alors on cherche $F$ et
$G$ tels que :
$$fG-gF=\mbox{R\'{e}s} (f,g).$$
En raisonnant comme dans la proposition pr\'{e}c\'{e}dente, on
obtient un syst\`{e}me non homog\`{e}ne ayant une solution non
nulle. Autrement dit, le d\'{e}terminant $\Delta$ utilis\'{e} dans
la preuve de la proposition pr\'{e}c\'{e}dente est nul ou ce qui
est \'{e}quivalent $f$ et $g$ n'ont pas de facteur commun, ce qui
est vrai par hypoth\`{e}se et ach\`{e}ve la d\'{e}monstration.
$\square$

Soient $\alpha_1,...,\alpha_m$ les $m$ racines du polyn\^{o}me $f$
compt\'{e}es avec multiplicit\'{e}. Le discriminant de $f$,
not\'{e} $Disc (f)$, est
\begin{eqnarray}
Disc (f)&=&a_0^{2m-2}(-1)^{\frac{m(m-1)}{2}}\prod_{i\neq
j}(\alpha_i-\alpha_j),\nonumber\\
&=&a_0^{2n-2}\prod_{1\leq j<i\leq
m}(\alpha_i-\alpha_j)^2.\nonumber
\end{eqnarray}
\begin{prop}
Le r\'{e}sultant de $f$ et de son polyn\^{o}me d\'{e}riv\'{e} $f'$
est
$$\mbox{R\'{e}s} (f,f')=(-1)^{\frac{m(m-1)}{2}}a_0Disc (f).$$
\end{prop}
\emph{D\'{e}monstration}: On a
$$f(x)=a_0\prod_{k=1}^m(x-\alpha_k),$$
et
$$f'(x)=a_0\sum_{k=1}^m\prod_{j\neq k}(x-\alpha_j).$$
En remplaçant dans cette derni\`{e}re \'{e}quation $x$ par
$\alpha_i$, on constate que tous les termes s'annulent sauf le
i-\`{e}me et d\`{e}s lors
$$f'(\alpha_i)=a_0\prod_{j\neq i}(\alpha_i-\alpha_j).$$
Par ailleurs, on sait que
\begin{eqnarray}
\mbox{R\'{e}s} (f,f')&=&a_0^{m-1}\prod_{i=1}^mf'(\alpha_i),\nonumber\\
&=&a_0^{2m-1}\prod_{j\neq i}(\alpha_i-\alpha_j).\nonumber
\end{eqnarray}
Notons que dans le produit ci-dessus, il y a $m(m-1)$ facteurs.
Comme chacun de ces derniers s'\'{e}crit sous la forme
$\alpha_i-\alpha_j$ et sous la forme $\alpha_j-\alpha_i$, alors
leur produit est $(-1)(\alpha_i-\alpha_j)^2$. En tenant compte du
fait qu'il y a $\frac{m(m-1)}{2}$ paires d'indices $i,j$ avec
$1\leq j<i\leq m$, alors
\begin{eqnarray}
\mbox{R\'{e}s}
(f,f')&=&(-1)^{\frac{m(m-1)}{2}}a_0^{2m-1}\prod_{1\leq
j<i\leq m}(\alpha_i-\alpha_j)^2,\nonumber\\
&=&(-1)^{\frac{m(m-1)}{2}}a_0Disc (f),\nonumber
\end{eqnarray}
et la proposition est d\'{e}montr\'{e}e. $\square$

On d\'{e}duit imm\'{e}diatement des propositions
pr\'{e}c\'{e}dentes le r\'{e}sultat suivant :
\begin{prop}
Le discriminant du polyn\^{o}me $f$ est nul si et seulement si les
polyn\^{o}mes $f$ et $f'$ ont un facteur en commun non constant ou
encore si et seulement si le polyn\^{o}me $f$ admet une racine
multiple.
\end{prop}

\end{document}